\documentclass[review, onefignum, onetabnum]{siamart220329}

\usepackage{setspace}

\nolinenumbers

\newtheorem{assumption}{Assumption}
\newtheorem{example}{Example}
\renewcommand{\d}{\ensuremath{\,\mathrm{d}}}

\usepackage{todonotes}

\usepackage{changes}
\definechangesauthor[name=XZ,
color=orange]{x}
\definechangesauthor[name=Bin,
color=yellow]{b}
\definechangesauthor[name=YF,
color=red]{y}

\graphicspath{{./figs/}{./figs/kde_plot/}}

\usepackage{lipsum}
\usepackage{amsfonts}
\usepackage{graphicx}
\usepackage{epstopdf}
\usepackage{algorithmic}
\ifpdf
  \DeclareGraphicsExtensions{.eps,.pdf,.png,.jpg}
\else
  \DeclareGraphicsExtensions{.eps}
\fi

\newsiamremark{remark}{Remark}
\newsiamremark{hypothesis}{Hypothesis}
\crefname{hypothesis}{Hypothesis}{Hypotheses}
\newsiamthm{claim}{Claim}

\headers{L\'{e}vy Score and Score-Based Particle Algorithm}{Y. HUANG, C. LIU and X. ZHOU}

\title{L\'{e}vy Score Function and Score-Based Particle Algorithm for Nonlinear L\'{e}vy--Fokker--Planck Equations\thanks{Submitted to the editors DATE.
\funding{This work was funded by the Hong Kong General Research Funds  (11308121, 11318522,   11308323),  and the NSFC/RGC Joint Research Scheme [RGC Project No. N-CityU102/20 and NSFC Project No.
12061160462.}}}

\author{
Yuanfei Huang\thanks{Asia Pacific Center for Theoretical Physics, Pohang 37673, Korea, and Department of Mathematics, City University of Hong Kong, Kowloon, Hong Kong SAR, 
  (\email{yuanfei.huang@apctp.org}).}
\and Chengyu Liu\thanks{Department of Data Science, City University of Hong Kong, Kowloon, Hong Kong SAR,
(\email{cliu687-c@my.cityu.edu.hk }).}
\and Xiang Zhou\thanks{Department of Mathematics, City University of Hong Kong, Kowloon, Hong Kong SAR, 
  (\email{xizhou@cityu.edu.hk})}}

\usepackage{amsopn}

\ifpdf
\hypersetup{
    pdftitle={L\'{e}vy Score Function and Score-Based Particle Algorithm for Nonlinear L\'{e}vy--Fokker--Planck Equations},
    pdfauthor={Yuanfei Huang, Chengyu Liu and Xiang Zhou}
}
\fi

\begin{document}
    \maketitle

    \begin{abstract}
        The score function for the diffusion process, also known as the gradient
        of the log-density, is a basic concept to characterize the probability flow
        with important applications in the score-based diffusion generative modelling
        and the simulation of It\^{o} stochastic differential equations. However,
        neither the probability flow nor the corresponding score function for the
        diffusion-jump process are known. This paper delivers mathematical
        derivation, numerical algorithm, and error analysis focusing on the
        corresponding score function in non-Gaussian systems with jumps and
        discontinuities represented by the nonlinear L\'{e}vy--Fokker--Planck equations.
        We propose the L\'{e}vy score function for such stochastic equations,
        which features a nonlocal double-integral term, and we develop its training
        algorithm by minimizing the proposed loss function from samples. Based
        on the equivalence of the probability flow with deterministic dynamics, we
        develop a self-consistent score-based transport particle algorithm to sample
        the interactive L\'{e}vy stochastic process at discrete time grid points.
        We provide error bound for the Kullback--Leibler divergence between the
        numerical and true probability density functions by overcoming the
        nonlocal challenges in the L\'{e}vy score. The full error analysis with
        the Monte Carlo error and the time discretization error is furthermore
        established. To show the usefulness and efficiency of our approach,
        numerical examples from applications in biology and finance are tested.
    \end{abstract}



    \section{Introduction}
    \label{sec1}
    

Stochastic processes, such as the  It\^{o} diffusion process 
and the L\'{e}vy jump process,  play prominent roles in a variety range of fields including  science, engineering, economics, and statistics. 
The temporal evolution of the probability distribution, satisfying the initial-value forward Kolmogorov equation, can be also viewed as a path or curve in the space of probability distributions. 
The path connecting  two probability density functions associated certain stochastic processes  is the crucial mathematical tool   the basis of   score-based generative algorithms \cite{song2020score}.  

The viewpoint of probability flow is to treat the time evolution of the probability density function
as the continuity equation whose velocity field contains the gradient of the density function in the self-consistent way.  In the setting of a diffusion process $X_t$,
the density function $p_t$ of $X_t$ is given by the Fokker-Planck equation
{\small\begin{equation}\label{eqn:BFPE}
    \begin{aligned}
        \frac{\partial p_t(x)}{\partial t}
        =&\ -\nabla\cdot\left[
        \underbrace{ \left(b(x,t) -\frac{1}{2}\nabla\cdot\Sigma(x,t)- \frac{1}{2}\Sigma(x,t)\nabla \log p_t(x)\right)}_{v(x,t)}p_t(x)\right],
    \end{aligned}
\end{equation}
}where  $b$ is the drift term and $\Sigma$ corresponds to the diffusion matrix. If the term $\nabla \log p_t$, denoted by the score function $s_t(x)$, is {\it known a prior}, then \eqref{eqn:BFPE}
takes the form of a  continuity equation corresponds to the  transportation between the density functions   $p_0$ and $p_T$ by 
a underlying  deterministic ODE flow  with the new vector field  $v(x,t):=b(x,t)-1/2\nabla\cdot\Sigma(x,t)-1/2\Sigma(x,t)s_t(x)$. This observation is also the theoretic foundation of reversing a Markov diffusion process in score-based diffusion models since the time reverse from $p_T$ to $p_0$ can implemented by the time-reversed deterministic flow. 
However,
the various numerical methods \cite{hyvarinen2005estimation, song_generative_2019,song2020score,song2021denoising} to find the score function $s_t$  for the generative modeling like image synthesis are   limited to the   constant drift or
the Orenstein-Ulenbeck process, because  these simple processes  suffice to inject Gaussian noise to the data.  
In such simple systems,  the transition probability is analytically available.

 Motivated by probability flow and score-based generative models, several studies \cite{maoutsa2020interacting, shen2022self,li2023selfconsistent,lu2024score,boffi2023probability,shen2024entropy} have recently been carried out to learn the \emph{underlying velocity field} $v(x,t)$ that drives the probability density function in   Fokker-Planck equations, instead of 
 directly simulating the diffusion SDE or solving the Fokker-Planck PDE.  These methods, usually referred to as  {\it score-based transport modeling} \cite{boffi2023probability},  deterministically push samples from the initial density onto samples from the solution at any later time by employing particle methods to learn the score function from the samples following the deterministic flow. By finding the score function appearing in the new velocity field, 
 this numerical strategy works well due to the self-consistency of the Fokker-Planck equation \cite{li2023selfconsistent} in the form of   fixed-point argument.   Compared with the traditional 
grid-based numerical methods for the Fokker--Planck equation,
these score-based methods can scale well in high dimensions. Compared to the traditional numerical SDE schemes, the deterministic flow is faster to compute and   offers the expression of the density function 
to facilitate  calculating many non-equilibrium statistical quantities like the probability current,  the entropy production rate and the heat \cite{tome2014book,boffi2024deep}.  
The theoretic  works related to these algorithms include the error analysis for Kullback-Leibler divergence due to the approximate score function \cite{shen2024entropy,shen2022self,lu2024score,boffi2023probability}.  There are also substantial extensions from the Fokker-Planck equation to the McKean-Vlaso equation   \cite{maoutsa2020interacting, shen2024entropy,lu2024score}
 and the Fokker-Planck-Landau equation\cite{Ilin2024TransportBP,huang2024score}.

However,  the existing   score-based transport modeling only  studied    the diffusion process with the focus on the conventional  score function $s_t=\nabla \log p_t$ and the score matching error 
$\int_{t_0}^T\int_\Omega | s^\theta(x,t)-\nabla \log p_t(x) |^2 p_t(x) \d x \d t$   as the loss function to train the score.
Despite the popularity in applications of  the diffusion model, the  non-Gaussian processes such as L\'{e}vy processes   are preferred  to more accurately represent complex stochastic events than conventional Gaussian models \cite{kanazawa2020loopy, barthelemy2008levy, song2018neuronal}. Recently, the L\'{e}vy-It\^{o} model with the {\it fractional} score function $s_t^{(\alpha)}(x)=(\triangle ^{\frac{\alpha-2}{2}} \nabla p_t(x))/p_t(x)$, $\alpha\in (1,2]$,  is used  in \cite{yoon2023score} for score-based generative models to generate more diverse  samples  by taking advantage of the heavy-tailed properties of the $\alpha$-stable  L\'{e}vy process.  In \cite{hu2024score}, the same fractional score function is computed for the L\'{e}vy–Fokker–Planck equation by training the fractional score-based  physics-informed neural networks from a large number of  sample stochastic trajectories as the training data.

In this paper, we systematically investigate the score-based transport modeling for the general Markov process of non-Gaussian type and develop the score-based theories and algorithms to simulate the underlying jump-diffusion process.
This paper delivers mathematical derivation, numerical algorithm, and error analysis focusing on the corresponding score function in non-Gaussian systems with jumps and discontinuities represented by the nonlinear L\'{e}vy--Fokker--Planck equations.  We propose the    \emph{L\'{e}vy score function} for such stochastic equations, which features a nonlocal double-integral term, and its training  algorithm  by minimizing the proposed loss function from samples.  Based on the equivalence of the probability flow with deterministic dynamics, we develop a self-consistent score-based transport particle algorithm to sample the interactive L\'{e}vy stochastic process at discrete time grid points.

{\bf Contributions}
\begin{enumerate}
    \item 
    For the nonlinear stochastic differential equation  \eqref{sdeBL}  driven by the Brownian and L\'{e}vy noise as well as the mean-field interaction,
    we derive the continuity equation for the nonlinear L\'{e}vy-Fokker-Planck equation  \eqref{eqn:FPE}   to identify the vector field in the probability flow. We propose the \emph{generalized score function}, $\widehat{S}_t$ in \eqref{eqn:GS}, consisting of both the conventional diffusion score function $\nabla \log p_t$ and the  \emph{L\'{e}vy score function}. The latter has the non-local double-integral terms involving the L\'{e}vy measure.

\item
We  show in Theorem \ref{theo:KL} that, 
the Kullback--Leibler (KL) divergence $\sup_{t\in[t_0,T]}D_{\mathrm{KL}}(p^\theta_t\| \widehat{p}_t)$ between the  numerical score-based $p^\theta_t$ and the true density function $\widehat{p}_t$  is  bounded by a loss function  we derived from the fixed-point  argument.
This is a significant extension to   \cite{shen2024entropy,shen2022self,boffi2023probability,lu2024score} since our analysis takes into account of both the non-Gaussian process and the mean-field interactions.

\item 
We propose the score-matching training algorithm to learn both the diffusion and L\'{e}vy score functions by one single neural network, then develop  the sequential score-based particle algorithm ({\bf Algorithm}~\ref{algorithm}) to solve the  L\'{e}vy–Fokker–Planck equation based on the idea of score-based transport modeling. We test several numerical examples to demonstrate the effectiveness of our approach. 

\item 
 We prove in Theorem \ref{theo:error} the numerical error of {\bf Algorithm}~\ref{algorithm}  for solving  \eqref{sdeBL}:  $ \frac{1}{T} \mathbb{E}  \left|X^*_{n \Delta t}-X^N_{n \Delta t}\right| \leq  \mathcal{O}(\varepsilon)+\mathcal{O}(\Delta t),$ with $N$ the number of particles, $\epsilon$ the bounds of score loss, and $\Delta t$ the time step size,   $ n\Delta t \le T$.   
\end{enumerate}

\smallskip
{\bf Related works:}

{\bf Score-based diffusion model and the error analysis.}
The score-based diffusion model is formulated as the forward and backward continuous  It\^{o} stochastic differential equations in \cite{song2020score}. 
The score matching to train the score function by neural networks includes \cite{hyvarinen2005estimation,song2020sliced,ho_2020_denoising,song2021denoising,gottwald2024stable1,liu2024diffusion}. For comprehensive reviews on related topics, we refer readers to \cite{gallon2024overview}.
There are also extensive works to analyze the convergence and error estimates, for instance   \cite{NEURIPS2021_940392f5,chen_improved_2022,cao2024exploring}.

{\bf  Score-based transport modeling (SBTM) for Fokker-Planck equation.}  These ~methods are all based on  \eqref{eqn:BFPE} -- the equivalence to  probability flow --  to  simulate  the time
evolution of Fokker--Planck solutions.
\cite{maoutsa2020interacting}   runs a large number of  interacting deterministic particles  
where the gradient–log–density (score) function is approximated in the reproducing kernel Hilbert space (RKHS).
\cite {shen2022self,shen2024entropy} formulated a fixed-point problem for  the self-consistency of the probability flow for the Fokker--Planck equation,
but the optimization for the diffusion score function is based on the adjoint method which is in general slow and challenging to implement. 
The score-based transport modeling in \cite{boffi2023probability} is to march the flow sequentially in time with the forward Euler scheme. \cite{lu2024score} extends this method to McKean--Vlasov SDEs with the assumption of bounded interaction potential. \cite {li2023selfconsistent} 
proposed to learn the velocity field  $v(t,x)$ as a whole for $t\in [0,T]$ by an iterative update of the velocity field.

{\bf Score-based generative
models with $\alpha$-stable L\'{e}vy processes}.
 \cite{yoon2023score} proposed  score-based generative
models with a simple   forward L\'{e}vy process $X_t=a(t)X_0+\gamma(t)\epsilon$, 
where $X_0$ follows the data distribution, and the noise $\epsilon$ follows the isotropic $\alpha$-stable distribution with the scalar parameter $1$. 
The  fractional score function $s_t^{(\alpha)}(x)=(\triangle ^{\frac{\alpha-2}{2}} \nabla p_t(x))/p_t(x)$ with  $1<\alpha \le 2$,
is proposed for the pdf  $p_t(x)$ in order to define the time-reversal L\'{e}vy process. This fractional score function is approximated by a neural network and simple enough to use the denoising score-matching similarly to the score-based diffusion models.

{\bf Fractional score-based physics-informed neural networks(PINN) for $\alpha$-stable L\'{e}vy processes .}
 \cite{hu2024score} proposed the score-fPINN  method for the log-likelihood $q_t:=\log p_t$ of an  $\alpha$-stable L\'{e}vy process
 by minimizing the PINN loss for the two nonlinear PDEs of $q_t$ 
 and   $s_t=\nabla q_t$ respectively, which involve the pairs $q$ and $s^{(\alpha)}$,    $s^{(\alpha)}$ and $s$, respectively, where $s^{(\alpha)}$ is the same fractional score defined as in \cite{yoon2023score}.  
 The fractional score $s^{(\alpha)}$ serves only as  an intermediate 
variable to compute from $s_t$ to $q_t$. 
By the  score-matching method applicable to {\it any} stochastic process \cite{song2020sliced}, the diffusion score $s_t\approx \nabla \log p_t$ is computed by the neural network, entirely based on the {\it simulated long trajectories} of the SDE beforehand.

{\bf Comparison with Monte Carlo method.}
The Monte Carlo (MC) method  via SDE simulation  is the classic   numerical approach to generating trajectory samples. 
Although the per-particle, forward-time cost of each MC step is computationally cheap, the method is primarily suitable for simple, sample-path-driven statistics.
For more advanced, density-centric tasks (density values, probability currents, entropy production \cite{huang2025entropy,boffi2024deep}),  
MC requires global density estimation or interpolation in high-dimensional spaces. This process is statistically and computationally expensive and often becomes inaccurate at scale.
In contrast, the Score-Based Transport Method (SBTM), including our present extension to interactive L\'{e}vy processes, 
   is typically more efficient overall  by  replacing  repeated high-dimensional density estimation in MC with a learned score function and deterministic flow, giving robust and principled KL control. Although SBTM incurs an optimization overhead that scales with network size and the number of training samples, it is particularly advantageous when the objective extends beyond sampling to understanding the underlying probability dynamics, such as long-time and non-equilibrium behaviors.

\smallskip

Finally, we provide some remarks in the following to further explain a few other aspects of our works.

\begin{itemize}
   
\item  
We consider  a broader range of L\'{e}vy processes than  the $\alpha$-stable L\'{e}vy process  studied  in literatures such as \cite{yoon2023score} and \cite{hu2024score}.
Our model  \eqref{sdeBL} here is far beyond this simple L\'{e}vy process, and the corresponding score function is much more intricate than the fractional score function.    
Refer to Section \ref{ss:Lg}.

\item Our algorithm and analysis cover the systems with mean-field interactions as in \cite{lu2024score,shen2024entropy}.  
Our theorems apply to the interaction kernels of   the bounded  interaction and the Biot-Savart interaction as in \cite{shen2024entropy} for the Fokker--Planck equation. But due to the challenges from the L\'{e}vy score, we have not yet prove for the Column interaction as the diffusion process in \cite{shen2024entropy}.

    \item     The  generalized  score function trained in {\bf Algorithm} \ref{algorithm}  contains the contributions of drift, interaction, diffusion and jump;  thus, the transport
    for the $N$ deterministic particles are independent of each other, while   \cite{lu2024score} still evolves  the $N$  interactive    particles  with the  mean-field interaction.
\item Even though our numerical algorithm to train the score is sequential in time   as in \cite{boffi2023probability}, 
the approach in \cite{li2023selfconsistent} of the self-consistent velocity matching  by learning from  the whole simulated trajectories until a  terminal time is applicable to 
the L\'{e}vy--Fokker--Planck equation here, because the same fixed-point formulation works with the new defined L\'{e}vy  score function in our work.

\end{itemize}

\medskip

The structure of this paper is as follows: Section \ref{sec:GSF} shows how to define the generalized score function (the diffusion score and the L\'{e}vy score) and establish the probability flow equivalence. Section \ref{sec:anaylis} introduces a particle system-based numerical method, presents bounds on the Kullback--Leibler divergence in Theorem \ref{theo:KL}, and analyzes discretization errors in Theorem \ref{theo:error}. Section \ref{sec:example} demonstrates the method through numerical experiments.
The conclusion part is summarized in Section \ref{sec:conclusions}.

    \section{Probability Flow of Jump-Diffusion and L\'{e}vy Score Function}
    \label{sec:GSF}
    \subsection{ Mckean--Vlasov SDE driven by L\'{e}vy noise}

We consider the following Mckean--Vlasov stochastic differential equation that
incorporates non-Gaussian noise: {\small\begin{equation}\label{sdeBL}\left\{\begin{aligned}\d X_{t}=&\ b(X_{t-},t)\d t + (K*\widehat{p}_{t})(X_{t-})\d t + \sigma (X_{t-},t)\d B_{t}+ \int_{|r|<1}F(r,t) \tilde{\mathcal{N}}(\d t,\d r) \\&\ +\int_{|r|\geq1}G(r,t) \mathcal{N}(\d t,\d r),\quad \widehat{p}_{t}=\mathrm{Law}(X_{t}),\quad X_{t}\in\Omega,\quad t> t_{0},\\ X_{t_0}\sim&\ \mu_{0}\in\mathcal{P}_{2}(\Omega),\end{aligned}\right.\end{equation} }where
$\mathcal{P}_{2}(\Omega)$ denotes the set of probability measures in $\Omega$
with finite quadratic moment and $\widehat{p}_{t}$ denotes the law of $X_{t}$. The
domain $\Omega$ is either $\mathbb{R}^{d}$ or the torus $\mathbb{T}^{d}$ (a
hypercube $[-L,L]^{d}$ with periodic boundary condition), $B$ is a standard
Brownian motion in $\mathbb{R}^{m}$ and $\mathcal{N}$ is an independent Poisson
random measure on $\mathbb{R}_{+}\times(\mathbb{R}^{d}\backslash\{0\})$ with the
associated compensator $\tilde{\mathcal{N}}$ and the intensity measure $\nu$
which is an arbitrary L\'{e}vy measure. Specifically, $\nu$ satisfies the condition $\int_{|r| > 0}\big(|r|^{2}\wedge 1\big)\nu(\mathrm{d}r) < \infty$, and the compensator is given by $\tilde{\mathcal{N}}(\mathrm{d}t, \mathrm{d}r) = \mathcal{N}(\mathrm{d}t, \mathrm{d}r) - \nu(\mathrm{d}r)\,\mathrm{d}t$. For instance, if $\mathcal{N}$ corresponds to a compound Poisson process with Poisson parameter $\lambda_{0}$ and jump sizes distributed according to $f_{A}$, then the associated L\'{e}vy measure is given by $\nu(\mathrm{d}r) = \lambda_{0}f_{A}(r)\,\mathrm{d}r$. 
The mappings $b_{i}:\mathbb{R}^{d}\times[t_{0},\infty)\to\mathbb{R}$, $\sigma_{ij}
:\mathbb{R}^{d}\times[t_{0},\infty)\to\mathbb{R}$,
$\sigma=(\sigma_{ij})_{d\times m}$, $\Sigma=\sigma\sigma^{\top}$,
$F_{i}:\mathbb{R}^{d}\times[t_{0},\infty)\to\mathbb{R}$, $G_{i}:\mathbb{R}^{d}\times
[t_{0},\infty)\to\mathbb{R}$, all assumed to be measurable for $1\leq i \leq d$
and $1\leq j\leq m$. $K:\mathbb{R}^{d}\to\mathbb{R}^{d}$ denotes some interaction
kernel and $*$ denotes the convolution operation, i.e.,
$(K*\widehat{p}_{t})(X_{t})=\int_{\Omega}K(X_{t}-y)\widehat{p}_{t}(y)\d y$. We introduce
the following assumptions on the drift term $b$, diffusion coefficient $\sigma$,
the jump noise coefficients $F,G$ and the L\'{e}vy measure $\nu$.
\begin{assumption}
    \label{ass:allcoefficient} \hfill
    \begin{enumerate}
        \item For every $t \geq t_{0}$, $b_{i}(\cdot, t)$ and
            $\sigma_{ij}(\cdot, t)$ are $C^{2}$-measurable functions for
            $1 \leq i \leq d$ and $1 \leq j \leq m$. There exists positive constants
            $C_{b},\ C_{\sigma}$ such that $|b(t,x)-b(t,y)|^{2}\leq\ C_{b}|x-y|^{2}
            ,\ (C_{\sigma})^{-1}\leq\ \Sigma_{ij}(x,t) \leq C_{\sigma}, \quad \forall
            x,y \in \Omega, \quad t \geq t_{0}, \quad i,j = 1, \ldots, d,$\label{subass:boundsigmaK}
            where $|\cdot|$ denotes the Euclidean norm in $\mathbb{R}^{d}$.

        \item The functions $F,G\in C(\mathbb{R}^{d}\times[t_{0},\infty),\mathbb{R}
            ^{d})$, the L\'{e}vy measure $\nu$ satisfies $\int_{\mathbb{R}\backslash\{0\}}
            (1\wedge|r|^{2})\nu(\d r) < \infty$, and there exist constants
            $C_{T}^{F}$ and $C_{T}^{G}$ dependent on $T$ such that
            $\int_{|r|<1}|F(r,t)|^{2}\nu(\d r) < C_{T}^{F},\quad \int_{|r|\geq 1}
            |G(r,t)|^{2}\nu(\d r) < C_{T}^{G},\quad \forall t\in[t_{0},T]$.
    \end{enumerate}
\end{assumption}
Our analysis will focus on the following two cases of the interaction kernel $K$.
\begin{enumerate}
    \item ({\bf Bounded interaction}) The kernel $K\in C^{2}(\mathbb{R}^{d})$ and
        there exists a positive constant $C_{K}$ such that $|K(x)| \leq C_{K},\ \forall
        x\in\Omega$.\label{subass:bounded}

    \item ({\bf Biot--Savart interaction}) The dimension $d=2$ and the interaction
        term $K$ is the Biot--Savart kernel, i.e., $K(x)=\frac{1}{2\pi}\left(-\frac{x_{2}}{
        |x |^{2}},\frac{x_{1}}{ |x |^{2}}\right)^{\top}$ where $x=(x_{1},x_{2})^{\top}$.\label{subass:BS}
        This kernel is used to describe the dynamics of electromagnetism \cite{golse2016dynamics}.
\end{enumerate}

For the bounded interaction case \ref{subass:bounded}, under Assumption
\ref{ass:allcoefficient}, the well-posedness of the SDE \eqref{sdeBL} can be
established by applying the interlacing technique to exclude large jumps (defined
as exceeding 1) \cite[Theorem 4]{liu2023stochastic}. For the Biot--Savart
interaction case \ref{subass:BS}, to the best of our knowledge, there is no
literature rigorously proving the well-posedness of \eqref{sdeBL}. However, we do
not focus on this issue in the present paper, and thus omit further discussion. To
manage possible singularities in the interaction kernel $K$, one can consider
the regularized form, for instance, for the potential, the modification $K^{\varepsilon}
(x)=\frac{1}{2\pi}\left(-\frac{x_{2}}{ |x |^{2}+\varepsilon},\frac{x_{1}}{ |x |^{2}+\varepsilon}
\right)^{\top}$ ensures $K^{\varepsilon}$ is bounded for every $\varepsilon>0$.

\subsection{ L\'{e}vy--Fokker-Planck equation and L\'{e}vy score function}
The probability flow serves as a crucial link between the forward Kolmogorov equation
and the continuity equation for the density function, as illustrated in \eqref{eqn:BFPE}.
It effectively handles (nonlinear) diffusion processes through an ODE flow
framework facilitated by a well-defined score function. In the following, we
investigate the counterpart to the conventional diffusion score function for
\eqref{sdeBL}, which we refer to as the L\'{e}vy score function. To identify
this score function, we first need to convert the L\'{e}vy-Fokker-Planck
equation into a continuity equation. Our formulation of the L\'{e}vy-Fokker-Planck
equation has appeared in our other work \cite{huang2025probability}. The result is
more general than the $\alpha$ stable L\'{e}vy process considered in
\cite{hu2024score,yoon2023score}, and has no explicit term of the fractional Laplacian.
We start from the infinitesimal generator\cite{applebaum2009levy,liang2021exponential}
for \eqref{sdeBL}: {\footnotesize\begin{equation}\label{generator}\begin{aligned}&(\mathcal{L}[\widehat{p}_{t}]f)(x)=\ \sum_{i=1}^{d}\left( b_{i}(x,t) + (K_{i}*\widehat{p})(x)\right)(\partial_{i}f)(x) + \int_{|r|\geq 1}\left[ f(x+G(r,t)) - f(x) \right]\nu(\d r)\\&+ \int_{|r|<1}\bigg[ f(x+F(r,t)) - f(x) - \sum_{i=1}^{d}F_{i}(x,z)(\partial_{i}f)(x)\bigg]\nu(\d r) + \frac{1}{2}\sum_{i,j=1}^{d}\Sigma_{ij}(x,t)(\partial_{i}\partial_{j}f)(x) ,\end{aligned}\end{equation} }for
each function $f \in C^{2}_{0}(\mathbb{R}^{d})$ and each point
$x \in \mathbb{R}^{d}$. We apply the Taylor's theorem to $f(x + F(r,t))$ and $f(x
+ G(r,t))$ at the point $x$:
$f(x+F(r,t))= f(x) + \sum_{i=1}^{d}F_{i}(r,t)\int_{0}^{1}(\partial_{i}f)(x+\lambda
F(r,t))\d \lambda$. Substitute it to \eqref{generator}, and compute the adjoint operator
of $\mathcal{L}[\widehat{p}_{t}]$, with more details in
\cite{huang2025probability}, we then have the following form of the L\'{e}vy--Fokker-Planck
equation associated with \eqref{sdeBL}: {\small\begin{equation}\label{eqn:FPE}\begin{aligned}\frac{\partial \widehat{p}_t(x)}{\partial t}=&-\sum_{i=1}^{d}\ \frac{\partial}{\partial x_i}\left[\left(b_{i}(x,t) + (K_{i}*\widehat{p}_{t})(x) - \int_{|r|<1}F_{i}(r,t)\nu(\d r)\right)\widehat{p}_{t}(x) \right] \\&- \int_{|r|<1}\int_{0}^{1}\sum_{i=1}^{d}\frac{\partial}{\partial x_i}\left( F_{i}(r,t) \widehat{p}_{t}(x-\lambda F(r,t) ) \right) \d \lambda \nu(\d r) \\&- \int_{|r|\geq1}\int_{0}^{1}\sum_{i=1}^{d}\frac{\partial}{\partial x_i}\left( G_{i}(r,t) \widehat{p}_{t}(x-\lambda F(r,t) ) \right)\d \lambda \nu(\d r) + \frac{1}{2}\sum_{i,j=1}^{d}\frac{\partial^2(\Sigma_{ij}(x,t)\widehat{p}_t(x)) }{\partial x_i\partial x_j}.\end{aligned}\end{equation} }
which can be rewritten as the {\emph{continuity equation}}:
{\small\begin{equation}\label{eqn:FPE-transpot}\begin{aligned}&\frac{\partial \widehat{p}_t(x)}{\partial t}=\ -\sum_{i=1}^{d}\ \frac{\partial }{\partial x_i}\left[ \left(b_{i}(x,t) + (K_{i}*\widehat{p}_{t})(x) - \int_{|r|<1}F_{i}(r,t)\nu(\d r) \right.\right. \\&\ + \int_{|r|<1}\int_{0}^{1}\frac{ F_i(r,t) \widehat{p}_t(x - \lambda F(r,t)) }{\widehat{p}_t(x)}\d \lambda \nu(\d r) + \int_{|r|\geq1}\int_{0}^{1}\frac{ G_i(r,t) \widehat{p}_t(x - \lambda G(r,t)) }{\widehat{p}_t(x)}\d \lambda \nu(\d r)\\&\ \left.\left. - \frac{1}{2}\sum_{j=1}^{d}\frac{\partial}{\partial x_j}\Sigma_{ij}(x,t) - \frac{1}{2}\sum_{j=1}^{d}\Sigma_{ij}(x,t)\frac{\frac{\partial}{\partial x_j} \widehat{p}_t(x)}{\widehat{p}_t(x)}\right)\widehat{p}_{t}(x) \right] =: - \nabla \cdot \left( \mathcal{V}[\widehat{p}_{t}](x,t) \widehat{p}_{t}(x) \right ).\end{aligned}\end{equation} }Based
on \eqref{eqn:FPE-transpot}, the solution density $\widehat{p}_{t}$ to the L\'{e}vy-Fokker-Planck
equation \eqref{eqn:FPE} can be interpreted as the pushforward of $\mu_{0}$ under
the flow map $X^{*}_{s,t}$ governed by {\small\begin{equation}\label{eqn:interactingODE}\begin{aligned}\frac{\d X^*_{s,t}(x)}{\d t}=\ \mathcal{V}[\widehat{p}_{t}] (X^{*}_{s,t}(x),t),\quad X^{*}_{s,s}(x)=x,\quad t\geq s\geq t_{0},\end{aligned}\end{equation} }referred
to as the {\it probability flow equation}. Here the vector field depends on the $\widehat
{p}_{t}$ in a self-consistent way:
{\small\begin{equation}\label{V}\begin{aligned}&\mathcal{V}[\widehat{p}_{t}](x,t):= b(x,t) + (K*\widehat{p}_{t})(x) - \int_{|r|<1}F(r,t)\nu(\d r) - \frac{\nabla \cdot\Sigma(x,t)}{2}- \sum_{j=1}^{d}\frac{\Sigma(x,t)}{2}\nabla\log \widehat{p}_{t}(x)\\&\ + \int_{|r|<1}\int_{0}^{1}\frac{ F(r,t) \widehat{p}_t(x - \lambda F(r,t)) }{\widehat{p}_t(x)}\d \lambda \nu(\d r) + \int_{|r|\geq1}\int_{0}^{1}\frac{ G(r,t) \widehat{p}_t(x - \lambda G(r,t)) }{\widehat{p}_t(x)}\d \lambda \nu(\d r) \\&=: b(x,t) + (K*\widehat{p}_{t})(x) - \int_{|r|<1}F(r,t)\nu(\d r) - \frac{1}{2}\nabla \cdot\Sigma(x,t) - \widehat{S}[\widehat{p}_{t}] (x,t).\end{aligned}\end{equation} }
In the above, we define the \emph{generalized} score function $\widehat{S}$ as
follows: {\small\begin{align}\label{eqn:GS}\widehat{S}[\widehat{p}_{t}](x,t) :=&\ \frac{1}{2}\sum_{j=1}^{d}\Sigma(x,t)\nabla\log \widehat{p}_{t}(x) +\widehat{S}_{L}[\widehat{p}_{t}](x,t) \\ \label{eqn:GSL}\widehat{S}_{L}[\widehat{p}_{t}](x,t) :=&- \int_{|r|<1}\int_{0}^{1}\frac{ F(r,t) \widehat{p}_{t}(x - \lambda F(r,t)) }{\widehat{p}_{t}(x)}\d \lambda \nu(\d r) \\&-\int_{|r|\geq1}\int_{0}^{1}\frac{ G(r,t) \widehat{p}_{t}(x - \lambda G(r,t)) }{\widehat{p}_{t}(x)}\d \lambda \nu(\d r). \notag\end{align} }The
gradient term $\nabla\log\widehat{p}_{t}$ in $\widehat{S}$ of \eqref{eqn:GS} is recognized
as the conventional score function derived from the Laplacian terms associated with
Gaussian noise, which will be referred to as {\bf diffusion score}. $\widehat{S}_{L}$
in \eqref{eqn:GSL} consist of non-local double integral terms, arising from jump
noise, which is called the {\bf L\'{e}vy score} function.

\subsection{Isotropic \texorpdfstring{$\alpha$}{Lg}-stable L\'{e}vy process}
\label{ss:Lg} The fractional score function
$s_{t}^{(\alpha)}(x)=((-\triangle)^{\frac{\alpha-2}{2}}\nabla p_{t}(x))/p_{t}(x),$
$\alpha\in (1,2]$, is proposed in L\'{e}vy-type score-based generative modeling
\cite{yoon2023score} and used in \cite{hu2024score} to solve the L\'{e}vy–Fokker–Planck
equation through physics-informed neural networks. This fractional score
function is the special case of our result above. Consider the special $\alpha$-stable
L\'{e}vy processes: $\d X_{t}= \ b(X_{t-},t)\d t + \sigma (X_{t-},t) \d B_{t}+ \sigma
_{L}(t)\d L^{\alpha}_{t}$ which fall within our diffusion-jump SDE \eqref{sdeBL}
with $K=0$, $F(r,t)=G(r,t)=r\sigma_{L}(t)$, and $\nu(\d r)$ is the $\alpha$-stable
L\'{e}vy measure $\nu_{\alpha}(\d r)\sim |r|^{-d-\alpha}\d r$; note that
$\int_{|r|<1}F(r,t)\nu (\d r)=0$. In this special case, the L\'{e}vy--Fokker-Planck
equation \eqref{eqn:FPE} reads {\footnotesize\begin{equation*}\begin{aligned}\frac{\partial \widehat{p}_t(x)}{\partial t}=&-\sum_{i=1}^{d}\ \frac{\partial\left[ b_i(x,t) \widehat{p}_t(x) \right]}{\partial x_i}+ \sum_{i,j=1}^{d}\frac{\partial^2(\Sigma_{ij}(x,t)\widehat{p}_t(x)) }{2\partial x_i\partial x_j}+ \int_{\mathbb{R}^d\backslash\{0\}}\left(p(x - \sigma_{L}(t)r) -p(x)\right)~\nu_{\alpha}(\d r) \\ =&-\sum_{i=1}^{d}\ \frac{\partial\left[ b_i(x,t) \widehat{p}_t(x) \right]}{\partial x_i}+ \sum_{i,j=1}^{d}\frac{\partial^2(\Sigma_{ij}(x,t)\widehat{p}_t(x))}{2\partial x_i\partial x_j}-\sigma_{L}^{\alpha}(t) (-\triangle)^{\frac{\alpha}{2}}\widehat{p}_{t}(x)\end{aligned}\end{equation*}}where
$(-\triangle)^{\frac{\alpha}{2}}$ denote the fractional Laplacian of order
$\frac{\alpha}{2}$ \cite{caffarelli2007extension}, i.e., $(-\triangle)^{\frac{\alpha}{2}}
\widehat{p}_{t}(x)= C_{d,\alpha}\int_{\mathbb{R}^d}\frac{\widehat{p}_{t}(y) - \widehat{p}_{t}(x)}{|x-y|^{d+\alpha}}
\d y,$ with a normalizing constant $C_{d,\alpha}$. Note that
\cite[Lemma C.1]{yoon2023score}
$(-\triangle)^{\frac{\alpha}{2}}\widehat{p}_{t}(x)= -\nabla\cdot \left( (-\triangle
)^{\frac{\alpha-2}{2}}\nabla \widehat{p}_{t}(x) \right)= - \nabla \cdot ( s^{(\alpha)_t(x)}
\widehat{p}_{t}(x)).$
So the L\'{e}vy score function we define in \eqref{eqn:GSL} reduces to
$\sigma_{L}^{\alpha}(t) s^{(\alpha)}_{t}(x)$.

\subsection{Probability flow for probability density function}
The remarkable property of the continuity equation \eqref{eqn:FPE-transpot} is the
following: The flow associated with the deterministic ordinary differential
equation \eqref{eqn:interactingODE} transports the initial distribution $\mu_{0}$
forward exactly along the density distribution $\widehat{p}_{t}$ of the SDE \eqref{sdeBL}.
To turn this idea into a feasible numerical algorithm certainly needs to obtain the
generalized score $\widehat{S}$ or the velocity $\mathcal{V}(x,t)$ in \eqref{V}
beforehand, which will be studied in next section. If $x$ is a sample from the
initial distribution $\mu_{0}$, then $X^{*}_{t_0,t}(x)$ in \eqref{eqn:interactingODE}
will be a sample from $\widehat{p}_{t}$. We use $\#$ to denote the push-forward operation,
then $\widehat{p}_{t}= X^{*}_{t_0,t}\# \mu_{0}$ can be determined at any
position using the change of variables formula
\cite{santambrogio2015optimal,villani2009optimal}: {\small\begin{equation}\label{eqn:fixpoint}\widehat{p}_{t}(x)=\mu_{0}(X^{*}_{t,t_0}(x))\exp\left(-\int_{t_0}^{t}\nabla\cdot\mathcal{V}[\widehat{p}_{t}](X^{*}_{t,s}(x),s)\d s\right),\end{equation} }where
$\nabla \cdot \mathcal{V}$ is the divergence of the velocity field defined in
\eqref{V}. It is straightforward to observe that $\widehat{p}_{t}$, the solution
of the L\'{e}vy--Fokker-Planck equation \eqref{eqn:FPE}, is the fixed point of the
map $p_{t}\mapsto \mu_{0}(X^{*}_{t,t_0}(x))\exp\left(-\int_{t_0}^{t}\nabla\cdot\mathcal{V}
[p_{t}](X^{*}_{t,s}(x),s)\d s\right).$

By evolving an ensemble of $N$ independent realizations of \eqref{eqn:interactingODE},
referred to as ``particles'', according to
\begin{equation}
    \label{eqn:ODEX*}\frac{\d X^{N,i}(t)}{\d t}=\ \mathcal{V}[\widehat{p}_{t}](X^{N,i}
    (t),t),\ i=1,\cdots,N,\quad X^{N,i}(t_{0})\sim\mu_{0},
\end{equation}
an empirical approximation to $\widehat{p}_{t}$ is obtained:
$\widehat{p}_{t}(x)\approx p^{N}_{t}(x)=\ \frac{1}{N}\sum_{i=1}^{N}\delta_{x-X^{N,i}(t)}
.$

In \cite{lu2024score}, the deterministic flow for $N$ interacting particles is formulated by explicitly including the interaction term $K$ : {\small\begin{equation}\label{eqn:IODEX*}\frac{\d X^{N,i}(t)}{\d t}=\ \mathcal{V}_{o}[p^{N}_{t}](X^{N,i}(t),t)+ \frac{1}{N}\sum_{j}K(X^{N,i}(t),X^{N,j}(t)),\ i=1,\cdots,N,\quad X^{N,i}(t_{0})\sim\mu_{0},\end{equation} }where $\mathcal{V}_{o}=\mathcal{V}- K*p^{N}_{t}$ is defined by subtracting the mean-field interaction from the velocity field $\mathcal{V}[\widehat{p}{t}]$ in our definition \eqref{V}. These two formulations, \eqref{eqn:ODEX*} and \eqref{eqn:IODEX*}, as the probability flow of stochastic interactive models, are mathematically equivalent in the mean field limit $N\to \infty$. The principal numerical distinction lies in how computational costs are allocated between the inference step (learning scores) and the transport step (propagating particles). We develop our theoretical results primarily based on formulation \eqref{eqn:ODEX*}, and provide a   discussion of the two resulting algorithms in Remark \ref{remark:algorithm} Section \ref{sec:algorithm}. 

    \section{Score-based Particle Approach: Theoretical and Numerical Analysis}
    \label{sec:anaylis}

Recall \eqref{eqn:fixpoint}, the general principle involves here addresses the following
fixed-point problem: For any given velocity field $V^{in}(x,t)$, the flow dictated
by the ODE \eqref{eqn:interactingODE} will transport the initial density $\mu$
to obtain $p_{t}$, and this transported $p_{t}$ furthermore induces the new velocity
field $V^{out}$ defined via \eqref{V}. It is evident that the true velocity field
$\mathcal{V}$ is the fixed point of this map $V^{in}\mapsto V^{out}$. Thus, if we
are provided with a set of vector fields $\{f^{\theta}\}_{\theta\in\Theta}$, and
obtain its corresponding probability flows $p^{\theta}$ via \eqref{eqn:interactingODE},
the ideal choice of these vector fields that approximate the true vector field
is the one that minimizes the following loss function with some samples from
$p^{\theta}_{t}$: $\int_{t_0}^{T}\int_{\Omega}|f^{\theta}(x,t)-\mathcal{V}[p^{\theta}
_{t}](x,t) |^{2}p^{\theta}_{t}(x)\d x\d t.$


Moving forward, we will first establish a theoretical foundation for our method
in Section \ref{sec:upperbound} using the above loss function. Following this,
we will develop a specific algorithm in Section \ref{sec:algorithm} and conduct
a detailed error analysis in Section \ref{sec:error}.

\subsection{The Kullback–Leibler divergence error}
\label{sec:upperbound}

For simplicity, we introduce the following notations
{\small\begin{align}\widehat{b}(x,t):=&\ b(x,t) - \int_{|r|<1}F(r,t)\nu(\d r) -\frac{1}{2}\nabla\cdot\Sigma(x,t),\label{notaionsbIA}\\ \mathcal{I}_{t}^{F,G}[\widehat{p}_{t}](x):=&\ \mathcal{I}_{t}^{F}[\widehat{p}_{t}](x)+\mathcal{I}_{t}^{G}[\widehat{p}_{t}](x),\label{notaionsIFG}\end{align}}
{\small\begin{align}\mathcal{I}_{t}^{F}[\widehat{p}_{t}](x):=&\ \int_{|r|<1}\int_{0}^{1}\frac{ F(r,t) \widehat{p}_{t}(x - \lambda F(r,t)) }{\widehat{p}_{t}(x)}\d \lambda \nu(\d r), \\ \mathcal{I}_{t}^{G}[\widehat{p}_{t}](x):=&\ \int_{|r|\geq1}\int_{0}^{1}\frac{ G(r,t) \widehat{p}_{t}(x - \lambda G(r,t)) }{\widehat{p}_{t}(x)}\d \lambda \nu(\d r).\end{align} }And
we recall the vector field \eqref{V}, which can be written as {\small\begin{equation}\label{def:mapA}\mathcal{V}[\mu](x,t):=\ \widehat{b}(x,t) + (K*\mu)(x) \\ \ - \frac{1}{2}\Sigma(x,t)\nabla\log \mu + \mathcal{I}_{t}(\mu,F,G).\end{equation} }

Now suppose we are given a family of time-varying hypothesis velocity fields
$\{f^{\theta},\ f^{\theta}:\mathbb{R}^{d}\times[t_{0},\infty)\to\mathbb{R}^{d}\}_{\theta\in\Theta}$
for some index set $\Theta$, and let $p^{\theta}$ be the solution to the continuity
equation {\small\begin{equation}\label{eqn:transport}\begin{aligned}\frac{\partial p^{\theta}_t(x)}{\partial t}=&\ -\nabla\cdot\left(f^{\theta}(x,t) p^{\theta}_{t}(x)\right),\quad p^{\theta}_{t_0}=\mu_{0}.\end{aligned}\end{equation} }Similar
with previous arguments, the solution density $p_{t}^{\theta}$ can be viewed as the
push forward of $\mu_{0}$ under the flow map $X^{\theta}_{s,t}$ of the ordinary differential
equation {\small\begin{equation}\label{eqn:ODEtheta}\begin{aligned}\frac{\d X^{\theta}_{s,t}(x)}{\d t}=\ f^{\theta}(X^{\theta}_{s,t}(x),t),\quad X^{\theta}_{s,s}(x)=x,\quad t,s\geq t_{0},\end{aligned}\end{equation} }and
thus the density function can be written as
{\small\begin{equation}\label{densitytheta}p^{\theta}_{t}(x)=\mu_{0}(X^{\theta}_{t,0}(x))\exp\left(-\int_{t_0}^{t}\nabla\cdot f^{\theta}(X^{\theta}_{t,s}(x),s)\d s\right),\end{equation} }as
discussed before. Let us make some assumptions on the the initial distribution
$\mu_{0}$, the L\'{e}vy-Fokker-Planck solution $\widehat{p}$, and the parametrized
vector fields $\{f^{\theta}\}_{\theta\in\Theta}$. In the remainder of this paper,
we will focus on the case where $\Omega$ is a torus. The method for extending
this to an unbounded domain is similar to the approach detailed in
\cite[Appendix G]{shen2024entropy}, which requires additional assumptions about the
regularity of the initial distribution $\mu_{0}$ and the true probability flow $\widehat
{p}_{t}$.

\begin{assumption}
    \label{ass:densitytheta} \hfill
    \begin{enumerate}
        \item The initial distribution $\mu_{0}$ is absolutely continuous with respect
            to the Lebesgue measure (we still denote its density as $\mu_{0}$) and
            there exists a positive constant $C_{0}$ such that $(C_{0})^{-1}\leq\mu
            _{0}(x)\leq C_{0},\quad \forall x\in\Omega,$ and $\mu_{0}\in C^{2}(\Omega
            )$;\label{subass:initiallaw}

        \item For every $T\geq0$, the solution to \eqref{eqn:FPE}$, \widehat{p}_{t}
            \in C^{3}(\Omega)$
            for every $t\geq t_{0}$ and there is a positive constant $C^{*}_{T}$
            such that
            $(C^{*}_{T})^{-1}\leq \widehat{p}_{t}(x)\leq C^{*}_{T},\ \forall x\in
            \Omega,\ t\in[t_{0},T]$.\label{subass:truedensity}

        \item The velocity field $f^{\theta}(\cdot,\cdot)\in C^{2,1}(\Omega\times
            \mathbb{R},\mathbb{R}^{d})$ for every $\theta\in\Theta$. For every
            $T \geq t_{0}$, there is a positive constant $C_{T}$ such that $\sup_{\theta\in\Theta}
            \sup_{(x,t)\in\Omega\times [t_0,T]}|\nabla \cdot f^{\theta}(x ,t)| \leq
            C_{T}$.\label{subass:KT}
    \end{enumerate}
\end{assumption}
\begin{remark}
    \label{eq:CTf} The Assumption \ref{ass:densitytheta}-\ref{subass:truedensity}
    can be achieved under some mild conditions. For further details, we refer to
    \cite{de2024multidimensional,olivera2024microscopic} as references. Under
    Assumption \ref{ass:densitytheta}-\ref{subass:KT}, for every $T\geq t_{0}$,
    it is easy to show that there exists a constant $C_{T}^{f}$ (dependent on $C_{T}$
    and $T$ only) such that
    $(C^{f}_{T})^{-1}\leq p^{\theta}_{t}(x)\leq C^{f}_{T},\ \forall x\in\Omega, t
    \in[t_{0},T],\ \theta\in\Theta$.
\end{remark}

We now give an upper bound for the KL divergence between $p^{\theta}$ and
$\widehat{p}$ for bounded kernel case \ref{subass:bounded} and Biot--Savart kernel
case \ref{subass:BS}.
\begin{theorem}
    \label{theo:KL} For interaction kernels \ref{subass:bounded} (bounded) and \ref{subass:BS}
    (Biot-Savart), let $p^{\theta}_{t}(x)$ denote the solution to the transport
    equation \eqref{eqn:transport}, and let $\widehat{p}_{t}(x)$ denote the
    solution to the L\'{e}vy--Fokker--Planck equation \eqref{eqn:FPE}. Suppose Assumption
    \ref{ass:allcoefficient} and \ref{ass:densitytheta} hold, then
    {\small\begin{equation}\label{ineqn:KLmix2}\begin{aligned}&\ \sup_{t\in[t_0,T]}D_{\mathrm{KL}}(p^{\theta}_{t}\| \widehat{p}_{t})\leq\  3\exp(\bar{C}T)C_{\sigma}\int_{t_0}^{T}\int_{\Omega}\left| f^{\theta}(x,t) - \mathcal{V}[p^{\theta}_{t}](x,t)\right|^{2}p^{\theta}_{t}(x)\d x\d t ,\end{aligned}\end{equation} }
    where $\bar{C}$ is a constant given in \eqref{eqn:barC}.
\end{theorem}

\subsection{Sequential L\'{e}vy score-based particle method}
\label{sec:algorithm}

By Theorem \ref{theo:KL}, we can minimize the following function to control the KL
divergence between the numerical $p^{\theta}_{t}$ and the true $\widehat{p}_{t}$:
$\min_{\theta}\int_{t_0}^{T}\int_{\Omega}\left| f^{\theta}(x,t) - \mathcal{V}[p^{\theta}
_{t}](x,t)\right|^{2}p^{\theta}_{t}(x)\d x\d t$, where $\mathcal{V}$ is defined in
\eqref{V} or \eqref{def:mapA}. By \eqref{eqn:interactingODE} and \eqref{V}, the velocity
field $f^{\theta}$ is given by {\small\begin{equation}\label{def:ftheta}f^{\theta}(x,t)=\  b(x,t) - \int_{|r|<1}F(r,t) \nu(\d r) -\frac{1}{2}\nabla\cdot\Sigma(x,t) - s^{\theta}(x,t),\end{equation} }where
$s^{\theta}$ is the numerical score function parametrized by $\theta$, {\small\begin{equation}\label{eqn:stheta}s^{\theta}\approx \widehat{S}[\widehat{p}_{t}] - K*\widehat{p}_{t},\end{equation} }where $\widehat{S}[\widehat{p}_{t}]$ is given in \eqref{eqn:GS}.
Now, the optimization problem turns to
{\small\begin{equation}\label{opt:origin}\begin{aligned}&\ \min_{\theta\in\Theta}\int_{t_0}^{T}\int_{\Omega}\left| f^{\theta}(x,t) - \mathcal{V}[p^{\theta}_{t}](x,t)\right|^{2}p^{\theta}_{t}(x)\d x\d t \\ =&\ \min_{\theta\in\Theta}\int_{t_0}^{T}\int_{\Omega}\left| s^{\theta}(x,t) + \left(K*p^{\theta}_{t}\right)(x) - \frac{1}{2}\Sigma(x,t)\nabla\log p^{\theta}_{t}(x) + \mathcal{I}_{t}^{F,G}[p^{\theta}_{t}](x)\right|^{2}p^{\theta}_{t}(x)\d x\d t.\end{aligned}\end{equation}}

This optimization over a score function $s^{\theta}:\mathbb{R}^{d}\times [t_{0}, T]\to \mathbb{R}^{d}$ can, in principle, be solved using the gradient-based method such as the adjoint method described in \cite{shen2024entropy} which may incur  substantial computational overhead.  \cite{boffi2023probability} proposed a sequential time-stepping procedure that  improves computational efficiency. In the following, we elaborate this sequential approach in details. Instead of directly learning $s^{\theta}(x,t):\mathbb{R}^{d}\times [t_{0}, T]\to \mathbb{R}^{d}$ on the whole time interval, we consider to learn $s^{\theta}(x,t)$ separately at each time sub-interval in a time-discrete way. Let $t_{0}< \cdots < t_{n}= T$ with $\Delta t = (T - t_{0})/n$. For {\it each} $[t_{i-1}, t_{i}]$, $1\le i\le n$, we use one neural network function of the variable $x$ {\it only}, denoted by $s^{\theta_i}(\cdot): \mathbb{R}^{d}\to\mathbb{R}^{d}$, to approximate the continuous-time score function $s^{\theta}(\cdot,t)$ {\it restricted} in this sub-interval.\footnote{One can construct the linear interpolation in time from the discrete time by $s^{\theta}(x,t)=(t-t_{i-1})\mathrm{NN}^{\theta_i}(x)+ s^{\theta_{i-1}}(x)$ on $t\in [t_{i-1},t_{i}]$. In addition, to learn $s^{\theta_i}(x)= \Delta t \times \mathrm{NN}^{\theta_i}(x)+ s^{\theta_{i-1}}(x)$ with a trained $s^{\theta_{i-1}}$ is equivalent to learn the neural network function $\mathrm{NN}^{\theta_i}(x)$}
Accordingly, $f^{\theta_i}(x)$, $1\le i\le n$, is defined via \eqref{def:ftheta} too. This piece-wise approximation in time domain introduces a time-discretization error which is proportional to $\Delta t$ and can be controlled with a small time step size. The discrete-time velocity fields $(f^{\theta_i})$, induce the transported probability density, denoted by $\tilde{p}_{t_i}$, as the approximation to the true probability density function $\widehat{p}_{t}(x)$ in the interval $(t_{i-1},t_{i})$. 
Specifically, set $\tilde{p}_{t_0}(x)$ as the initial distribution $\mu_{0}(x)$ in \eqref{sdeBL}, and $\tilde{p}_{t_i}$, $i=1,2,\ldots$, as the solution $q_{\Delta t}$ of the following transport equation driven by $f^{\theta_i}$ in a short time interval $\Delta t$ with the initial $\tilde{p}_{t_{i-1}}$: {\footnotesize\begin{equation}\label{eq:qf}\frac{\partial q_{t}(x)}{\partial t}=-\nabla \cdot\left(f^{\theta_{i}}(x) q_{t}(x)\right), \quad q_{0}(x)= \tilde{p}_{t_{i-1}}(x).\end{equation} }We now consider the following first order approximation in time to \eqref{opt:origin}: {\small \begin{equation}\label{eq:lossEuler}\begin{aligned}&\  \int_{t_0}^{T}\int_{\Omega}\left| f^{\theta}(x, t) - \mathcal{V}[p^{\theta}_{t}](x, t) \right|^{2}p^{\theta}_{t}(x) \,\mathrm{d}x \,\mathrm{d}t \\ \approx&\sum_{i=1}^{n}\int_{\Omega}\left| f^{\theta_{i}}(x) - \mathcal{V}[\tilde{p}_{t_{i-1}}](x, t_{i-1}) \right|^{2}\tilde{p}_{t_{i-1}}(x)\,\mathrm{d}x \,\Delta t.\end{aligned}\end{equation}}Note that $\tilde{p}_{t_{i-1}}$ only depends on the past velocity fields up to $f^{\theta_{i-1}}$, which allow us to apply the following iterative scheme: for $i=1,2,\ldots$, we sequentially solve \begin{equation}\label{eq:seqloss}\theta_{i}^{\star}:=\underset{{\theta_{i}}}{\operatorname{argmin}}~ \int_{\Omega}\left| f^{\theta_{i}}(x) - \mathcal{V}[\tilde{p}_{t_{i-1}}](x, t_{i-1}) \right|^{2}\tilde{p}_{t_{i-1}}(x)\,\mathrm{d}x.\end{equation} with the initial $\tilde{p}_{t_0}=\mu_{0}$, and for $i\ge 2$, $\tilde{p}_{t_{i-1}}$ is the image density of $\tilde{p}_{t_{i-2}}$ transported by the velocity field $f^{\theta^\star_{i-1}}$ within the short time interval $[t_{i-2}, t_{i-1}]$. To solve \eqref{eq:seqloss} for a given $i$, note that $\tilde{p}_{t_{i-1}}$ is known as samples from the transport maps obtained up to $f^{\theta_{i-1}}$, then we can derive the similar objective function analogous to Hyvärinen score-matching:

{\small \begin{equation*}\begin{aligned}&\int_{\Omega}\left| f^{\theta_i}(x,t_{i}) - \mathcal{V}[\tilde{p}_{t_{i-1}}](x,t_{i})\right|^{2}\tilde{p}_{t_{i-1}}\d x \\ =&\int_{\Omega}\left| s^{\theta_i}(x) + \left(K*\tilde{p}_{t_{i-1}}\right)(x) - \frac{1}{2}\Sigma(x,t_{i})\nabla\log \tilde{p}_{t_{i-1}}+ \mathcal{I}^{F,G}_{t_i}[\tilde{p}_{t_{i-1}}](x)\right|^{2}\tilde{p}_{t_{i-1}}\d x \\ =&\int_{\Omega}\left| s^{\theta_i}(x)\right|^{2}\tilde{p}_{t_{i-1}}\d x + \int_{\Omega}\left| \left(K*\tilde{p}_{t_{i-1}}\right)(x) - \frac{1}{2}\Sigma(x,t_{i})\nabla\log \tilde{p}_{t_{i-1}}+ \mathcal{I}_{t_i}^{F,G}[\tilde{p}_{t_{i-1}}](x)\right|^{2}\\&\times \tilde{p}_{t_{i-1}}\d x + 2\int_{\Omega}s^{\theta_i}(x)\cdot\left[\left( K*\tilde{p}_{t_{i-1}}\right)(x) - \frac{1}{2}\Sigma(x,t_{i})\nabla\log \tilde{p}_{t_{i-1}}+ \mathcal{I}_{t_i}^{F,G}[\tilde{p}_{t_{i-1}}](x)\right]\tilde{p}_{t_{i-1}}\d x.\end{aligned}\end{equation*}}

 Note that the second term contains no  $s^{\theta_i}$, so we can drop it during this optimization sub-problem at iteration $i$. For the remaining terms, we express them by the expectation w.r.t. $\tilde{p}_{t_{i-1}}$. That is, $\int_{\Omega}\left| s^{\theta_i}(x )\right|^{2}\tilde{p}_{t_{i-1}}(x) \d x = \mathbb{E}_{X_{t_i}\sim \tilde{p}_{t_{i-1}}}\left| s^{\theta_i}(X_{t_i}) \right |^{2}$, and $\int_{\Omega}s^{\theta_i}(x ) \cdot \left( K*\tilde{p}_{t_{i-1}}\right)(x) ~\tilde{p}_{t_{i-1}}(x)\d x = \mathbb{E}_{X_{t_i}\sim \tilde{p}_{t_{i-1}}}\mathbb{E}_{\widetilde {X}_{t_i}\sim \tilde{p}_{t_{i-1}}}\left(s^{\theta_i}(X_{t_i})\cdot K(X_{t_i}-\widetilde{X}_{t_i}) \right)$. And by integral by part, {\small \begin{equation*}\begin{split}&\int_{\Omega}s^{\theta_i}(x )\cdot \left[ - \frac{1}{2}\Sigma(x,t_{i})\nabla\log \tilde{p}_{t_{i-1}}(x) \right]~\tilde{p}_{t_{i-1}}(x) \d x \\&=\ -\frac{1}{2}\int_{\Omega}\left[\Sigma(x,t_{i})s^{\theta_i}(x )\right] \cdot\nabla \tilde{p}_{t_{i-1}}(x)\d x \\&= \frac{1}{2}\mathbb{E}_{X_{t_i}\sim \tilde{p}_{t_{i-1}}}\left[\nabla\cdot\left(\Sigma(X_{t_i},t_{i}) )s^{\theta_i}(X_{t_i})\right)\right].\end{split}\end{equation*} } And for the non-local terms by \eqref{notaionsbIA} we have, {\small\begin{equation*}\begin{aligned}&\int_{\Omega}s^{\theta_i}(x ) \mathcal{I}^{F,G}_{t_i}[\tilde{p}_{t_{i-1}}](x) ~\tilde{p}_{t_{i-1}}(x) \d x = \int_{\Omega}s^{\theta_i}(x )\cdot\int_{|r|<1}F(r,t_{i}) \int_{0}^{1}\tilde{p}_{t_{i-1}}(x-\lambda F(r,t_{i}))\\&\times\d \lambda\,\nu(\d r) \d x + \int_{\Omega}s^{\theta_i}(x )\cdot\int_{|r|\geq1}G(r,t_{i}) \int_{0}^{1}\tilde{p}_{t_{i-1}}(x-\lambda G(r,t_{i}))\d \lambda\,\nu(\d r) \d x\\ =&\mathbb{E}_{X_{t_i}\sim \tilde{p}_{t_{i-1}}}\left(\int_{|r|<1}\int_{0}^{1}\left(s^{\theta_i}(X_{t_i}+\lambda F(r,t_{i}) )\cdot F(r,t_{i})\right)\d \lambda\,\nu(\d r) \right)\\&+ \mathbb{E}_{X_{t_i}\sim \tilde{p}_{t_{i-1}}}\left( \int_{|r|\geq1}\int_{0}^{1}\left(s^{\theta_i}(X_{t_i}+\lambda G(r,t_{i}) )\cdot G(r,t_{i})\right) \d\lambda\,\nu(\d r)\right).\end{aligned}\end{equation*}}

In summary, for each $i=1,2,\ldots, n$ sequentially, to minimize \eqref{eq:seqloss} is equivalent to minimizing the following loss function of $s^{\theta_i}(\cdot )$ {\small\begin{equation}\label{eqn:losst}\begin{aligned}&Loss(t_{i}) = \mathbb{E}_{X_{t_i}\sim \tilde{p}_{t_{i-1}}}\left(\left|s^{\theta_i}(X_{t_i}) \right|^{2}\right) + \mathbb{E}_{X_{t_i}\sim \tilde{p}_{t_{i-1}}}\left[\nabla\cdot\left(\Sigma(X_{t_i},t_{i})s^{\theta_i}(X_{t_i})\right)\right] +2\mathbb{E}_{X_{t_i}\sim \tilde{p}_{t_{i-1}}}\bigg[ \\&\mathbb{E}_{X'_{t_i}\sim \tilde{p}_{t_{i-1}}}\left(s^{\theta}(X_{t_i},t_{i}) \cdot K(X_{t_i}-X'_{t_i}) \bigg) \bigg] + 2 \mathbb{E}_{X_{t_i}\sim \tilde{p}_{t_{i-1}}}\left(\int_{|r|<1}\int_{0}^{1}s^{\theta_i}(X_{t_i}+\lambda F(r,t_{i}))\cdot F(r,t_{i}) \right. \right.\\&\times\d \lambda\nu(\d r)\bigg) + 2\mathbb{E}_{X_{t_i}\sim \tilde{p}_{t_{i-1}}}\left( \int_{|r|\geq1}\int_{0}^{1}s^{\theta_i}(X_{t_i}+\lambda G(r,t_{i}))\cdot G(r,t_{i}) \d\lambda\nu(\d r)\right).\end{aligned}\end{equation}  }

Given the current samples of $X_{t_i}$ from law
 $\tilde{p}_{t_{i-1}}$ at any time $t_{i}$, we can
obtain $s^{\theta_i}(\cdot)$ via direct minimization of objective
in \eqref{eqn:losst}. Given $s^{\theta_i}(\cdot)$, we propagate
 $X^{\theta}_{t_0,t_i}$ forward in time up to time
 $t_{i}+\Delta t$ via \eqref{eqn:ODEtheta}. The resulting
procedure, which alternates between self-consistent score estimation and sample propagation,
is presented in Algorithm \ref{algorithm} for the choice of a forward-Euler
integration routine as  $t_{i}$ to
 $t_{i}+\Delta t$.

\begin{algorithm}
    \caption{Sequential  mean-field L\'{e}vy score-based transport
    modeling }
    \label{algorithm}
    \begin{algorithmic}
        \STATE{Input: An initial time $t_{0}$. A set of $N$ samples $\{x^{(i)}\}_{i=1}^{N}$ from the initial distribution $p_{t_0}$. A time step $\Delta t$ and the number of steps $N_{T}= [(T-t_{0})/\Delta t]$. Initialize sample locations $X^{(i)}_{t_0}=x^{(i)}$ for $i=1,\cdots,N$.}
        \FOR{$k=0:N_{T}-1$} \STATE{Optimize \eqref{eqn:losst}:} \STATE{{\small$s^{\theta_k}(\cdot)=\underset{s^{\theta_k}(\cdot)}{\operatorname{argmin}}~ \frac{1}{N}\sum_{i=1}^{N}\left[ \left|s^{\theta_k}(X_{t_k}^{(i)})\right|^{2}+ 2\left(\frac{1}{N}\sum_{j=1}^{N}s^{\theta_k}(X_{t_k}^{(i)})K(X_{t_k}^{(i)}-X_{t_k}^{(j)}) \right) \right.$\  $+ \nabla\cdot \left( \Sigma(X_{t_k}^{(i)},t_{k})s^{\theta_k}(X_{t_k}^{(i)}) \right) + 2\int_{|r|<1}\int_{0}^{1}s^{\theta_k}(X_{t_k}^{(i)}+\lambda F(r,t_{k}) )\cdot F(r,t_{k})\d\lambda\nu(\d r)$\
$\left. + 2 \int_{|r|\geq 1}\int_{0}^{1}s^{\theta_k}(X_{t_k}^{(i)}+ \lambda G(r,t_{k}) )\cdot G(r,t_{k})\d \lambda\nu(\d r) \right]$;}}
        \STATE{Propagate the samples for $i=1,\cdots,N$: \\ $X^{(i)}_{t_{k+1}}=X^{(i)}_{t_{k}}+ \Delta t \big[ b(X^{(i)}_{t_{k}},t_{k}) - \int_{|r|<1}F(r,t)\nu(\d r) -\frac{1}{2}\nabla\cdot\Sigma(X^{(i)}_{t_{k}},t_{k}) - s^{\theta_k}(X^{(i)}_{t_{k}}) \big];$}
        \STATE{Set $t_{k+1}=t_{k}+\Delta t$ ; } \ENDFOR \RETURN \STATE{Output: Samples $\{X^{(i)}_{t_k}\}_{i=1}^{N}$ from $p_{t_k}$ and $\{s^{\theta_k}(X^{(i)}_{t_k})\}_{i=1}^{N}$ for all $\{t_{k}\}_{k=0}^{N_T}$.}
    \end{algorithmic}
\end{algorithm}
Here, the one-dimensional integral of $\int_{0}^{1}\bullet \d \lambda$ is discretized
by any quadrature scheme, like the trapezoidal rule. The integral w.r.t.
$\nu(\d r )$ in $\mathbb{R}^{d}$ is approximated by the quadrature scheme in low
dimension or by the Monte Carlo average in high dimension. $\int_{|r|<1}F(r,t)\nu
(\d r)$ is pre-computed as a function of $t$ only.

In parallel, Algorithm \ref{algorithm2} is based on the interaction flow \eqref{eqn:IODEX*} and adopts the approach from \cite{lu2024score}. This strategy learns the score function while excluding the interaction term, thereby delegating the computation of the mean-field interaction entirely to the transport step. 

{ { \begin{algorithm}\caption{Sequential L\'{e}vy score-based mean-field transport modeling } \label{algorithm2} \begin{algorithmic}\STATE{Input: An initial time $t_{0}$. A set of $N$ samples $\{x^{(i)}\}_{i=1}^{N}$ from the initial distribution $p_{t_0}$. A time step $\Delta t$ and the number of steps $N_{T}$. Initialize sample locations $X^{(i)}_{t_0}=x^{(i)}$ for $i=1,\cdots,N$.} \FOR{$k=0:N_{T}-1$} \STATE{Optimize:} \STATE{{\small$s^{\theta_k}(\cdot )=\underset{s^{\theta_k}(\cdot )}{\operatorname{argmin}}~ \frac{1}{N}\sum_{i=1}^{N}\left[ \left|s^{\theta_k}(X_{t_k}^{(i)})\right|^{2}+ \nabla\cdot \left( \Sigma(X_{t_k}^{(i)},t_{k})s^{\theta_k}(X_{t_k}^{(i)}) \right) \right.$\  $+2\int_{|r|<1}\int_{0}^{1}s^{\theta_k}(X_{t_k}^{(i)}+\lambda F(r,t_{k}) )\cdot F(r,t_{k})\d\lambda\nu(\d r) \qquad \qquad \quad \quad \quad \quad\quad \quad \quad\quad \quad \quad\quad$ $\left. + 2 \int_{|r|\geq 1}\int_{0}^{1}s^{\theta_k}(X_{t_k}^{(i)}+ \lambda G(r,t_{k}) )\cdot G(r,t_{k})\d \lambda\nu(\d r) \right]$;}} \STATE{Propagate the samples for $i=1,\cdots,N$: \\ $X^{(i)}_{t_{k+1}}=X^{(i)}_{t_{k}}+ \Delta t \big[ b(X^{(i)}_{t_{k}},t_{k}) + \frac{1}{N}\sum_{j=1}^{N}K(X^{(i)}_{t_k}- X^{(j)}_{t_k}) - \int_{|r|<1}F(r,t)\nu(\d r) -\frac{1}{2}\nabla\cdot\Sigma(X^{(i)}_{t_{k}},t_{k}) - s^{\theta_k}(X^{(i)}_{t_{k}}) \big];$} \STATE{Set $t_{k+1}=t_{k}+\Delta t$ ; } \ENDFOR \RETURN \STATE{Output: Samples $\{X^{(i)}_{t_k}\}_{i=1}^{N}$ from $p_{t_k}$ and $\{s^{\theta_k}(X^{(i)}_{t_k})\}_{i=1}^{N}$ for all $\{t_{k}\}_{k=0}^{N_T}$.}\end{algorithmic}\end{algorithm} } }

\begin{remark}
    \label{remark:algorithm} The difference between Algorithm \ref{algorithm} and Algorithm \ref{algorithm2} lies in the treatment of the interaction term, as shown in their respective names. Algorithm~\ref{algorithm} requires higher computational resources during the training phase due to the double sum $\sum_{i}\sum_{j}$ of interaction in the score computation, which has a complexity of $O(N^{2})$, while Algorithm~\ref{algorithm2} avoids this by excluding the interaction term in the score computation, resulting in a complexity of $O(N)$. However, during the propagation step, Algorithm~\ref{algorithm2} still requires the double sum $\sum_{i}\sum_{j}$ to compute the interaction term in \eqref{eqn:ODEX*}, leading to a complexity of $O(N^{2})$. In contrast, Algorithm~\ref{algorithm} incorporates the interaction term directly into the learned velocity field, allowing the propagation step to be performed independently for each particle as in \eqref{eqn:IODEX*}, with a complexity of $O(N)$. But in practice, the training phase is often more time-consuming than the propagation phase since the training involves multiple iterations to optimize the neural network parameters, while propagation requires only a single pass through the data. In this sense, Algorithm~\ref{algorithm2} is recommended in practice for challenging problems with a large number of particles. 
\end{remark}
\bigskip
If the initial density $\widehat{p}_{0}$ has the analytical expression like the Gaussian
distribution, then the initial score function $s^{\theta_0}(\cdot )$ is computed
by the direct minimization of the following loss due to \eqref{eqn:stheta}:
{\small\[\frac{\int_{\Omega}\left\| s^{\theta_0}(\cdot ) - \widehat{S}[\widehat{p}_{0}](x,t_{0}) + K* \widehat{p}_{0}(x)\right\|^{2}\widehat{p}_{0}(x) \, dx}{\int_{\Omega}\left\| \widehat{S}[\widehat{p}_{0}](x,t_{0}) - K*\widehat{p}_{0}(x)\right\|^{2}\widehat{p}_{0}(x) \, dx}\approx \frac{\sum_{i=1}^{N}\left\| s^{\theta_0}(X_{0}^{(i)}) - \widehat{S}[\widehat{p}_{0}](X_{0}^{(i)},t_{0}) + K* \widehat{p}_{0}(x)\right\|^{2}}{\sum_{i=1}^{N}\left\| \widehat{S}[\widehat{p}_{0}](X_{0}^{(i)},t_{0}) - K* \widehat{p}_{0}(x)\right\|^{2}},\] }
where $\widehat{S}[\widehat{p}_{0}]$ is defined in \eqref{eqn:GS}.

\subsection{ Error analysis of Algorithm \ref{algorithm}}
\label{sec:error} Now we focus on analyzing the discrete error associated with Algorithm
\ref{algorithm}. The following theorem demonstrates that the discrete numerical
error of Algorithm \ref{algorithm} is on the order of $1/\sqrt{N}$, where $N$ represents
the number of particles utilized in the algorithm.
\begin{theorem}
    \label{theo:error} Let $N$, $N_{T}$, $\Delta t$, $t_{0}$ and $p_{t_0}$ be
    defined as in Algorithm \ref{algorithm} and denote
    $\varepsilon:=\sup_{x\in\Omega,t_0\leq t\leq N_T\Delta t}\left|s(x,t) + (K*\widehat
    {p}_{t})(x) - \frac{\Sigma(x,t)}{2}\nabla\log \widehat{p}_{t}(x) + \mathcal{I}
    _{t}^{F,G}[\widehat{p}_{t}](x)\right|$, as the approximation error of \eqref{opt:origin}
    by a neural network function $s(x,t)$, with $\widehat{p}_{t}$ the solution of
    the L\'{e}vy--Fokker--Planck equation \eqref{eqn:FPE}. Denote
    $X^{N}_{t_0,t}(x)$ the numerical transport map obtained from Algorithm
    \ref{algorithm} by stacking the transport map in each sub-interval starting
    from $x$. Let $X^{*}_{t_0,t}(x)$ starting from $x$ solve the deterministic probability
    flow equation $\frac{\d}{\d t}X^{*}_{t_0,t}(x)=\ \mathcal{V}[\widehat{p}_{t}]
    (X^{*}_{t_0,t}(x),t) = \widehat{b}(X^{*}_{t_0,t}(x),t) + (K*\widehat{p}_{t})(
    X^{*}_{t_0,t}( x)) \ -\frac{1}{2}\Sigma(X^{*}_{t_0,t}(x),t)\nabla\log \widehat
    {p}_{t}(X^{*}_{t_0,t}(x)) + \mathcal{I}_{t}^{F,G}[\widehat{p}_{t}](X^{*}_{t_0,t}
    (x),t)$, corresponding to ODE \eqref{eqn:ODEX*}, $\widehat{b}$ and
    $\mathcal{I}_{t}^{F,G}[\widehat{p}_{t}]$ are given in \eqref{notaionsbIA},
    \eqref{notaionsIFG}. Under Assumption \ref{ass:densitytheta}, for any $t\in\{
    t_{0}+n\Delta t\}_{n=0}^{N_T}$,
    {\small\begin{equation}\label{numerror}\begin{aligned}\mathbb{E}_{x\sim p_{t_0}}\left|X^{*}_{t_0,t}(x)-X^{N}_{t_0,t}(x)\right|\leq \left[ \mathcal{O}(\varepsilon)+\mathcal{O}(\Delta t)\right](t-t_{0}),\end{aligned}\end{equation} }as
    $\varepsilon$ and $\Delta t$ all tend to 0.
\end{theorem}

\begin{remark}
 Algorithm \ref{algorithm} and \ref{algorithm2} in fact adopt the sequential training approach as explained in Section \ref{sec:algorithm}. But if the first order approximation of the forward Euler scheme in \eqref{eq:lossEuler} holds uniformly for all possible functions $s^{\theta}$, $\theta\in\Theta$, then, when our sequential minimization algorithms make the sequential loss \eqref{eq:seqloss} used in practice smaller than $\varepsilon$, the continuous loss \eqref{opt:origin} is bounded by $\varepsilon+C_{t}\Delta t$ with a constant $C_{t}$ which may depend on $t$. The error estimate \eqref{numerror} is still valid by absorbing the extra term $C_{t}\Delta t$ into $\mathcal{O}(\Delta t)$.    
\end{remark} 

\subsection{Proof of Theorem \ref{theo:KL}}

We first present a lemma 
about the  bounds of  the squared $L^2$-norm of two probability density functions  by their KL divergence, which will be used in our main proof.
\begin{lemma}\label{L2KL}
    Suppose $p$ and $q$ are two probability densities on $\Omega$, and there exists a positive constant $\tau$ such that
  $0<p(x),q(x)<\tau,\ \forall x\in\Omega$. 
    Then we have
    {\footnotesize\begin{equation*}
        \int_{\Omega}|p(x)-q(x)|^2\d x\leq \frac{2\tau}{1-\log2}D_{\mathrm{KL}}(p\|q).
    \end{equation*}}
\end{lemma}
\begin{proof}
Define $\zeta(x):=\ \frac{q(x)-p(x)}{p(x)},\ \forall x\in\Omega$. Then $D_{\mathrm{KL}}(p\|q)=\ \int_{\Omega}p(x)\log\frac{p(x)}{q(x)}\d x=\ -\int_{\Omega}p(x)\log(1+\zeta(x))\d x$. Define two Borel sets as follows, $A:= \{x | \zeta(x)>1\}=\ \{x | q(x)>2p(x)\},\ B:= \{x | \zeta(x)\leq 1\}=\ \{x | q(x)\leq 2p(x)\}$. We obtain that for $x\in A$, $1+\zeta(x)\leq e^{\alpha\zeta(x)}$ where $\alpha=\log2$; for $x\in B$, $1+\zeta(x)\leq e^{\zeta(x)-\beta\zeta(x)^2}$ where $\beta=1-\log2$. Note that $p$ and $q$ are two probabilities density on $\Omega$, that is $\int_{\Omega}p(x)\zeta(x)\d x=\ \int_{\Omega}(q(x)-p(x))\d x=0$, which implies that $\int_{A}p(x)\zeta(x)\d x=\ - \int_{B}p(x)\zeta(x)\d x$. Thus 
    {\footnotesize\begin{equation*}
        \begin{aligned}
            D_{\mathrm{KL}}(p\|q) =&\ -\int_{A}p(x)\log(1+\zeta(x))\d x  -\int_{B}p(x)\log(1+\zeta(x))\d x\\
            \geq&\ -\alpha\int_{A}p(x)\zeta(x)\d x -\int_{B}p(x)\zeta(x)\d x + \beta\int_{B}p(x)\zeta(x)^2\d x\\
            =&\ (1-\alpha)\int_{A}p(x)\zeta(x)\d x + \beta\int_{\Omega}p(x)\zeta(x)^2\d x\\
            =&\ (1-\log2)\left( \int_{A}|q(x)-p(x)|\d x + \int_{B}p(x) \left(\frac{q(x)-p(x)}{p(x)}\right)^2\d x  \right).
        \end{aligned}
    \end{equation*}
    }For the first term in right hand side of the above equality, we have $\int_{A}|q(x)-p(x)|\d x\geq\ \frac{1}{2\tau}\int_{A}|q(x)-p(x)|^2\d x$. For the second term, we have $\int_{B}p(x) \left(\frac{q(x)-p(x)}{p(x)}\right)^2\d x \geq\ \frac{1}{2\tau}\int_{B}  \left| q(x)-p(x) \right|^2\d x$. Finally, we have $D_{\mathrm{KL}}(p\|q)  
            \geq\ \frac{1-\log2}{2\tau}\int_{\Omega}|q(x)-p(x)|^2\d x$.
\end{proof}

\medskip 
\begin{proof}[Proof of Theorem \ref{theo:KL}]

   The proof is inspired by the methodologies used in \cite[Proposition 1]{boffi2023probability} and \cite[Appendix E]{shen2024entropy}. However, dealing with the L\'{e}vy term introduces substantial complexity, requiring the introduction of novel techniques for a thorough analysis.
    
    Firstly, according to the definition of KL divergence, we derive that
    {\footnotesize\begin{equation*}
        \begin{aligned}
            &\ \frac{\d}{\d t}D_{\mathrm{KL}}(p^{\theta}_t\|\widehat{p}_t)=\ \frac{\d}{\d t}\int_{\Omega}\log \left(\frac{p^{\theta}_t(x)}{\widehat{p}_t(x)}\right)p^{\theta}_t(x)\d x\\
            =&\ -\int_{\Omega}\frac{p^{\theta}_t(x)}{\widehat{p}_t(x)}\partial_t\widehat{p}_t(x)\d x + \int_{\Omega}\log\left(\frac{p^{\theta}_t(x)}{\widehat{p}_t(x)}\right)\partial_tp^{\theta}_t(x)\d x + \int_{\Omega}\partial_tp^{\theta}_t(x) \d x\\
            =&\ -\int_{\Omega}\frac{p^{\theta}_t(x)}{\widehat{p}_t(x)}\partial_t\widehat{p}_t(x)\d x + \int_{\Omega}\log\left(\frac{p^{\theta}_t(x)}{\widehat{p}_t(x)}\right)\partial_tp^{\theta}_t(x)\d x\\
            =&\ -\int_{\Omega}\mathcal{V}[\widehat{p}_t](x,t)\cdot\nabla\left(\frac{p^{\theta}_t(x)}{\widehat{p}_t(x)}\right)\widehat{p}_t(x)\d x + \int_{\Omega}f^\theta(x,t)\cdot\nabla\log\left(\frac{p^{\theta}_t(x)}{\widehat{p}_t(x)}\right) p^{\theta}_t(x)\d x\\
            =&\ -\int_{\Omega}\left(\mathcal{V}[\widehat{p}_t](x,t)-f^\theta(x,t) + \mathcal{V}[p^\theta_t](x,t) - \mathcal{V}[p^\theta_t](x,t)\right)\cdot \nabla\log \left(\frac{p^{\theta}_t(x)}{ \widehat{p}_t(x) }\right)p^{\theta}_t(x)\d x,
        \end{aligned}
    \end{equation*}
}where we used  $\int\partial_t p^{\theta}_t(x)\d x =0$ and the integration by parts. 

With the expression $\mathcal{V}[p^\theta_t]$ of \eqref{notaionsbIA}, we  decompose  the derivative of the KL divergence  as follows,
 {\footnotesize \begin{equation*}
        \begin{aligned}        
            &\ \frac{\d}{\d t}D_{\mathrm{KL}}(p^\theta_t\|\widehat{p}_t)=\ \underbrace{\int_{\Omega}\left( f^\theta(x,t) - \mathcal{V}[p^\theta_t](x,t)\right)\cdot \nabla\log \left(\frac{p^{\theta}_t(x)}{ \widehat{p}_t(x) }\right)p^{\theta}_t(x)\d x}_{\mbox{Perturbation}}\\
            &\ + \underbrace{\int_{\Omega}\left(\int_{\Omega}K(x-y)(p^{\theta}_t(y)-\widehat{p}_t(y))\d y\right)\cdot\nabla\log \left(\frac{p^{\theta}_t(x)}{ \widehat{p}_t(x) }\right)p^{\theta}_t(x)\d x }_{\mbox{Interaction}} \underbrace{-\int_{\Omega}\frac{1}{2}\left|\nabla\log \left(\frac{p^{\theta}_t(x)}{ \widehat{p}_t(x) }\right)\right|^2_{\Sigma(x,t)}p^{\theta}_t(x)\d x}_{\mbox{Diffusion}}\\
            &\ + \underbrace{\int_{\Omega}\left(\mathcal{I}_t(p^\theta_t,F,G)(x)-\mathcal{I}_t(\widehat{p}_t,F,G)(x) \right)\cdot\nabla\log \left(\frac{p^{\theta}_t(x)}{ \widehat{p}_t(x) }\right)p^\theta_t(x)\d x}_{\mbox{L\'{e}vy}},
        \end{aligned}
    \end{equation*}
    }where $|\cdot|^2_{\Sigma(x,t)}=\langle\cdot,\Sigma(x,t)\cdot\rangle$.The bounds for these four terms are studied below respectively in four steps.
    
{\bf (Step 1.)} For the perturbation part, by the Cauchy–Schwarz inequality we have
    {\footnotesize\begin{equation*}
        \begin{aligned}
            &\ \int_{\Omega}\left( f^\theta(x,t) - \mathcal{V}[p^\theta_t](x,t)\right)\cdot \nabla\log \left(\frac{p^{\theta}_t(x)}{ \widehat{p}_t(x) }\right)p^{\theta}_t(x)\d x\\
            \leq&\ \frac{1}{12C_\sigma}\int_{\Omega}\left|\nabla\log \left(\frac{p^{\theta}_t(x)}{ \widehat{p}_t(x) }\right)\right|^2 p^\theta_t(x)\d x + 3 C_\sigma\int_{\Omega}\left| f^\theta(x,t) - \mathcal{V}[p^\theta_t](x,t)\right|^2p^\theta_t(x)\d x.
        \end{aligned}
    \end{equation*}
}

{\bf (Step 2.)} For the interaction part, there are two cases \ref{subass:bounded} and \ref{subass:BS} for the interaction kernel $K$. We discuss each case respectively.

{\bf (Bounded interaction)} By the Cauchy–Schwarz inequality we have
 {\footnotesize \begin{equation*}
        \begin{aligned}
            &\ \int_{\Omega}\left(\int_{\Omega}K(x-y)(p^{\theta}_t(y)-\widehat{p}_t(y))\d y\right)\cdot\nabla\log \left(\frac{p^{\theta}_t(x)}{ \widehat{p}_t(x) }\right)p^{\theta}_t(x)\d x\\
            \leq&\ \frac{1}{12C_\sigma}\int_{\Omega}\left|\nabla\log \left(\frac{p^{\theta}_t(x)}{ \widehat{p}_t(x) }\right)\right|^2 p^\theta_t(x)\d x+ 3C_\sigma\int_{\Omega}\left| (K*(p^{\theta}_t-\widehat{p}_t))(x)\right|^2p^\theta_t(x)\d x \\
            \leq&\ \frac{1}{12C_\sigma}\int_{\Omega}\left|\nabla\log \left(\frac{p^{\theta}_t(x)}{ \widehat{p}_t(x) }\right)\right|^2 p^\theta_t(x)\d x+ 3C_\sigma (C_K)^2 \left(\int_{\Omega}|p^{\theta}_t(y)-\widehat{p}_t(y)|\d y\right)^2\\
            \leq&\ \frac{1}{12C_\sigma}\int_{\Omega}\left|\nabla\log \left(\frac{p^{\theta}_t(x)}{ \widehat{p}_t(x) }\right)\right|^2 p^\theta_t(x)\d x+ 6C_\sigma (C_K)^2 D_{\mathrm{KL}}(p^\theta_t\| \widehat{p}_t),
        \end{aligned}
    \end{equation*}
}where we use the Csisz\'{a}r–Kullback–Pinsker inequality \cite{villani2009optimal} for the last inequality.

{\bf (Biot--Savart interaction)} 
For the Biot--Savart interaction kernel, as noted by \cite{jabin2018quantitative}, it can be written as:
{\footnotesize $$
        K^\top=\nabla \cdot U \quad\mbox{with}\quad U(x)=\frac{1}{2\pi} \left(
\begin{matrix}
  -\arctan(\frac{x_1}{x_2}) & 0 \\
  0 & \arctan(\frac{x_2}{x_1})
\end{matrix} \right).$$
}
We follow the statement in \cite[Appendix E.1]{shen2024entropy} to have 
{\footnotesize\begin{equation*}
    \begin{aligned}
        &\ \int_{\Omega}\left(\int_{\Omega}K(x-y)(p^{\theta}_t(y)-\widehat{p}_t(y))\d y\right)\cdot\nabla\log \left(\frac{p^{\theta}_t(x)}{ \widehat{p}_t(x) }\right)p^{\theta}_t(x)\d x\\
        \leq&\ \frac{1}{12C_\sigma}\int_{\Omega}\left|\nabla\log \left(\frac{p^{\theta}_t(x)}{ \widehat{p}_t(x) }\right)\right|^2 p^\theta_t(x)\d x+  3C_\sigma\int_{\Omega}\left|  (\nabla\log\widehat{p}_t(x))^\top U*(p^\theta_t -\widehat{p}_t )(x) \right|^2p^\theta_t(x)\d x\\
        &\ + 4\|U\|_\infty\left\|\frac{\nabla^2\widehat{p}_t}{\widehat{p}_t}\right\|_\infty D_{\mathrm{KL}}(p^\theta_t\|\widehat{p}_t)\\
         \leq&\ \frac{1}{12C_\sigma}\int_{\Omega}\left|\nabla\log \left(\frac{p^{\theta}_t(x)}{ \widehat{p}_t(x) }\right)\right|^2 p^\theta_t(x)\d x+  3C_\sigma  \left\| \nabla\log\widehat{p}_t \right\|^2_\infty \left\|U\right\|^2_\infty \left(\int_{\Omega}\left|p^\theta_t -\widehat{p}_t )(x) \right| \d x\right)^2\\
        &\ +  4\|U\|_\infty\left\|\frac{\nabla^2\widehat{p}_t}{\widehat{p}_t}\right\|_\infty D_{\mathrm{KL}}(p^\theta_t\|\widehat{p}_t)\\
        \leq&\  \int_{\Omega}\left|\nabla\log \left(\frac{p^{\theta}_t(x)}{ \widehat{p}_t(x) }\right)\right|^2 \frac{p^\theta_t(x)}{12C_\sigma} \d x   + \left( 6C_\sigma  \left\| \nabla\log\widehat{p}_t \right\|^2_\infty \left\|U\right\|^2_\infty + 4\|U\|_\infty\left\|\frac{\nabla^2\widehat{p}_t}{\widehat{p}_t}\right\|_\infty\right) D_{\mathrm{KL}}(p^\theta_t\|\widehat{p}_t),
    \end{aligned}
\end{equation*}
}we use the Csisz\'{a}r–Kullback–Pinsker inequality to obtain the last inequality.

{\bf (Step 3.)} For the Diffusion part, due to Assumption \ref{ass:allcoefficient}-\ref{subass:boundsigmaK} we simply have that
{\footnotesize  \begin{equation*}
        \begin{aligned}
            -\int_{\Omega}\frac{1}{2}\left|\nabla\log \left(\frac{p^{\theta}_t(x)}{ \widehat{p}_t(x) }\right)\right|^2_{\Sigma(x,t)}p^{\theta}_t(x)\d x\leq -\int_{\Omega}\frac{1}{2C_\sigma}\left|\nabla\log \left(\frac{p^{\theta}_t(x)}{ \widehat{p}_t(x) }\right)\right|^2 p^\theta_t(x)\d x.
        \end{aligned}
    \end{equation*}}

{\bf (Step 4.)}  The last part is about the  L\'{e}vy terms $\mathcal{I}$.  We have the following estimate:
{\footnotesize \begin{equation*}
        \begin{aligned}
            &\ \int_{\Omega}\left(\mathcal{I}^F_t[p^\theta_t](x)-\mathcal{I}^{F}_t[\widehat{p}_t](x) \right)\cdot\nabla\log \left(\frac{p^{\theta}_t(x)}{ \widehat{p}_t(x) }\right)p^\theta_t(x)\d x \\
            \leq&  \int_{\Omega}\int_{|r|<1}\int_0^1   \left| \frac{F(r,t)}{p^\theta_t(x)\widehat{p}_t(x)} \right| \bigg|p^\theta_t(x - \lambda F(r,t))\widehat{p}_t(x)-p^\theta_t(x - \lambda F(r,t))p^\theta_t(x) \\
            &\ + p^\theta_t(x - \lambda F(r,t))p^\theta_t(x) - \widehat{p}(x - \lambda F(r,t))p^\theta_t(x)\bigg| \d \lambda\nu(\d r) \left| \nabla\log \left(\frac{p^{\theta}_t(x)}{ \widehat{p}_t(x) }\right)\right| p_t^\theta(x) \d x\\
            \leq&\   \int_{\Omega}\int_{|r|<1}\int_0^1   \left| \frac{F(r,t)p^\theta_t(x - \lambda F(r,t))}{p^\theta_t(x)\widehat{p}_t(x)} \right| \left| \widehat{p}_t(x)- p^\theta_t(x)\right| \d\lambda \nu(\d r) \left| \nabla\log \left(\frac{p^{\theta}_t(x)}{ \widehat{p}_t(x) }\right)\right| p_t^\theta(x) \d x \\
            &\  +  \int_{\Omega}\left( \int_{|r|<1} \int_0^1 \left| \frac{F(r,t)p^\theta_t(x)}{p^\theta_t(x)\widehat{p}_t(x)} \right| \left| p^\theta_t(x - \lambda F(r,t))  - \widehat{p}(x - \lambda F(r,t)) \right| \d \lambda\nu(\d r)\right) \left| \nabla\log \left(\frac{p^{\theta}_t(x)}{ \widehat{p}_t(x) }\right)\right| p_t^\theta(x) \d x  \\
            \overset{(1)}{\leq}&\   6C_\sigma\int_{\Omega} \left(\int_{|r|<1}\int_0^1   \left|\frac{F(r,t) p^\theta_t(x - \lambda F(r,t)) }{p^\theta_t(x)\widehat{p}_t(x) }\right| \left|\widehat{p}_t(x)- p^\theta_t(x)\right| \d\lambda \nu(\d r)\right)^2 p^\theta_t(x)\d x \\
            &\ + \frac{1}{12C_\sigma}\int_{\Omega}\left| \nabla\log \left(\frac{p^{\theta}_t(x)}{ \widehat{p}_t(x) }\right)\right|^2 p_t^\theta(x) \d x \\
            &\ +  6C_\sigma\int_{\Omega}\left(\int_{|r|<1} \int_0^1 \left|\frac{F(r,t)}{\widehat{p}_t(x)}\right| \left| p^\theta_t(x - \lambda F(r,t)) - \widehat{p}(x - \lambda F(r,t)) \right| \d \lambda\nu(\d r) \right)^2p_t^\theta(x) \d x \\
             \overset{(2)}{\leq}&\   6C_\sigma (C_T^f)^3(C_T^*)^2\int_{|r|<1}  \left| F(r,t) \right|^2\nu(\d r)  \int_{\Omega}\left|\widehat{p}_t(x)- p^\theta_t(x)\right|^2  \d x + \frac{1}{12C_\sigma}\int_{\Omega}\left| \nabla\log \left(\frac{p^{\theta}_t(x)}{ \widehat{p}_t(x) }\right)\right|^2 p_t^\theta(x) \d x \\
            &\ +  6C_\sigma C_T^f(C_T^*)^2\int_{|r|<1} \int_0^1 \left| F(r,t) \right|^2 \int_{\Omega} \left| p^\theta_t(x - \lambda F(r,t)) - \widehat{p}(x - \lambda F(r,t)) \right|^2 \d x  \d \lambda\nu(\d r)  \\
            =&\ 6C_\sigma \left( (C_T^f)^3(C_T^*)^2 + C_T^f(C_T^*)^2\right)\int_{|r|<1}  \left| F(r,t) \right|^2\nu(\d r)  \int_{\Omega}\left|\widehat{p}_t(x)- p^\theta_t(x)\right|^2  \d x + \frac{1}{12C_\sigma}\int_{\Omega}\left| \nabla\log \left(\frac{p^{\theta}_t(x)}{ \widehat{p}_t(x) }\right)\right|^2 p_t^\theta(x) \d x,
        \end{aligned}
    \end{equation*}}where the inequality ``$ \overset{(1)}{\leq}$'' arises from   the Cauchy–Schwarz inequality, and   ``$ \overset{(2)}{\leq}$'' comes from Assumption\ref{ass:densitytheta}-\ref{subass:truedensity}.    
    $C^f_T$ is the constant given in Remark  \ref{eq:CTf}. Similarly,  we have
{\footnotesize\begin{equation*}
    \begin{aligned}
        &\ \int_{\Omega}\left(\mathcal{I}_t^G[p^\theta_t](x)-\mathcal{I}_t^G[\widehat{p}_t](x)\right) \cdot\nabla\log \left(\frac{p^{\theta}_t(x)}{ \widehat{p}_t(x) }\right) p^\theta_t(x)\d x\\
        \leq &\ 6C_\sigma \left( (C_T^f)^3(C_T^*)^2 + C_T^f(C_T^*)^2\right)\int_{|r|\geq1}  \left| G(r,t) \right|^2\nu(\d r)  \int_{\Omega}\left|\widehat{p}_t(x)- p^\theta_t(x)\right|^2  \d x\\
            &\ + \frac{1}{12C_\sigma}\int_{\Omega}\left| \nabla\log \left(\frac{p^{\theta}_t(x)}{ \widehat{p}_t(x) }\right)\right|^2 p_t^\theta(x) \d x.
    \end{aligned}
\end{equation*}
}By Lemma \ref{L2KL}, the squared $L^2$-norm of $p^\theta_t-\widehat{p}_t$ can be controlled by the KL divergence between $p^\theta_t$ and $\widehat{p}_t$, thus we know that for any $t\in[t_0,T]$,
{\footnotesize   \begin{equation*}
        \begin{aligned}
            &\ \int_{\Omega}\left(\mathcal{I}^{F,G}_t[p^\theta_t](x)-\mathcal{I}_t^{F,G}[\widehat{p}_t](x) \right)\cdot\nabla\log \left(\frac{p^{\theta}_t(x)}{ \widehat{p}_t(x) }\right)p^\theta_t(x)\d x\\
            \leq &\ \frac{1}{6C_\sigma}\int_{\Omega}\left|\nabla\log \left(\frac{p^{\theta}_t(x)}{ \widehat{p}_t(x) }\right)\right|^2 p^\theta_t(x)\d x + 6C_\sigma C^{\mathrm{nG}}_TD_{\mathrm{KL}}(p^\theta_t\|\widehat{p}_t),
        \end{aligned}
    \end{equation*}
}where   
  {\footnotesize  \begin{equation}
        \begin{aligned}
            C^{\mathrm{nG}}_T=&\ \frac{2\left(C_T^f\wedge C_T^*\right)\left(C^f_TC^*_T\right)^2}{1-\log2}\left[ \int_{|r|<1}|F(r,t) |^2 \nu(\d r)+ \int_{|r|\geq1}|G(r,t) |^2 \nu(\d r)\right].
        \end{aligned}
    \end{equation}
    }Finally, we aggregate  the above four bounds to get the upper bound for the derivative of the KL divergence between $p^\theta_t$ and $\widehat{p}_t$. It is treated for two cases of interaction kernels as below.

{\bf (Bounded interaction)} 
{\footnotesize\begin{equation*}
        \begin{aligned}        
            &\ \frac{\d}{\d t}D_{\mathrm{KL}}(p^\theta_t\|\widehat{p}_t)\leq
            - \frac{1}{6C_\sigma}\int_{\Omega}\left|\nabla\log \left(\frac{p^{\theta}_t(x)}{ \widehat{p}_t(x) }\right)\right|^2 p^\theta_t(x)\d x+ 6C_\sigma (C_K)^2 D_{\mathrm{KL}}(p^\theta_t\| \widehat{p}_t)\\
            &\  + 3C_\sigma C^{\mathrm{nG}}_TD_{\mathrm{KL}}(p^\theta_t\|\widehat{p}_t) + 3 C_\sigma\int_{\Omega}\left| f^\theta(x,t) - \mathcal{V}[p^\theta_t](x,t)\right|^2p^\theta_t(x)\d x\\
            \leq&\ 3\left(2C_\sigma (C_K)^2 +  C_\sigma C^{\mathrm{nG}}_T\right)D_{\mathrm{KL}}(p^\theta_t\| \widehat{p}_t)+  3C_\sigma\int_{\Omega}\left| f^\theta(x,t) - \mathcal{V}[p^\theta_t](x,t)\right|^2p^\theta_t(x)\d x.
        \end{aligned}
    \end{equation*}}

{\bf (Biot--Savart interaction)} 

{\footnotesize\begin{equation*}
        \begin{aligned}        
            &\ \frac{\d}{\d t}D_{\mathrm{KL}}(p^\theta_t\|\widehat{p}_t)\leq -\frac{1}{6C_\sigma}\int_{\Omega}\left|\nabla\log \left(\frac{p^{\theta}_t(x)}{ \widehat{p}_t(x) }\right)\right|^2 p^\theta_t(x)\d x + \bigg( 6C_\sigma  \left\| \nabla\log\widehat{p}_t \right\|^2_\infty \left\|U\right\|^2_\infty \\
            &\ \left. + 4\|U\|_\infty\left\|\frac{\nabla^2\widehat{p}_t}{\widehat{p}_t}\right\|_\infty\right)  D_{\mathrm{KL}}(p^\theta_t\| \widehat{p}_t) + 3C_\sigma C^{\mathrm{nG}}_TD_{\mathrm{KL}}(p^\theta_t\|\widehat{p}_t)\\
            &\  + 3 C_\sigma\int_{\Omega}\left| f^\theta(x,t) - \mathcal{V}[p^\theta_t](x,t)\right|^2p^\theta_t(x)\d x\\
            \leq&\  \left( 6C_\sigma  \left\| \nabla\log\widehat{p}_t \right\|^2_\infty \left\|U\right\|^2_\infty + 4\|U\|_\infty\left\|\frac{\nabla^2\widehat{p}_t}{\widehat{p}_t}\right\|_\infty  + 3C_\sigma C^{\mathrm{nG}}_T\right)D_{\mathrm{KL}}(p^\theta_t\| \widehat{p}_t)\\
            &\ +  3C_\sigma\int_{\Omega}\left| f^\theta(x,t) - \mathcal{V}[p^\theta_t](x,t)\right|^2p^\theta_t(x)\d x.
        \end{aligned}
    \end{equation*}
    }By Gronwall inequality, we  obtain that
{\footnotesize\begin{equation*}
    \begin{split}      &\ \sup_{t\in[t_0,T]}D_{\mathrm{KL}}(p^\theta_t\| \widehat{p}_t)\leq\   \left(3 C_\sigma\int_{t_0}^t\int_{\Omega}\left| f^\theta(x,t) - \mathcal{V}[p^\theta_t](x,t)\right|^2p^\theta_t(x)\d x\d t\right)\exp\left\{\bar{C}T\right\},
    \end{split}
\end{equation*}
}where
    {\footnotesize\begin{equation}\label{eqn:barC}
   \bar{C}=\left\{ \begin{aligned}
        3\left(2C_\sigma (C_K)^2 + C_\sigma C^{\mathrm{nG}}_T\right),\quad \mbox{Interaction}\ \ref{subass:bounded},\\
        \max_{t\in[t_0,T]}\left(6C_\sigma  \left\| \nabla\log\widehat{p}_t \right\|^2_\infty \left\|U\right\|^2_\infty + 4\|U\|_\infty\left\|\frac{\nabla^2\widehat{p}_t}{\widehat{p}_t}\right\|_\infty\right)  + 3C_\sigma C^{\mathrm{nG}}_T,\quad \mbox{Interaction}\ \ref{subass:BS},
    \end{aligned}\right.
\end{equation}
}
\end{proof}

\subsection{Proof of Theorem \ref{theo:error}}
\begin{proof}[Proof of Theorem \ref{theo:error}]
For convenience and simplicity, we write $X^*_{t_0,t}$, $X^N_{t_0,t}$ as $X^*_{t}$, $X^N_{t}$. For $t_0\leq t\leq t_0+(N_T-1)\Delta t$, using Taylor expansion and the smoothness of
the velocity $\mathcal{V}[\widehat{p}_t]$, we have
{\small\begin{equation}\label{differencestar}
    \begin{aligned}
        X^*_{t+\Delta t}(x)=&\ X^*_t(x) + \mathcal{V}[\widehat{p}_t](X^*_t(x),t)\Delta t + \frac{\d^2 X^*_t (x)}{2\d t^2}\bigg|_{t=\tau}(\Delta t)^2\\
        =&\ X^*_t(x) + \frac{\d^2 X^*_t (x)}{2\d t^2}\bigg|_{t=\tau}(\Delta t)^2 + \left( b(X^*_t,t) + \int_{\Omega}K(X^*_t(x)-y)\widehat{p}_t(y)\d y \right.\\
        &\ -\int_{\Omega}F(r,t)\nu(\d r) - \frac{1}{2}\Sigma(X^*_t(x),t)   - \frac{1}{2}\Sigma(X^*_t(x),t)\nabla \log \widehat{p}_t(X^*_t(x)) \\
        &\ \left. + \int_{|r|<1}F(r,t)\frac{\int_0^1\widehat{p}_t(X^*_t(x)-\lambda F(r,t))\d \lambda}{\widehat{p}_t(X^*_t(x))} \nu(\d r)\right)\Delta t,
    \end{aligned}
\end{equation}
}for some $\tau\in[t,t+\Delta t]$. According to Algorithm \ref{algorithm}, we have
{\small\begin{equation}\label{differencedelta}
    \begin{aligned}
       X^N_{t+\Delta t}(x)=&\ X^N_t(x) + \Delta t \bigg[ b(X^N_t(x),t) - s(X^N_t(x),t) - \int_{|r|<1}F(r,t)\nu(\d r)   - \nabla\cdot\frac{\Sigma(X^N_t(x),t)}{2} \bigg],
    \end{aligned}
\end{equation}
}where $\{x^{(i)}\}_{i=1}^N$ is a set of $N$ samples from initial density $p_{t_0}$. Note that $s(x,t)$ is the total ``score'' term, we could separate it into three parts that correspond to the interaction, diffusion and L\'{e}vy losses respectively. For simplicity, we assume the numerical errors of these three components are also bounded by $\varepsilon$ respectively. Let
{\small\begin{equation}\label{sKBL}
    s(x,t)= -s_K(x,t) + \frac{1}{2}\Sigma(x,t)s_B(x,t) - s_L(x,t).
\end{equation}
}Subtracting equation \eqref{differencestar} by \eqref{differencedelta} we have
{\small\begin{equation}\label{diferencestar-N}
    \begin{aligned}
        &\ X^*_{t+\Delta t}(x) - X^N_{t+\Delta t}(x) =\ X^*_{t }(x) - X^N_{t }(x) + \frac{\d^2 X^*_t (x)}{2\d t^2}\bigg|_{t=\tau}(\Delta t)^2 + \nabla b(X_{t,\xi},t) \\
        &\ \cdot(X^*_{t }(x) - X^N_{t }(x))\Delta t - \frac{1}{2}\nabla\Sigma(X_{t,\xi'})\cdot(X^*_{t }(x) - X^N_{t }(x))\Delta t + e_K\Delta t -\frac{1}{2}e_B \Delta t  + e_L\Delta t,
    \end{aligned}
\end{equation}
}where $e_K,\ e_B$ and $e_L$ are given in the following. First, based on the law of large numbers and the assumption that the kernel function is twice-diﬀerentiable, for the bounded interaction \eqref{subass:bounded},
{\small\begin{equation*}
    \begin{aligned} \ 
\mathbb{E}_{x}|e_K|&:=\mathbb{E}_{x}\left|\int_{\Omega}K(X^*_t(x)-y)\widehat{p}_t(y)\d y - s_K(X^N_t(x),t) \right| 
\\
       & \leq \ \mathbb{E}_{x}\left|\int_{\Omega}K(X^*_t(x)-y)\widehat{p}_t(y)\d y - \int_{\Omega}K(X^N_t(x)-y)\widehat{p}_t(y)\d y \right| \\
        &\ + \mathbb{E}_{x}\left|\int_{\Omega}K(X^N_t(x)-y)\widehat{p}_t(y)\d y  -s_K(X^N_t(x),t) \right| 
        \leq\ C_1\mathbb{E}_x\left|X^*_t(x)-X^N_t(x)\right| + \varepsilon,
    \end{aligned}
\end{equation*}
}where $C_1=\ 2\sup_{x\in\Omega}|\nabla \cdot K(x)|<\infty$.  
For the Biot--Savart interaction \eqref{subass:BS},
{\small\begin{equation*}
    \begin{aligned}
       \mathbb{E}_{x}\left|e_K \right| 
        \leq &\ \mathbb{E}_{x}\left|\int_{\Omega}K(X^*_t(x)-y)\widehat{p}_t(y)\d y - \int_{\Omega}K(X^N_t(x)-y)\widehat{p}_t(y)\d y \right| \\
        &\ + \mathbb{E}_{x}\left|\int_{\Omega}K(X^N_t(x)-y)\widehat{p}_t(y)\d y  -s_K(X^N_t(x),t) \right| \\
        =&\ \mathbb{E}_{x}\left|\int_{\Omega}K(z)\widehat{p}_t(z + X^*_t(x) )\d z - \int_{\Omega}K(z)\widehat{p}_t(z + X^N_t(x))\d z \right| \\
        &\ + \mathbb{E}_{x}\left|\int_{\Omega}K(X^N_t(x)-y)\widehat{p}_t(y)\d y  -s_K(X^N_t(x),t) \right| \\
        =&\ \mathbb{E}_{x}\left|\int_{\Omega}-U(z)\cdot\nabla\widehat{p}_t(z + X^*_t(x) )\d z + \int_{\Omega}U(z)\cdot\nabla\widehat{p}_t(z + X^N_t(x))\d z \right| \\
        &\ + \mathbb{E}_{x}\left|\int_{\Omega}K(X^N_t(x)-y)\widehat{p}_t(y)\d y  -s_K(X^N_t(x),t) \right| \\
        \leq&\ C_2\mathbb{E}_x\left|X^*_t(x)-X^N_t(x)\right| + \varepsilon,
    \end{aligned}
\end{equation*}
}where $C_2:= \sup_{t\in[t_0,T]}\left\|\nabla_x\left( \int_{\Omega}U(z)\cdot\nabla\widehat{p}_t(z + x)\d z\right)\right\|_{\infty} <\infty$.

Second,
{\small\begin{equation*}
    \begin{aligned}
       \left|e_B\right| :=&\ \left|\Sigma(X^*_t(x),t)\nabla \log \widehat{p}_t(X^*_t(x)) - \Sigma(X^N_t(x),t)s_B(X^N_t(x),t)\right|\\
        \leq&\ \left|\Sigma(X^*_t(x),t)\nabla \log \widehat{p}_t(X^*_t(x)) - \Sigma(X^*_t(x),t)\nabla \log \widehat{p}_t(X^N_t(x))\right|\\
        &\ + \left|\Sigma(X^*_t(x),t)\nabla \log \widehat{p}_t(X^N_t(x)) - \Sigma(X^N_t(x),t)\nabla \log \widehat{p}_t(X^N_t(x))\right|\\
        &\ + \left| \Sigma(X^N_t(x),t)\nabla \log \widehat{p}_t(X^N_t(x))- \Sigma(X^N_t(x),t)s_B(X^N_t(x),t)\right|\\
        \leq&\ C_3|X^*_t-X^N_t| + \varepsilon,
    \end{aligned}
\end{equation*}
}where $C_3= \sup_{t\in[t_0,T]}\sup_{x\in\Omega}\left(C_\sigma\left|\nabla^2\log \widehat{p}_t(x) \right| + \|\nabla\log \widehat{p}_t\|_\infty|\nabla\cdot\Sigma(x,t)|\right)<\infty$. 

Third,
{\small\begin{equation*}
    \begin{aligned}
        |e_L|:=&\ \left| \int_{|r|<1}F(r,t)\frac{\int_0^1\widehat{p}_t(X^*_t(x)-\lambda F(r,t))\d \lambda}{\widehat{p}_t(X^*_t(x))} \nu(\d r) \right.\\
        &\ \left.+ \int_{|r|\geq1}G(r,t)\frac{\int_0^1\widehat{p}_t(X^*_t(x) -\lambda G(r,t))\d \lambda}{\widehat{p}_t(X^*_t(x))} \nu(\d r)- s_L(X^N_t(x),t)  \right|\\
        \leq&\ \left| \mathcal{I}_t(\widehat{p}_t,F)(X^*_t(x)) - \mathcal{I}_t(\widehat{p}_t,F)(X^N_t(x))  \right| + \left|\mathcal{I}_t(\widehat{p}_t,G)(X^*_t(x)) - \mathcal{I}_t(\widehat{p}_t,G)(X^N_t(x))  \right|\\
        &\ + \left|  \mathcal{I}_t(\widehat{p}_t,F)(X^*_t(x)) - \mathcal{I}_t(\widehat{p}_t,F,G)(X^N_t(x))  - s_L(X^N_t(x),t) \right|\\
        \leq&\ C_4|X^*_t(x)-X^N_t(x)| + \varepsilon,
    \end{aligned}
\end{equation*}
}where  $C_4=\  \sup_{t\in[t_0,T]}C^f_TC^*_T\left\|\nabla \widehat{p}_t\right\|_{\infty} \left(\int_{|r|<1}|F(t)r|\d r + \int_{|r|\geq1}|G(t)r|\d r\right) <\infty$.

Now from \eqref{diferencestar-N}, we have, for some positive constant $C$,
{\small\begin{equation*}
    \begin{aligned}
        \mathbb{E}_x\left|X^*_{t+\Delta t}-X^N_{t+\Delta t}\right|\leq&\  \mathbb{E}_x\left|X^*_{t}-X^N_{t}\right|  + C\mathbb{E}_x\left|X^*_{t}-X^N_{t}\right|\Delta t +  3\varepsilon \Delta t +\mathcal{O}((\Delta t)^2),\\
        C=&\ \left\{\begin{aligned}
             \|\nabla b\|_\infty + \frac{1}{2}\|\Sigma(\cdot,t)\|_\infty + C_1 + C_3+ C_4 ,\quad \mbox{Interaction}\ \ref{subass:bounded},\\
              \|\nabla b\|_\infty + \frac{1}{2}\|\Sigma(\cdot,t)\|_\infty + C_2 + C_3 + C_4 ,\quad \mbox{Interaction}\ \ref{subass:BS}.
        \end{aligned}\right. 
    \end{aligned}
\end{equation*}
}Therefore $\mathbb{E}_x\left|X^*_{t}-X^N_{t}\right|\leq e^{C(t-t_0)}\mathbb{E}_x\left|X^*_{t_0}-X^N_{t_0}\right| + \left[ \mathcal{O}(\varepsilon)+\mathcal{O}(\Delta t)\right](t-t_0)$, as $\varepsilon$ and $\Delta t$ all tend to 0. With the same initial condition, we immediately obtain \eqref{numerror}.
\end{proof}

    \section{Numerical Experiments}
    \label{sec:example}
    We examine several examples to demonstrate the effectiveness of our algorithm \ref{algorithm}.
The time interval \([0, T]\) is uniformly partitioned  into \( N_T \) sub-intervals \([t_k, t_{k+1}]\), where   \( t_k = k\frac{T}{N_T} \) for \( k = 0, 1, \dots, N_T \). 
On each sub-interval \([t_k, t_{k+1}]\), the transport map is approximated by a neural network \( s^{\theta_{k}}(\cdot,t_k): \mathbb{R}^d \to \mathbb{R}^d \), modeled as a multi-layer perceptron (MLP)  with $3$ hidden layers, $32$ neurons per layer, and the $Swish$ activation function. 
The algorithm is implemented with the following parameter settings: the time step size of \(\Delta t = T/N_T=t_{k+1}-t_{k}  = 10^{-3}\), the time horizon of \(T=1\), and the sample size of \(N=4000\). 

The initial condition $\widehat{p}_0$ of \eqref{sdeBL} in all examples is set as the   Gaussian distribution for its simplicity in generating initial samples \( \{ X_0^{(i)} \}_{i=0}^{N} \) (unless otherwise specified, it is assumed to be the standard normal distribution). 
For each time step in training the score function $s^{\theta_{k+1}}(\cdot )$, we use the warm start for the optimization by initializing the neural network parameter \( \theta_{{k+1}} \) by 
the obtained parameters \( \theta^*_{k} \) from the previous step, followed by  the standard   the Adam optimizer with a learning rate of \( 10^{-4} \) to optimize  \( \theta_{{k+1}} \). 

To evaluate our method, we 
use   the total variation (TV) distance of  the generated samples to compare out method  with the Monte Carlo simulation method.
At each time \(t_k\), we identify the smallest rectangular domain \(\Omega_{t_k}\) covering all sample points and discretize it into  uniform grid cells \(\{\Delta_i\}\). 
The distribution for each method is then approximated by a histogram:  \(P(\Delta_i) = \frac{\#(\text{samples} \in \Delta_i)}{\#\text{samples}}\). 
We denote these binned empirical distributions as \(P^{\mathrm{MC}}\) for Monte Carlo and \(P^{\mathrm{S}}\) for our method, and the TV distance between these two distributions is then numerically computed by {\small\begin{equation*}
    d_{\mathrm{TV}}(P^{\mathrm{MC}}, P^{\mathrm{S}}) = \sum_{i} \left|P^{\mathrm{MC}}(\Delta_i) - P^{\mathrm{S}}(\Delta_i)\right|.
\end{equation*}}

\begin{example}[1D Jump-diffusion with  finite jump activity]\label{example1}
To test the efficiency of our proposed method, we consider the following one-dimensional SDE:{\small\begin{equation}\label{exa:1dfinite}
        \begin{aligned}
            \d X_t = \kappa(\eta - X_{t-})\d t + \sigma\d B_t + J_t\d N_t,
        \end{aligned}
    \end{equation} 
}which has been applied in \cite{chevallier2017estimation} to model \( \mathrm{CO}_2 \) emissions and fuel-switching strategies. 
In this example, we set \( \kappa = 1 \), \( \eta = 1 \), and \( \sigma = 2 \). The process \( N_t \) is assumed to be a Poisson process with intensity \( \Lambda = 30 \), while the jump sizes \( J_t \) are drawn from a Gaussian distribution \( \mathrm{N}_{m, v^2} \) with mean \( m = 0.1 \) and standard deviation \( v = {1}/{24} \). To compare with \eqref{sdeBL}, the jump term \( J_t \d N_t \) is equivalently written as \( \int_{\mathbb{R}} r \mathcal{N}(\d t, \d r) \), i.e., \( F(r,t) = G(r,t) = r \), where \( \mathcal{N} \) is a Poisson random measure with L\'{e}vy measure \( \nu(\d r) = \Lambda \mathrm{N}_{m, v^2}(\d r) \), and the numerical integration with respect to 
$\mathrm{N}_{m, v^2}(\d r)$ is based on a quadrature formula  within a truncated interval 
\([m-3v, m+3v]\).



   
Figure \ref{fig:example1-probability} illustrates the temporal evolution of probability flows and probability density functions for Equation \eqref{exa:1dfinite}, obtained via both the Monte Carlo simulation and the proposed method. Specifically, the Monte Carlo simulation employs the following Euler–Maruyama discretization scheme:  
\(
X_{t+\Delta t}^{(i)} = X_t^{(i)} + \kappa(\eta - X_t^{(i)})\Delta t + \sigma \xi_t + \sum_{k=1}^{N_{\Delta t}} J_k, \quad i = 1, \dots, N,
\)
where \( \xi_t \sim \mathrm{N}_{0, \Delta t} \), \( N_{\Delta t} \sim \text{Po}(\Lambda \Delta t) \) is a Poisson random variable with rate \( \Lambda \Delta t \), and \( J_k \)'s are i.i.d.\ random variables for jump sizes distributed as $\mathrm{N}_{m,v^2}$.

\begin{figure}[ht] 
    \centering
   \includegraphics[width=0.49\textwidth]{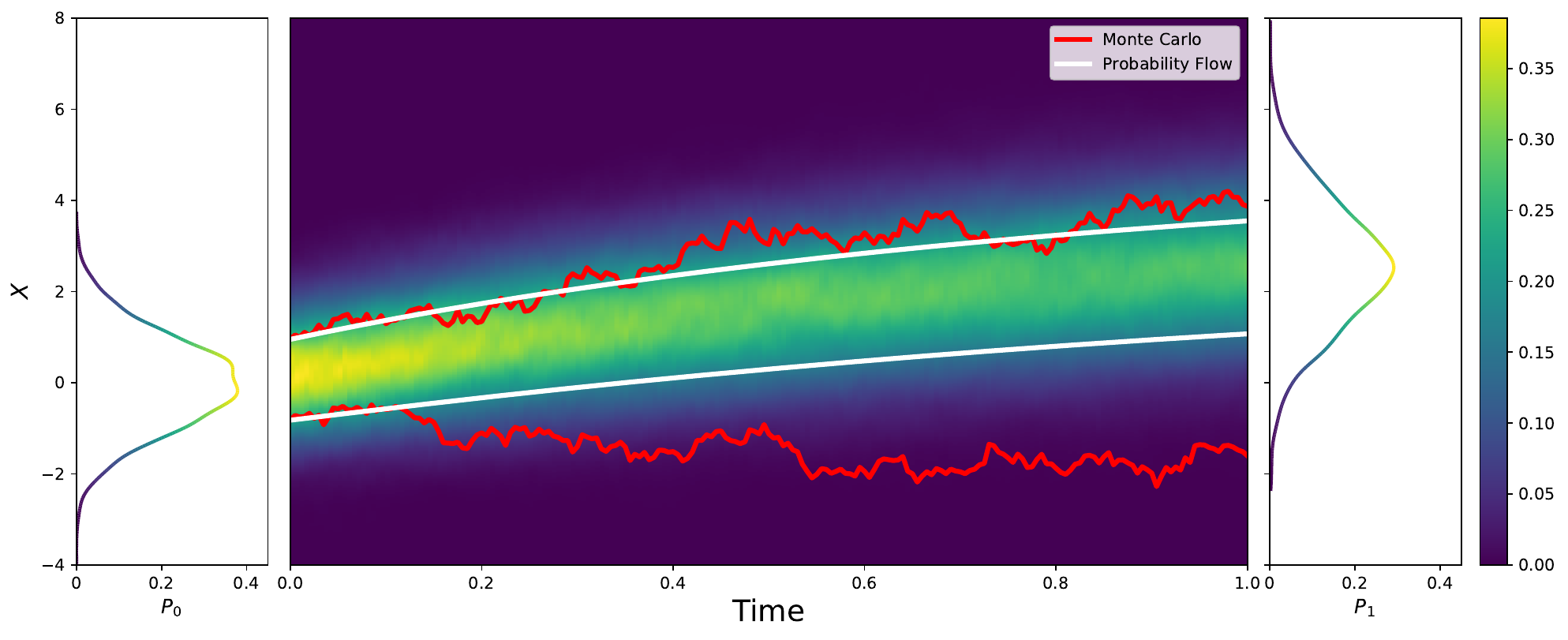}
   \centering
   \includegraphics[width=0.49\textwidth]{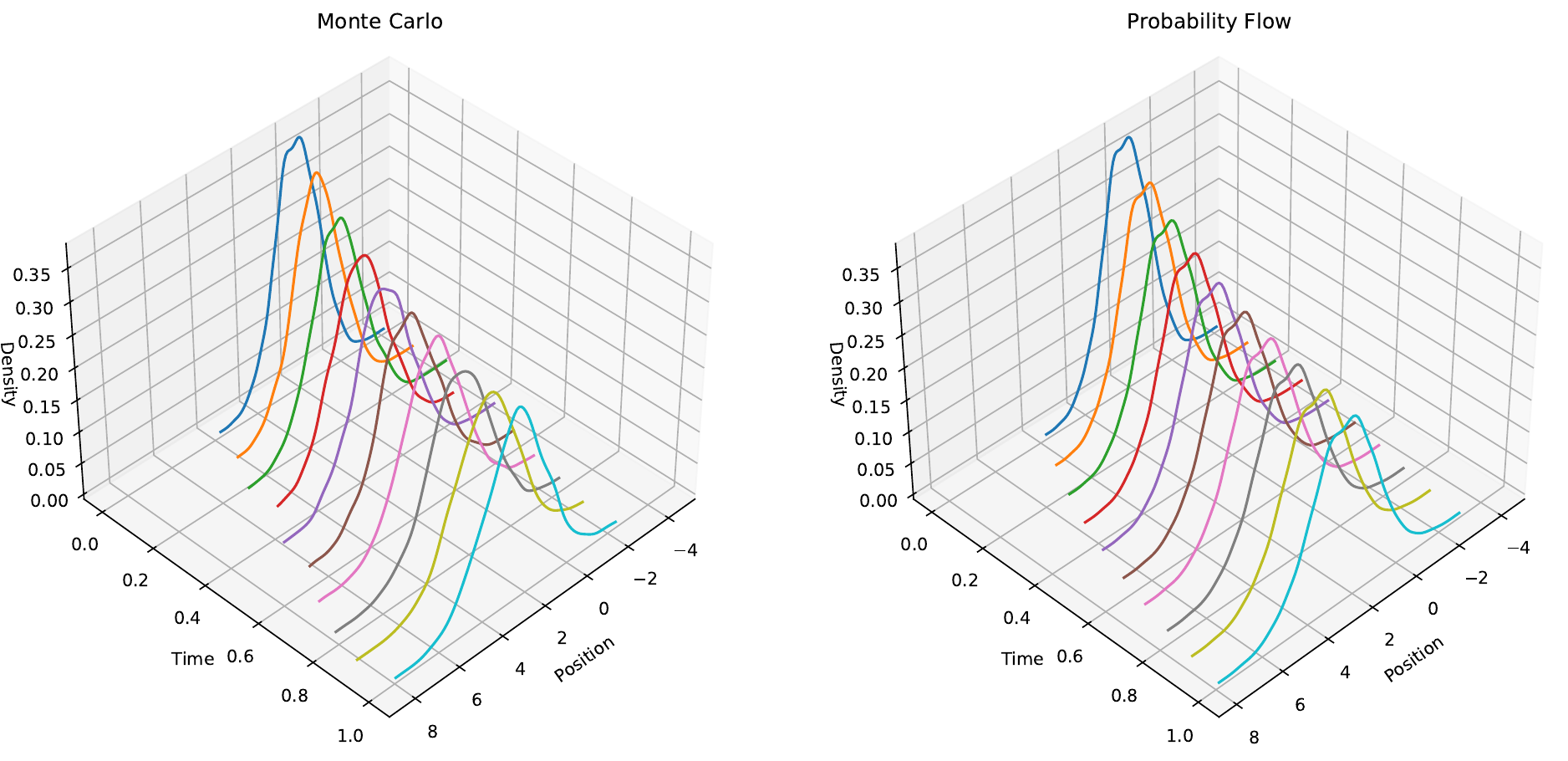}
   \caption{[Example \ref{example1}] Probability flows of \eqref{exa:1dfinite}. The left panel illustrates the temporal evolution of the probability distribution as a heat map, overlaid with two selective stochastic trajectories based on the Monte Carlo simulation (red) and the two deterministic trajectories based on the transport map (white). The middle and right panels compare the probability distributions $\hat{p}_t(x)$ from the Monte Carlo simulation and the proposed method in the time-state space.}
\label{fig:example1-probability}
\end{figure}

Figure \ref{fig:example1-TV} illustrates the TV distance between \( P^{\mathrm{MC}} \) and \( P^{\mathrm{S}} \). 
The TV distance remains consistently on the order of \( 10^{-2} \), which demonstrates the relative accuracy of our method across time.

\begin{figure}[ht]
   \centering
   \includegraphics[width=0.55\textwidth]{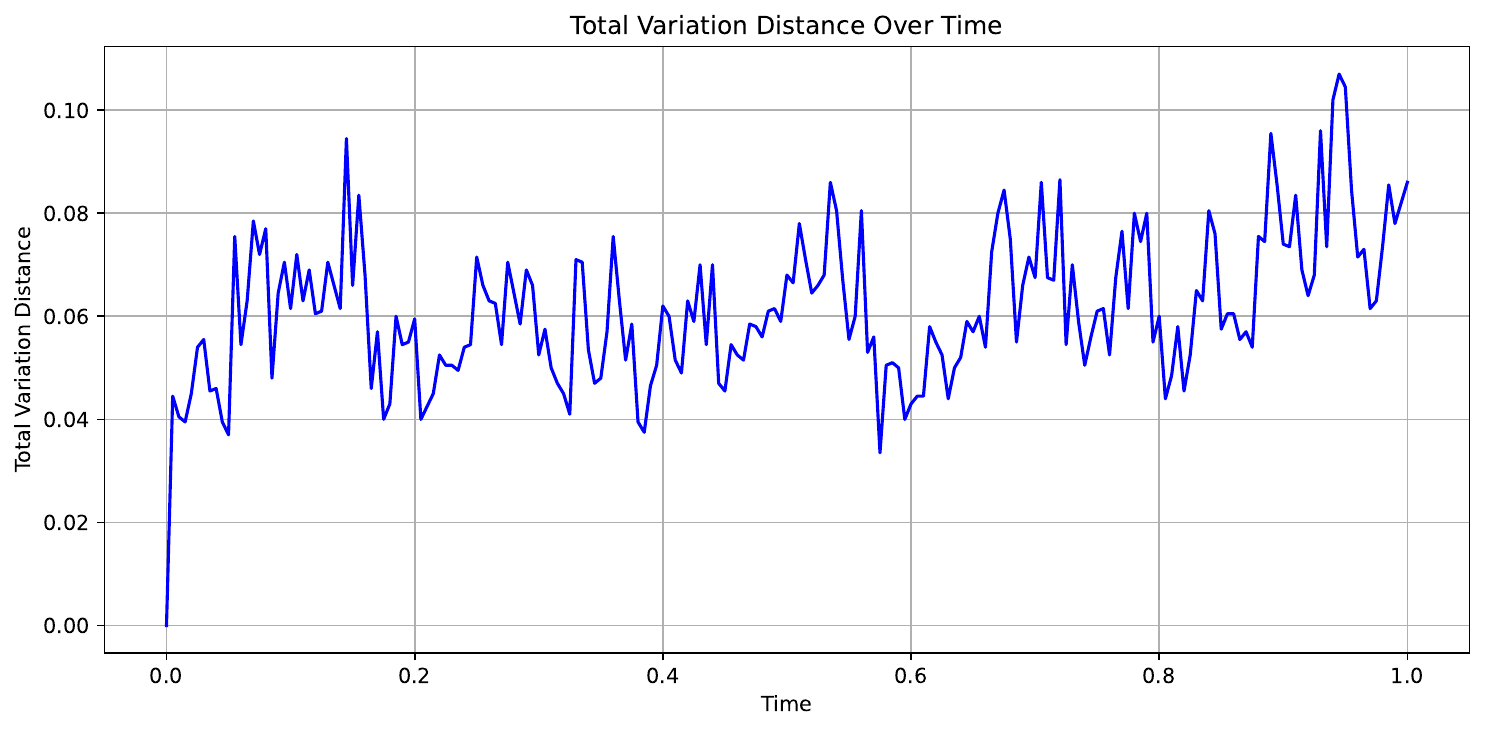}
    \caption{[Example \ref{example1}] The TV distance between $P^{\mathrm{MC}}$ and $P^{\mathrm{S}}$ for \eqref{exa:1dfinite}.}
    \label{fig:example1-TV}
\end{figure}

\end{example}

\begin{example}[Jump-diffusion with $\alpha$-stable L\'{e}vy noise]\label{example2}
In this example, we consider the following one-dimensional SDE:
    {\small\begin{equation}\label{eqn:exampe2}
        \begin{aligned}
            \d X_t = \kappa(\eta - X_{t-})\d t + \sigma\d B_t +  \d L^\alpha_t,
        \end{aligned}
    \end{equation}
   }which has been widely used to capture dynamics with both continuous fluctuations and heavy-tailed, discontinuous shocks, making them relevant in fields such as finance, physics \cite{ariga2021noise}, and generative modeling \cite{yoon2023score}.
   This system captures mean reversion, diffusion, and  an $\alpha$-stable L\'{e}vy motion $L^\alpha_t$ with L\'{e}vy measure: 
   {\small \begin{equation*}
        \nu(r)= c_\alpha\frac{1}{|r|^{\alpha+1}},\quad c_\alpha=\frac{\alpha}{2^{1-\alpha}\sqrt{\pi}}\frac{\Gamma(\frac{1+\alpha}{2})}{\Gamma(1-\frac{\alpha}{2})},\quad \Gamma(z):= \int_0^\infty t^{z-1}e^{-t}\d t.
    \end{equation*}
    }

    \begin{figure}[ht]
    \centering
   \includegraphics[width=0.495\textwidth]{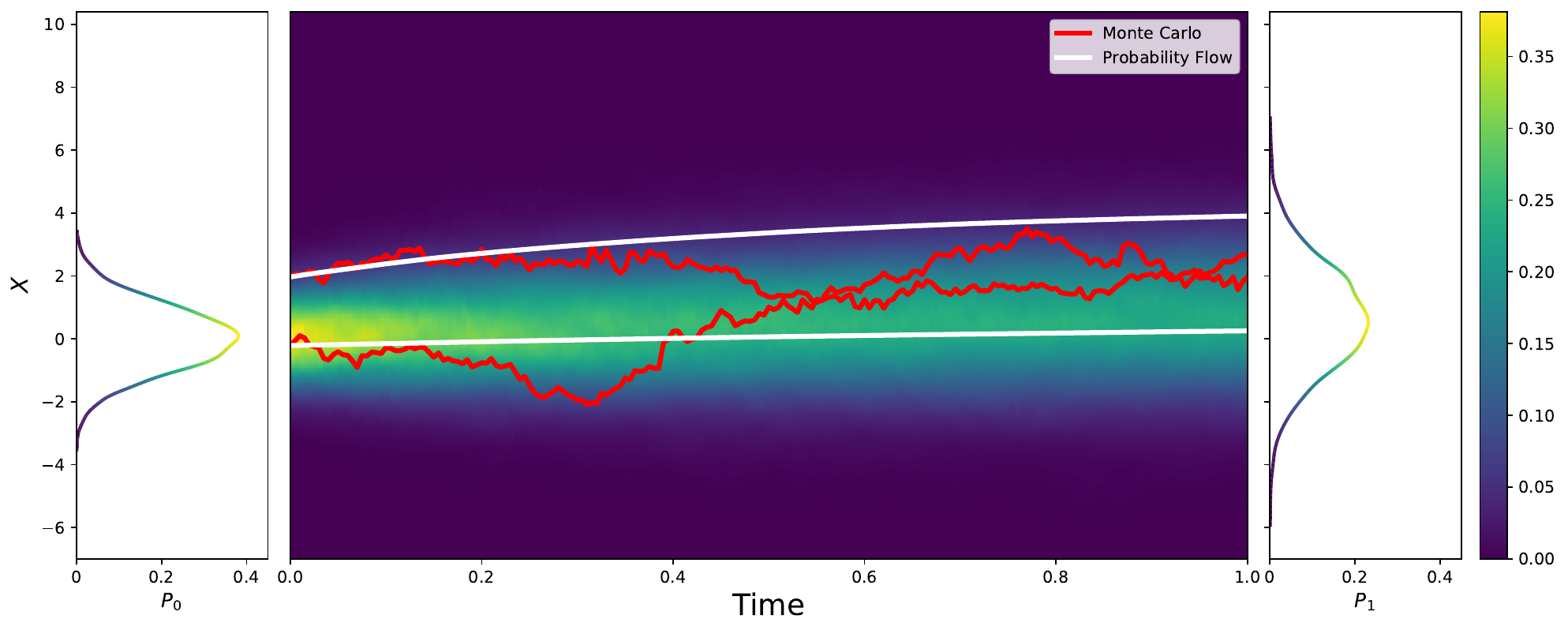}
   \centering
   \includegraphics[width=0.495\textwidth]{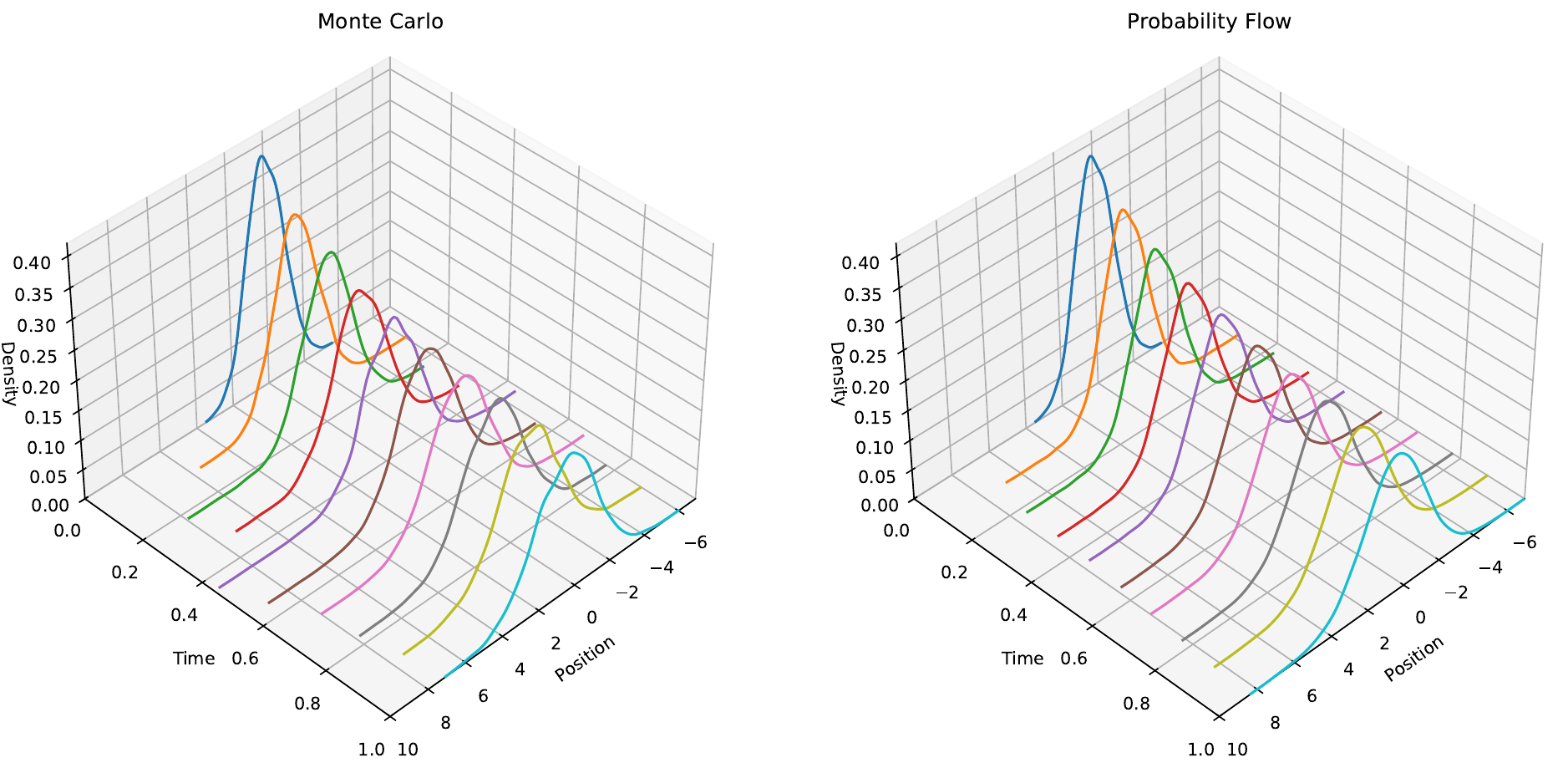}
   \caption{[Example \ref{example2}] Probability flows of \eqref{eqn:exampe2}.} 
    \label{fig:example2-probability}
\end{figure}

We set \(\alpha = 1.5\) while keeping all other hyperparameters the same as in Example \ref{example1}. To avoid the singularity of the intensity measure \(\nu \) at $0$, the numerical integration with respect to \(\nu(r) \) is restricted to the interval \([ -5, -0.01 ] \cup [ 0.01, 5 ]\). 
In this example, the Monte Carlo simulation is based on the following Euler–Maruyama scheme: $X_{t+\Delta t}^{(i)} = X_t^{(i)} + \kappa(\eta - X_t^{(i)})\Delta t  + \sigma \xi_t + \Delta L^\alpha_{t},\ i=1,\cdots,N$, where the jump increments \(\Delta L_t^\alpha = L_{t+\Delta t}^\alpha - L_t^\alpha\) follow a stable distribution \(S_\alpha(  \Delta t^{\frac{1}{\alpha}},0,0)\).

Figure \ref{fig:example2-probability} and Figure \ref{fig:example2-TV} depict the evolutions of probability flow and the TV distance between \( P^{\mathrm{MC}} \) and \( P^{\mathrm{S}} \). The proposed method demonstrates convincing performance in handling \(\alpha\)-stable L\'{e}vy processes and singular intensity measures, maintaining a TV distance on the order of \( 10^{-2} \). 



\begin{figure}[ht]
   \centering
   \includegraphics[width=0.55\textwidth]{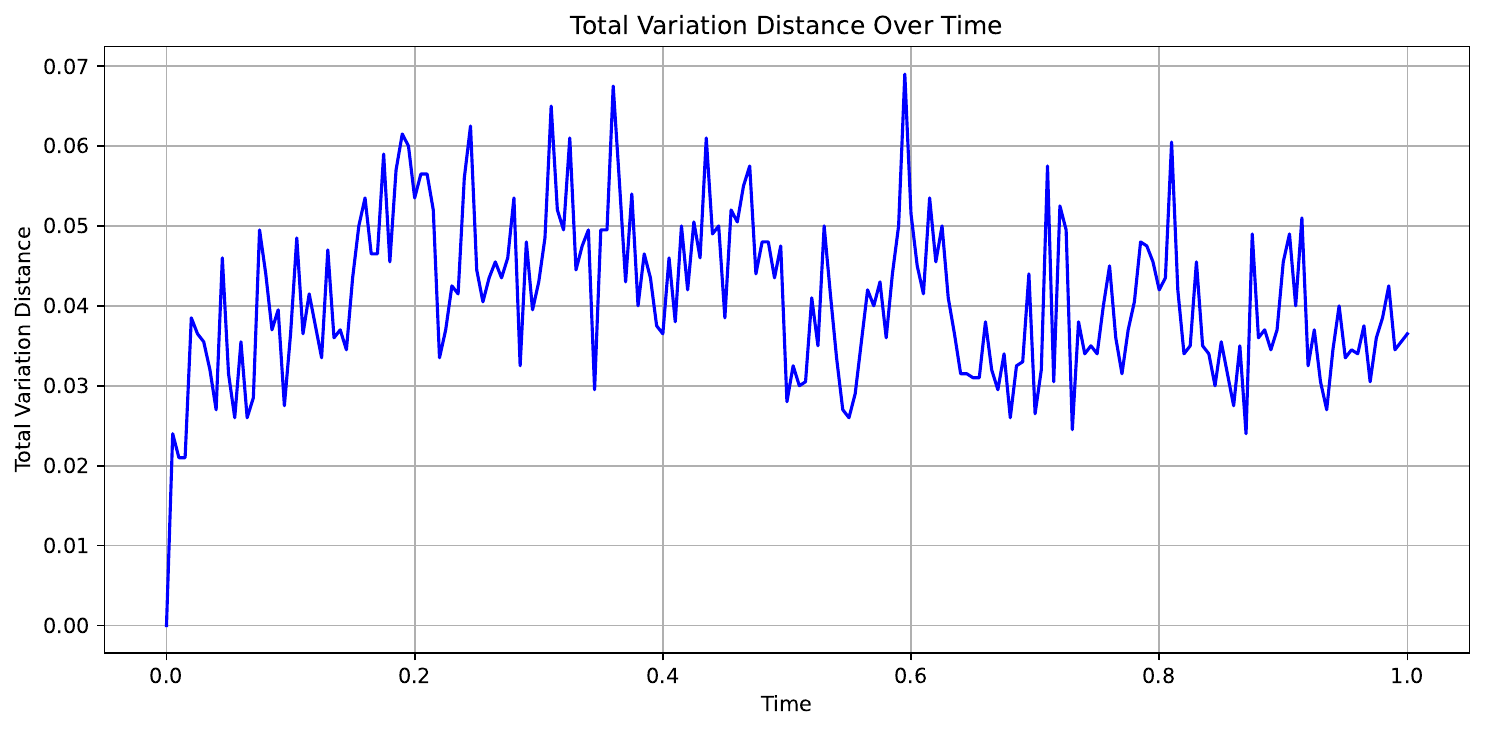}
    \caption{[Example \ref{example2}] The TV distance between $P^{\mathrm{MC}}$ and $P^{\mathrm{S}}$ for \eqref{eqn:exampe2}.}
    \label{fig:example2-TV}
\end{figure}
\end{example}

\begin{example}[Double-well system with interaction force]\label{example3}
    Next, we consider the following  two-dimensional interactive SDE:
    {\small\begin{equation}\label{exa:doublewell}
        \begin{aligned}
            \d X_t=&\ (X_{t-}-X_{t-}^3)\d t + K*\widehat{p}_{t}(X_{t-})\d t + \sigma\d B^1_t +  J_t\d N_t,\\
            \d Y_t=&\ Y_{t-}\d t + K*\widehat{p}_{t}(Y_{t-})\d t + \sigma\d B^2_t,
        \end{aligned} 
    \end{equation}
   }where the interaction kernel is $K(x) = x$.  We set the parameters $\sigma$ and $J_t$ to be the same as those in Example \ref{example1}, and the Monte Carlo simulation is based on the following Euler--Maruyama method:
{\footnotesize\begin{equation}\label{MC:example2dim}
    \begin{aligned}
        X_{t+\Delta t}^{(i)} =&\ X_t^{(i)} + (X_t^{(i)} - (X_t^{(i)})^3)\Delta t  + \frac{1}{N}\sum_{j=1}^N K(X_t^{(i)}-X_t^{(j)})\Delta t+ \sigma \xi^1_t + \sum_{k=1}^{N_{\Delta t}} J_k,\\
        Y_{t+\Delta t}^{(i)} =&\ Y_t^{(i)} + Y_t^{(i)}\Delta t  + \frac{1}{N}\sum_{j=1}^N K(Y_t^{(i)}-Y_t^{(j)})\Delta t+ \sigma \xi^2_t ,\quad i=1,\cdots,N,
    \end{aligned}
\end{equation}
}where $\xi^1_t$ and $\xi^2_t$ are independent random variables sampled from the Gaussian distribution $\mathrm{N}_{0,\Delta t}$.

    \begin{figure}[htbp]
    \centering
   \includegraphics[scale=0.18]{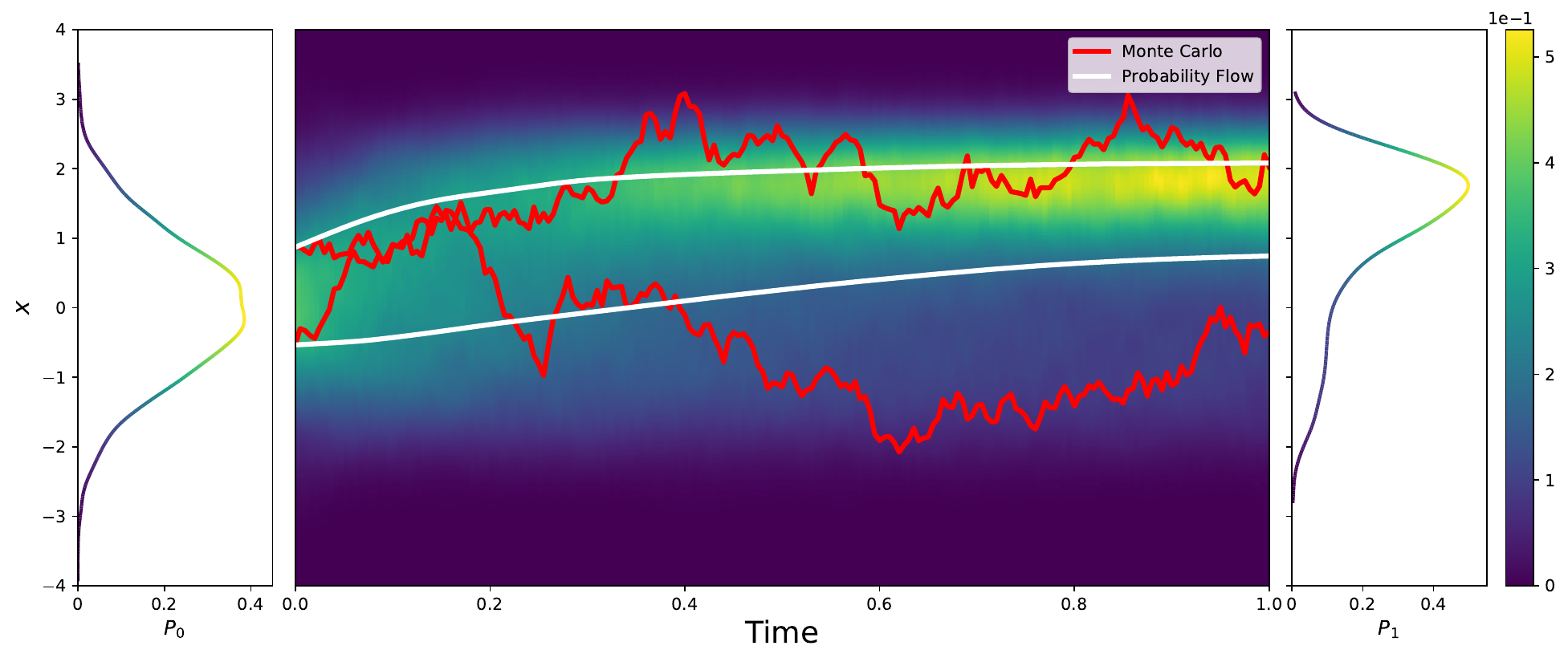}
    \centering
   \includegraphics[scale=0.18]{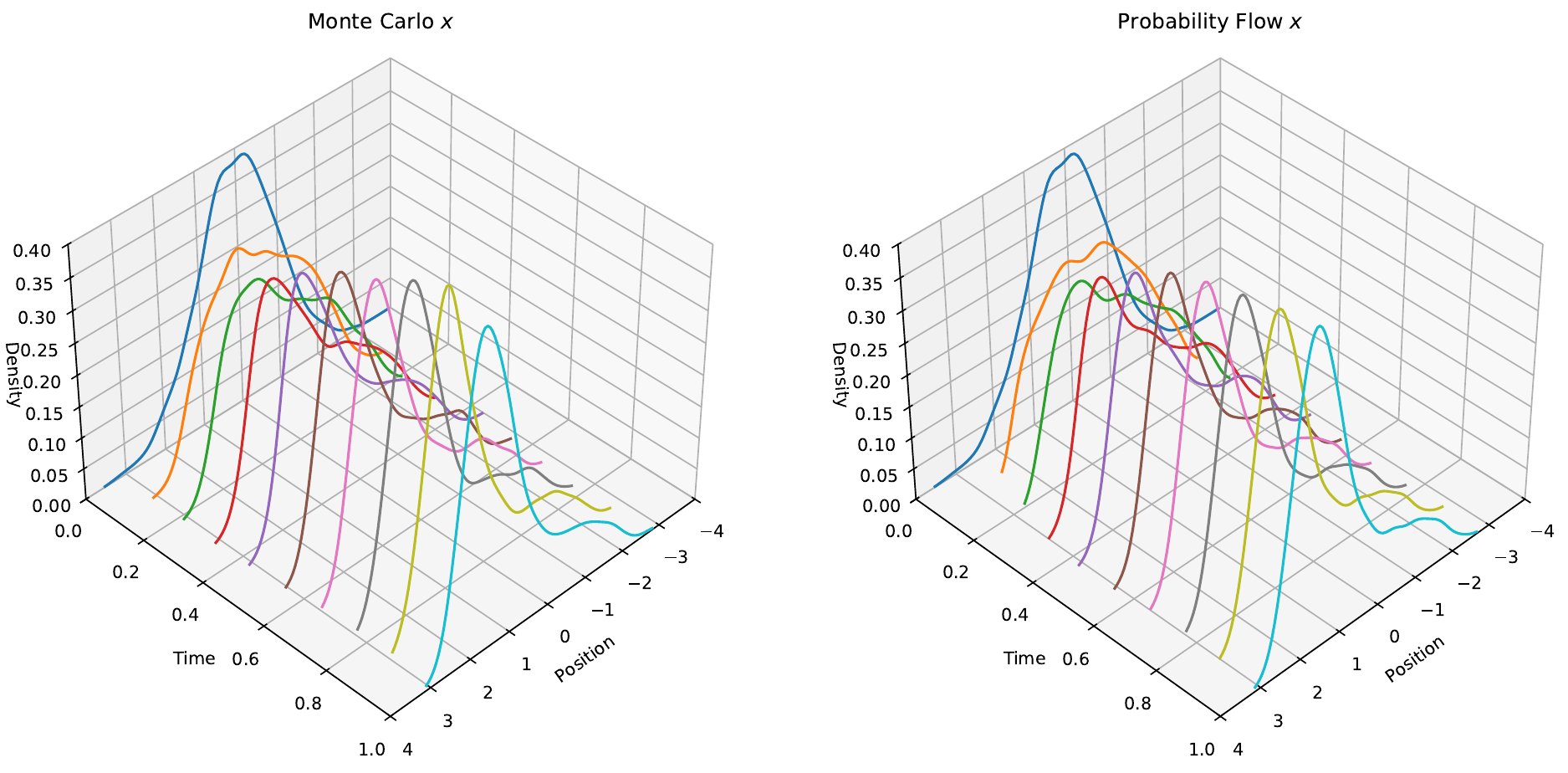}
    \centering
   \includegraphics[scale=0.18]{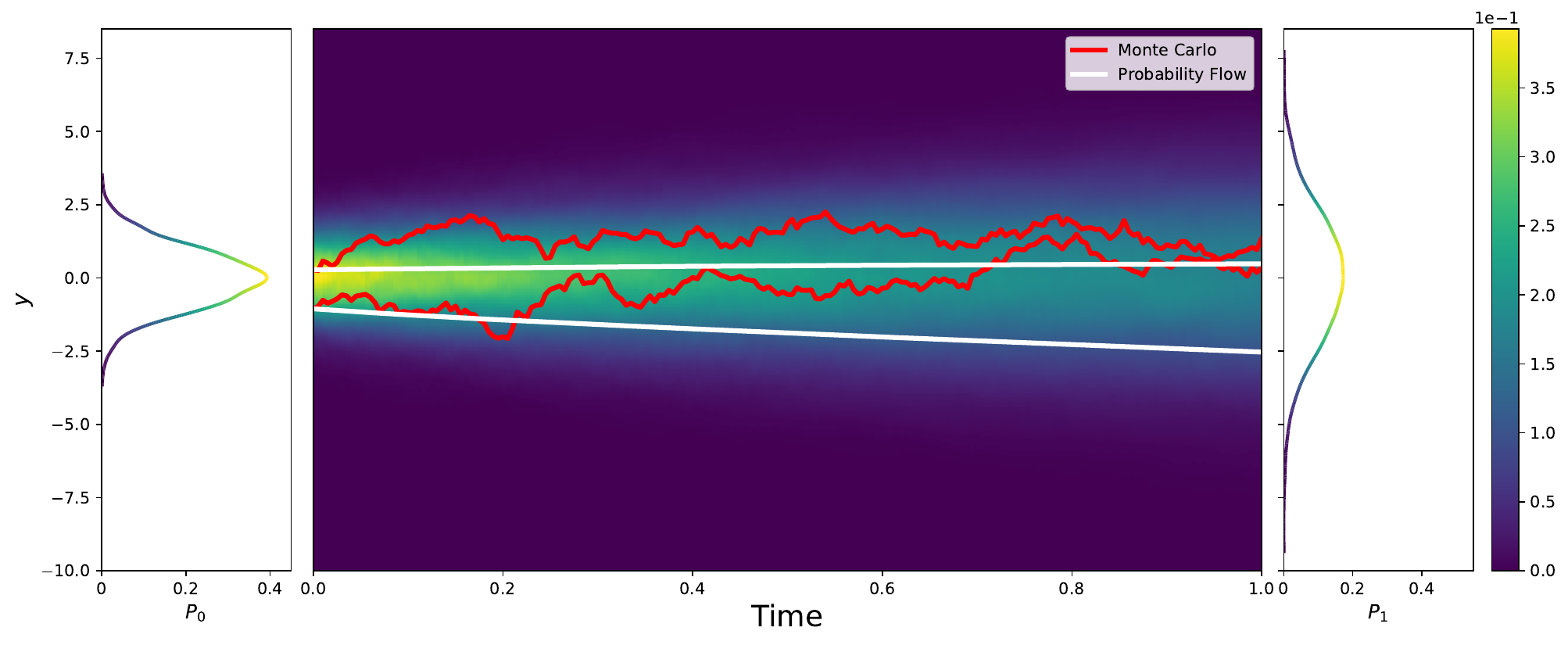}
   \centering
   \includegraphics[scale=0.18]{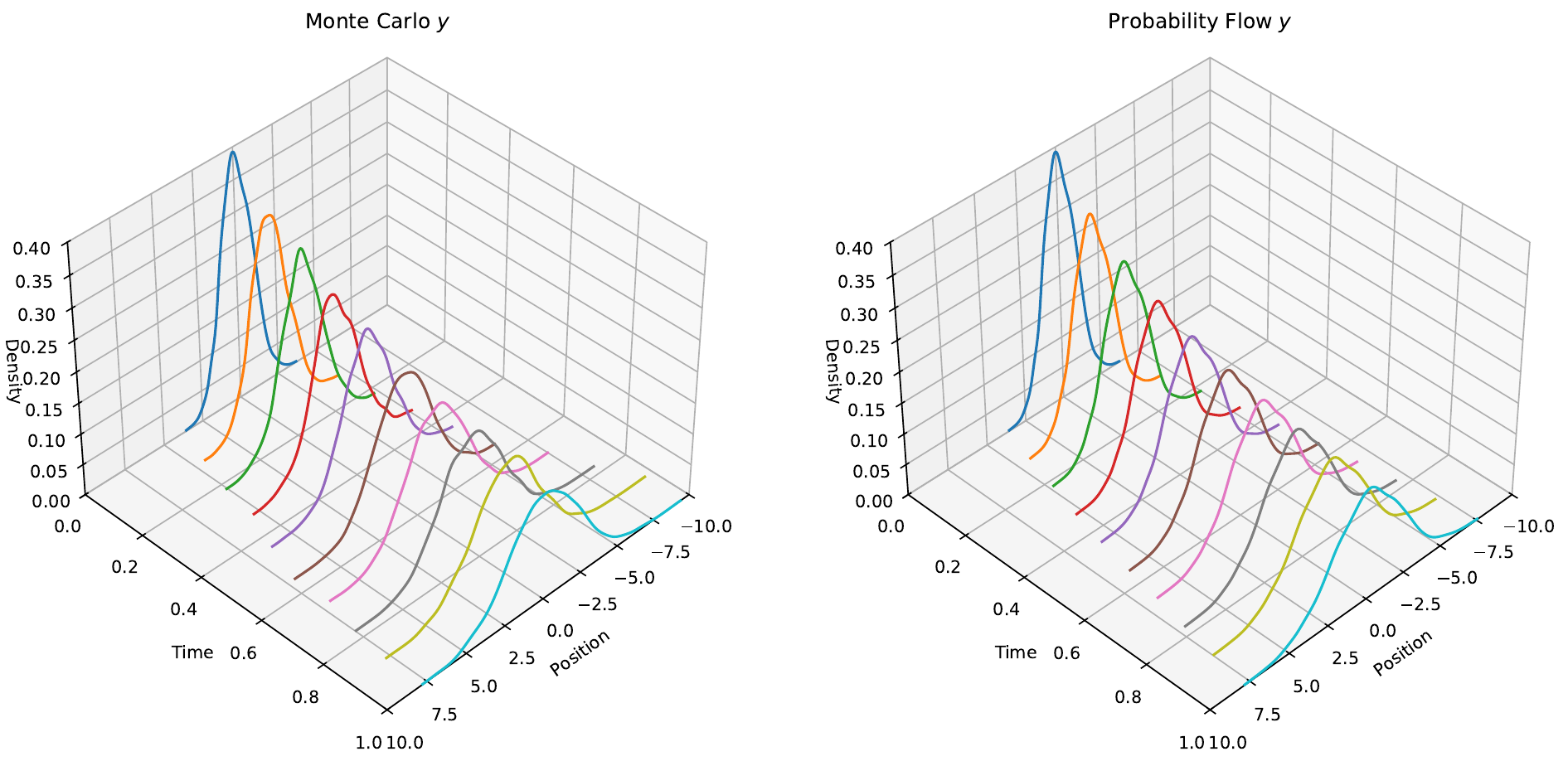}
   \caption{[Example \ref{example3}] The Probability flows of \eqref{exa:doublewell} along the axes: $x$ (top), and $y$ (bottom).}
\label{fig:example2dim-probability}
\end{figure}

\begin{figure}[htbp]
    \centering
   \includegraphics[scale=0.19  ]{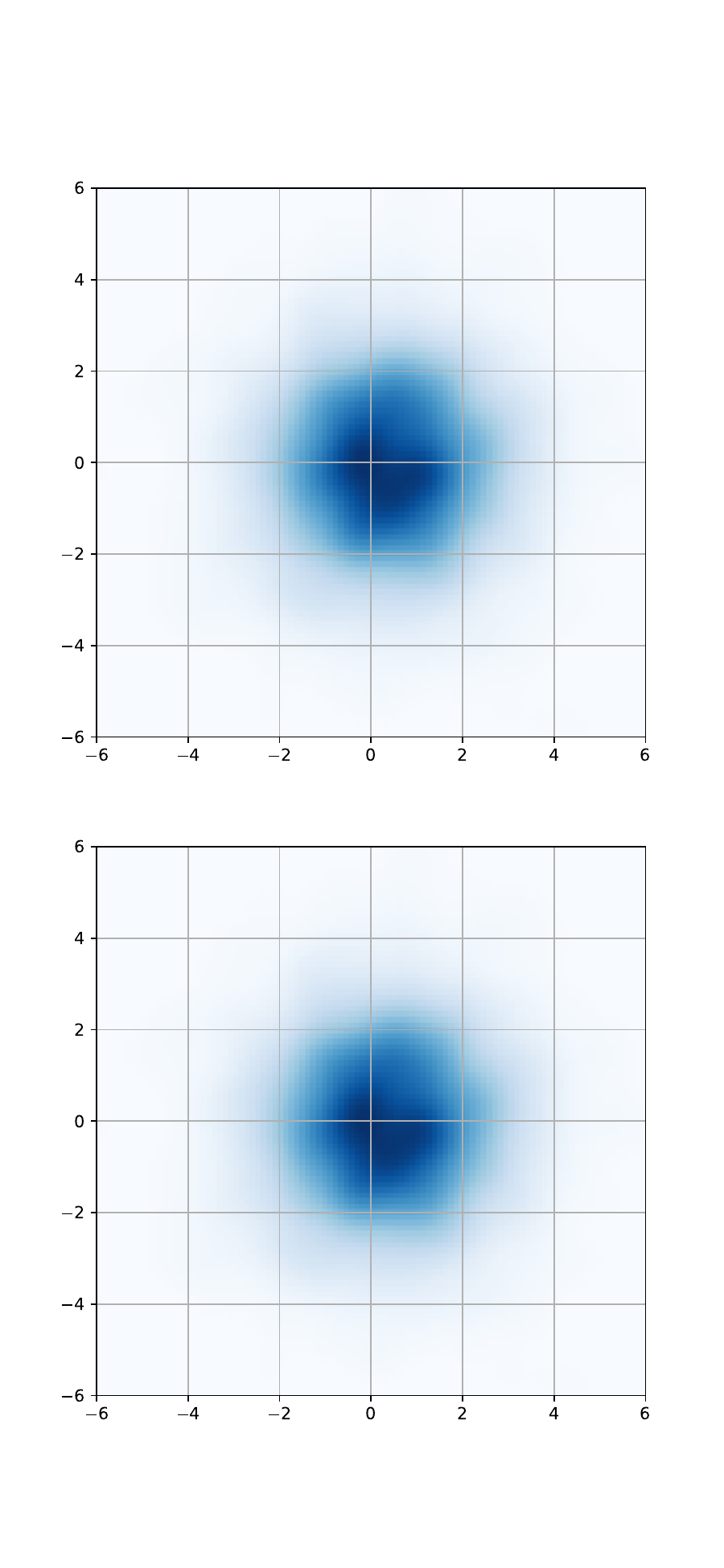}
    \centering
   \includegraphics[scale=0.19 ]{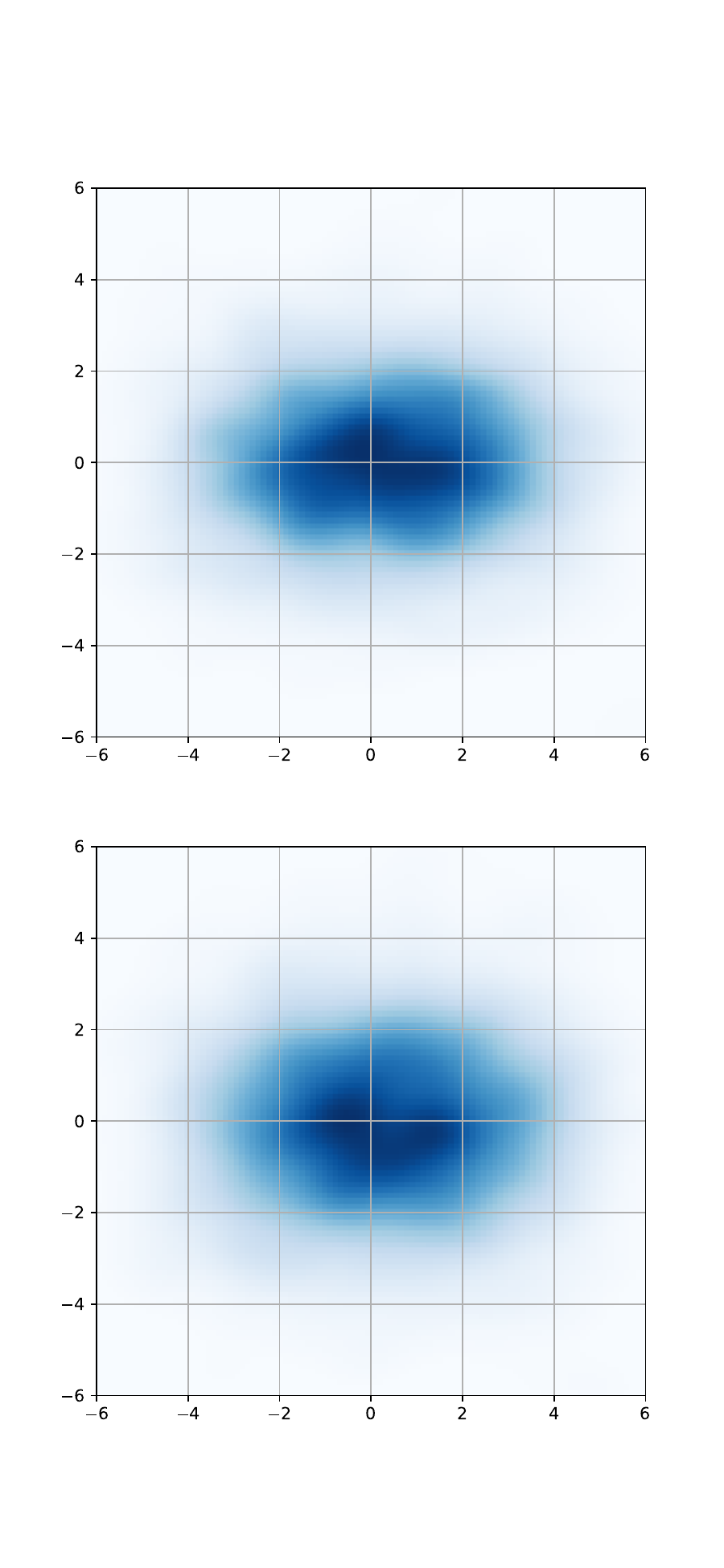}
   \centering
   \includegraphics[scale=0.19 ]{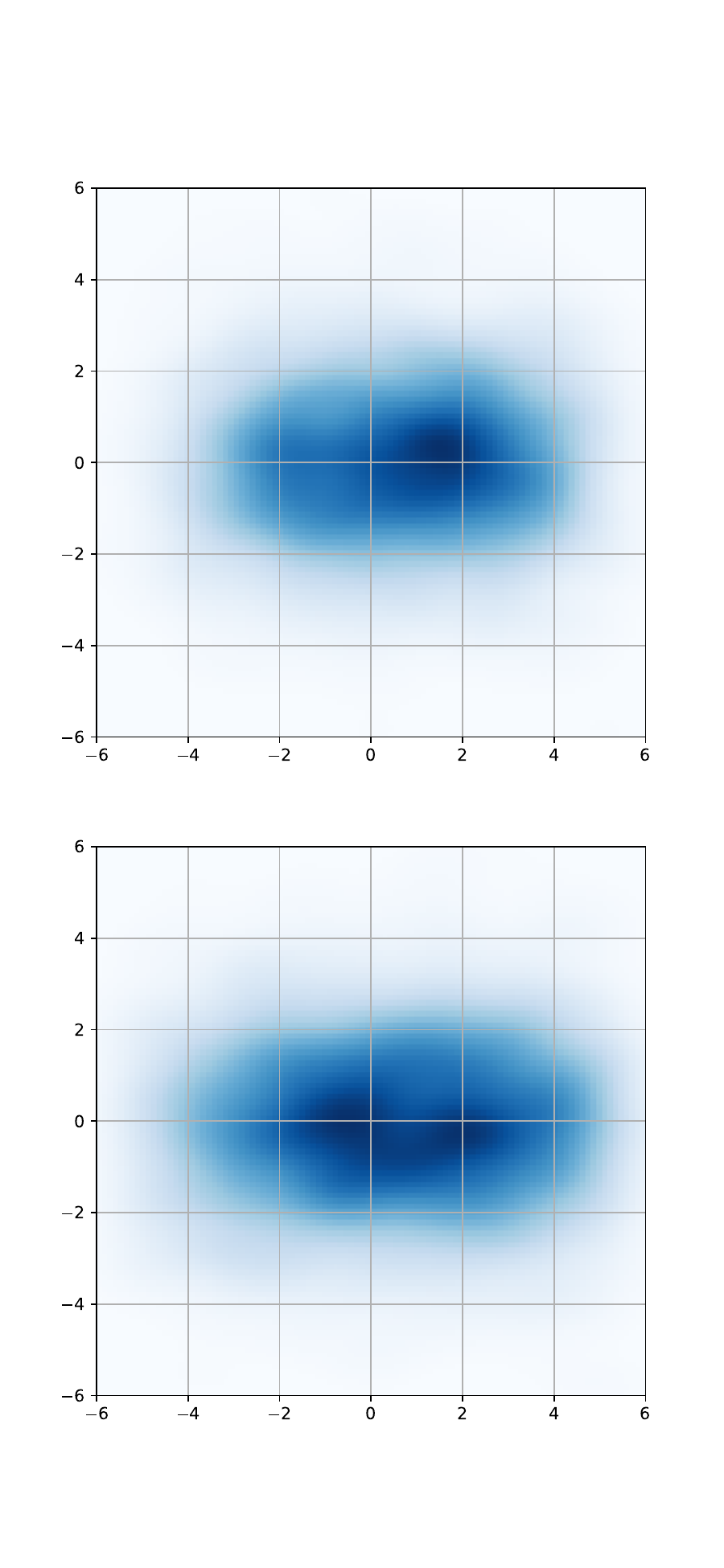}
    \centering
   \includegraphics[scale=0.19 ]{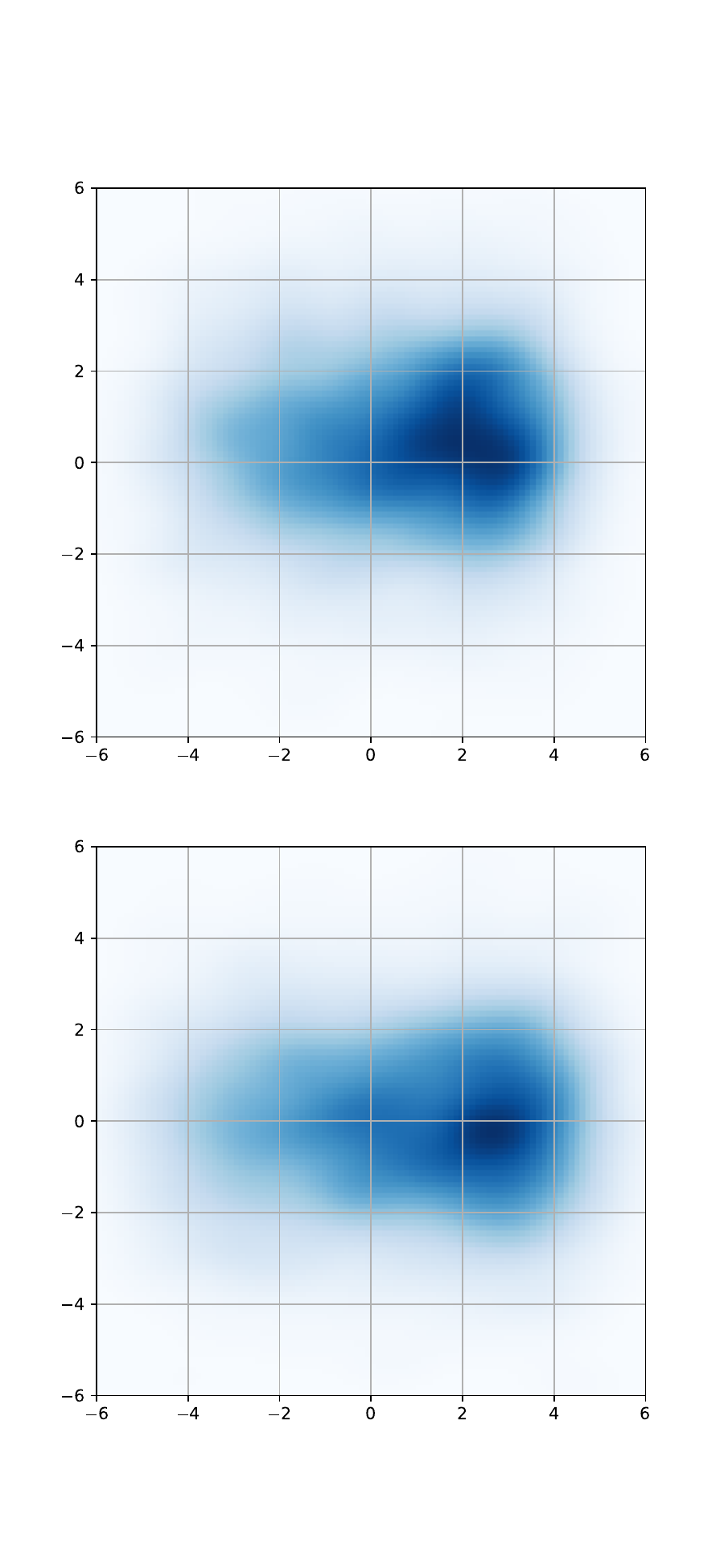}
   \caption{[Example \ref{example3}] Kernel density estimation (KDE) of distributions for $4000$ samples in the $xy$-plane. Top: Monte Carlo simulation \eqref{MC:example2dim}; Bottom: results from Algorithm \ref{algorithm}. From left to right: the distributions are shown at the $0$th, $40$th, $100$th, and $250$th time steps.}
\label{fig:example2dim-probability2d}
\end{figure}

Figure \ref{fig:example2dim-probability} illustrates the evolution of probability flows obtained from Equation \eqref{exa:doublewell}, computed using both the Monte Carlo simulation and the proposed approach. Figure \ref{fig:example2dim-probability2d} depicts the probability flow evolution on the $xy$-plane at different time steps for both methods. Finally, Figure \ref{fig:example2dim-TV} presents the TV distance between \( P^{\mathrm{MC}} \) and \( P^{\mathrm{S}} \), consistently demonstrating values on the order of \( 10^{-2} \). These results confirm the accuracy of the proposed method in approximating the probability flows for systems with interaction terms.

\begin{figure}[ht]
   \centering
   \includegraphics[width=0.55\textwidth]{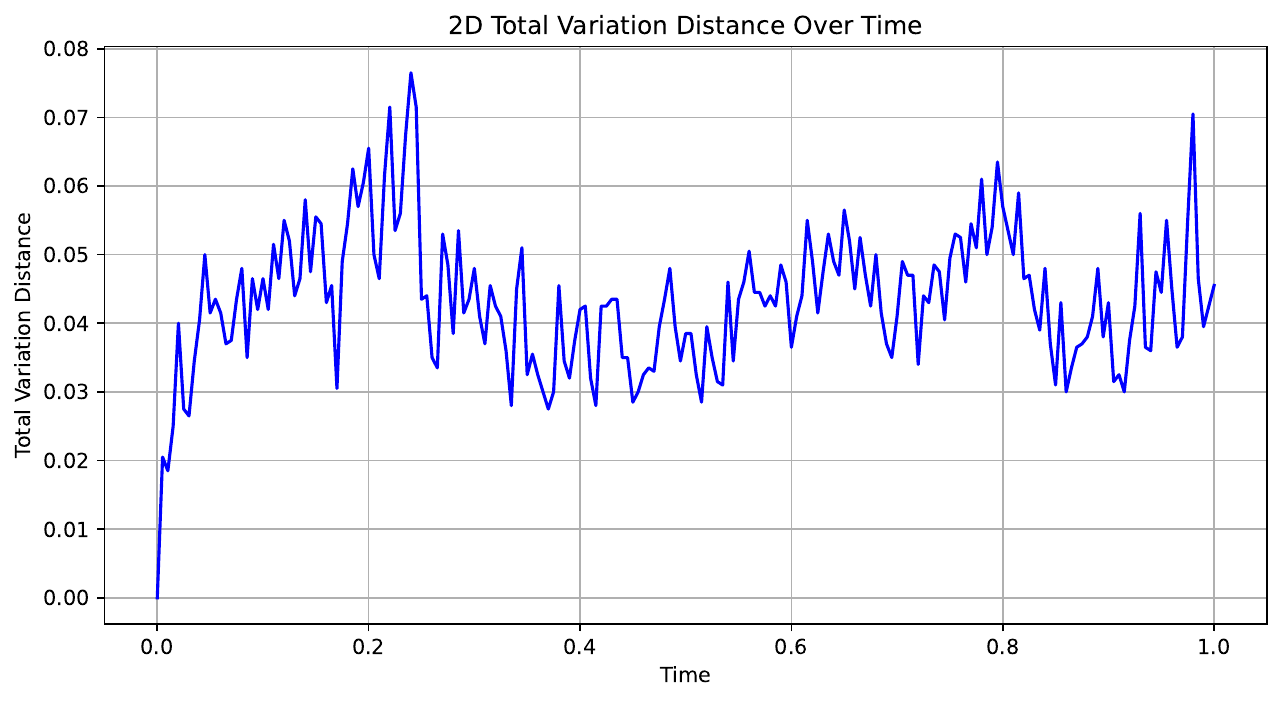}
    \caption{[Example \ref{example3}] The TV distance between $P^{\mathrm{MC}}$ and $P^{\mathrm{S}}$ for \eqref{exa:doublewell}.}
    \label{fig:example2dim-TV}
\end{figure}
\end{example}

\begin{example}[State-dependent jump-diffusion]\label{example4}
   To conclude our numerical study, we consider the following three-dimensional SDE which describes the asset price dynamics  under both jump and stochastic volatility conditions\cite{kristensen2024closed}
    {\small\begin{equation}\label{eqn:example3}
        \begin{aligned}
            \d s_t =&\ (\gamma - \delta - v_{t-}/2 - \Lambda(v_{t-},m_{t-})\bar{J} )\d t + \sqrt{|v_{t-}|} \d B_t^1 + \log(J^s_{t} + 1)\d N_t,\\
            \d v_t =&\  \kappa_v(m_{t-} - v_{t-})\d t + \sigma_v \sqrt{|v_{t-}|} (\d B^1_t + \sqrt{1-\rho^2}\d B_t^2 ) + J^v_t\d N_t,\\
            \d m_t =&\ \kappa_m (\alpha_m - m_{t-} )\d t + \sigma_m\sqrt{|m_{t-}|}\d B^3_t.
        \end{aligned}
    \end{equation}
}
Here \(B_t = (B_t^1, B_t^2, B_t^3)\) is the 
standard Brownian motion in \(\mathbb{R}^3\) and \(N_t\) is an independent Poisson process with the state-dependent rate \(\Lambda(v, m) = \Lambda_0 + \Lambda_1 v + \Lambda_2 m\). The jump component \(J_t^v\) follows the exponential distribution with the expectation \(\mu^v_J  \). Meanwhile, the jump component for the stock price, \(J_t^s\), follows a shifted log-normal distribution, where  \(J_t^s + 1\) is characterized by the parameters \(\mu_s\) and \(\sigma_s\), leading to \(\mathbb{E}(J^s_t)=\bar{J} = \exp(\mu_s + \sigma_s^2/2) - 1\).
The parameters are from  \cite{kristensen2024closed} and  specified as follows: $
     (\gamma, \delta, \kappa_v, \kappa_m, \alpha_m, \sigma_v, \sigma_m, \rho)$ = $(0.04, 0.015, 3.1206, 3.3168, 0.1125, 0.394, 0.0835, -0.688)$,  $(\Lambda_0, \Lambda_1, \Lambda_2) = (2.096, 21.225, 0)$, and  \( (\mu_s, \sigma_s, \mu_J^v) = (-0.012, 0.043, 0.002)\).



We set the initial condition as a three-dimensional Gaussian distribution with the mean vector $(5, 5, 5)^\top$ and the covariance matrix 
of a $3\times3$ identity matrix.
The numerical integration for $r$ is performed over a two-dimensional uniform grid defined on  $[-0.17, 0.17]\times[0, 0.0015]$, which is designed to capture the critical density features required for accurate numerical evaluation.
The Monte Carlo approximation is easy to achieve by noting that the numerical approximation of $\log(J^s_t+1)\d N_t$ is $\sum_{k=1}^{N_{\Delta t}}\xi^{J^s}_k$ where $\{\xi^{J^s}_k\}_{k\in\mathbb{N}}$ are i.i.d samples from the law of $\log(J_t^s+1)$ and \( N_{\Delta t} \sim \text{Po}(\Lambda \Delta t) \). Similarly $J_t^v\d N_t$ is approximated by $\sum_{k=1}^{N_{\Delta t}}\xi^{J^v}_k$ where $\{\xi^{J^v}_k\}_{k\in\mathbb{N}}$ are i.i.d samples from the law of $J_v^s$.

Let $X_t = (s_t, v_t, m_t)^\top$, the SDE \eqref{eqn:example3} can be rewritten as
{\footnotesize\begin{equation*}
    \begin{aligned}
        \d X_t =&\  \left[ \left(\begin{matrix}
            \gamma - \delta - \lambda_0\bar{J} \\
            0 \\
            \kappa_m \alpha_m \\
        \end{matrix}\right)  + \left(\begin{matrix}
            0 & -1/2 - \Lambda_1\bar{J} & -\Lambda_2\bar{J}\\
            0 &  -\kappa_v  &  \kappa_v\\
            0 & 0 & -\kappa_m
        \end{matrix}\right) X_t \right]\d t \\
        &\ + \left(\begin{matrix}
            \sqrt{|v_t|} & 0 & 0\\
            \sigma_v\sqrt{|v_t|} & \sigma_v\sqrt{|v_t|}\sqrt{1-\rho^2} & 0\\
            0 & 0 & \sigma_m\sqrt{|m_t|}
        \end{matrix}\right)\d B_t + \int_{\mathbb{R}^3}\left(\begin{matrix}
            \log (r_1 + 1)\\
            r_2\\
            0
        \end{matrix}\right)\mathcal{N}(\d t,\d r),
    \end{aligned}
\end{equation*}
}where $\mathcal{N}$ denotes a random Poisson measure on $\mathbb{R}_+\times\mathbb{R}^3$ with L\'{e}vy measure $\d t\times\Lambda(v_t,m_t)\nu_J(\d r)$, and the jump-size measure $\nu_J(\d r) = \nu_{J^s}(\d r_1)\times\nu_{J^v}(\d r_2)\times\nu_{J^m}(\d r_3)$. In this expression, $F(r,t)=G(r,t)=(\log(r_1+1),r_2,0)^\top$ .  

 \begin{figure}[ht]
    \centering
   \includegraphics[scale=0.19]{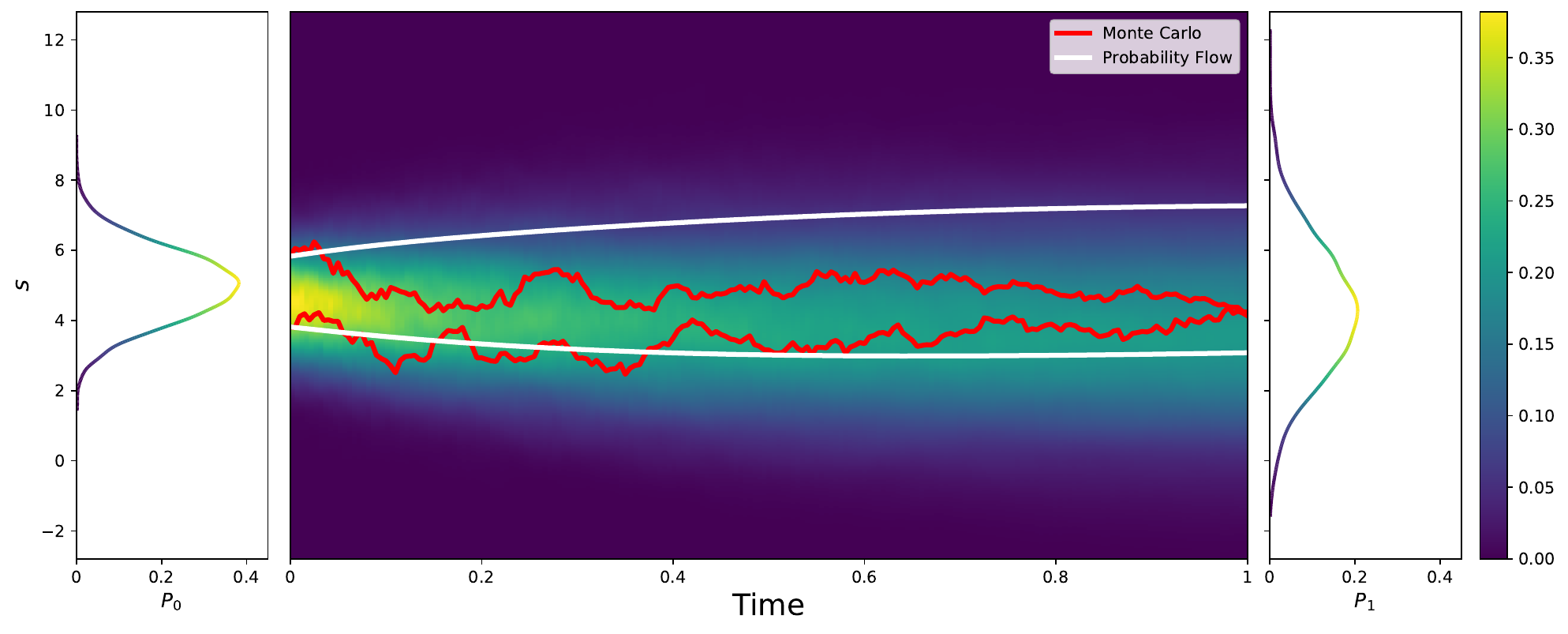}
   \centering
   \includegraphics[scale=0.19]{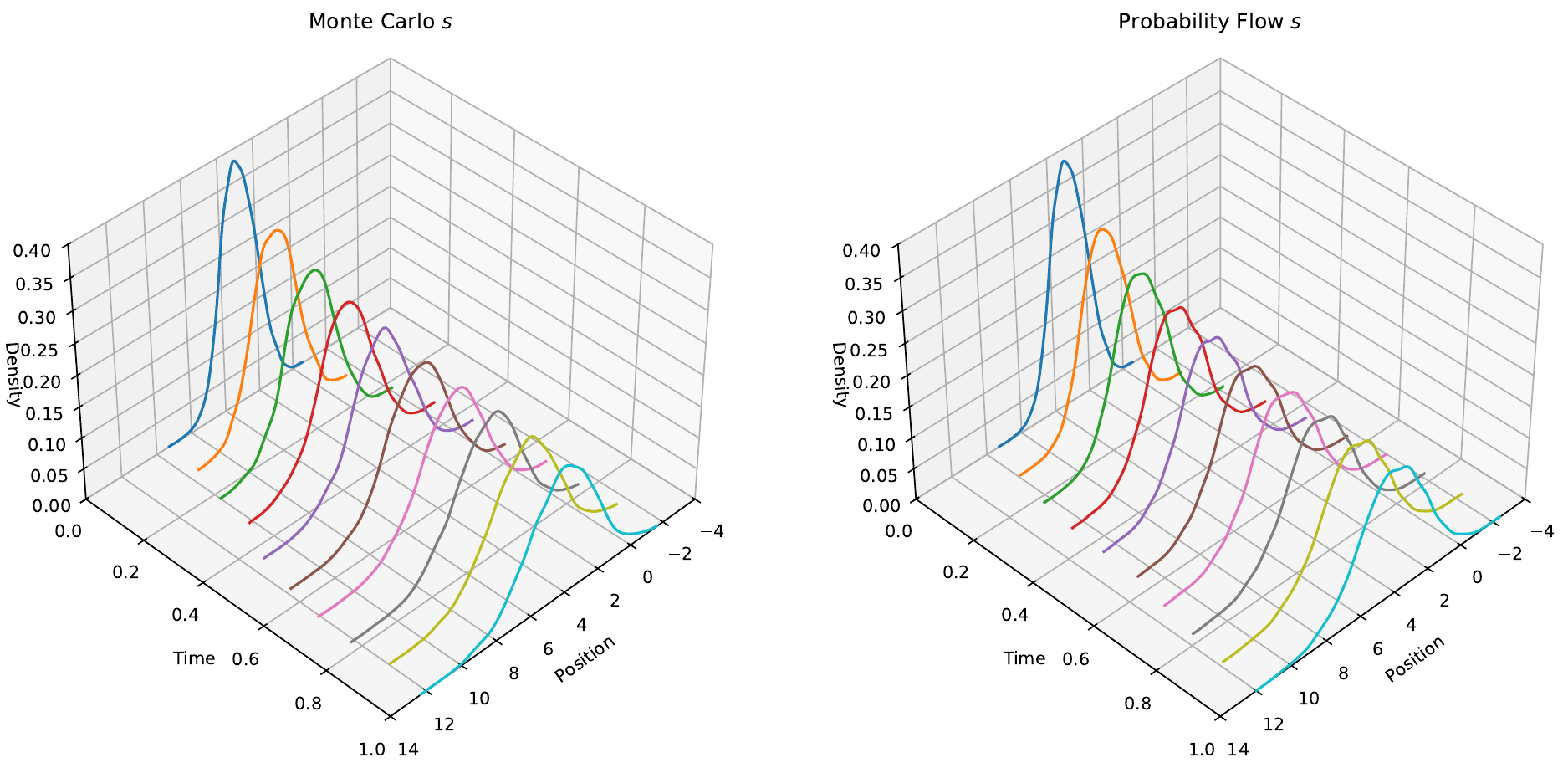}
   \centering
   \includegraphics[scale=0.19]{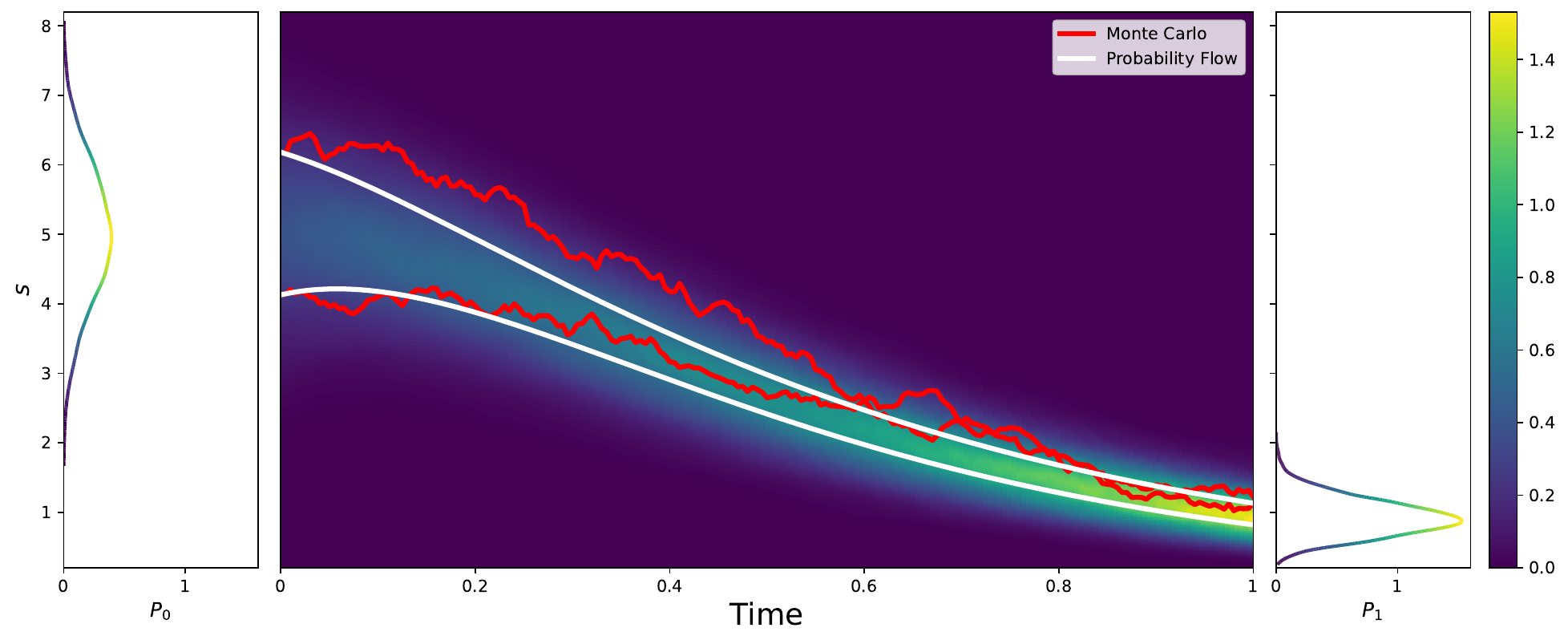}
   \centering
   \includegraphics[scale=0.19]{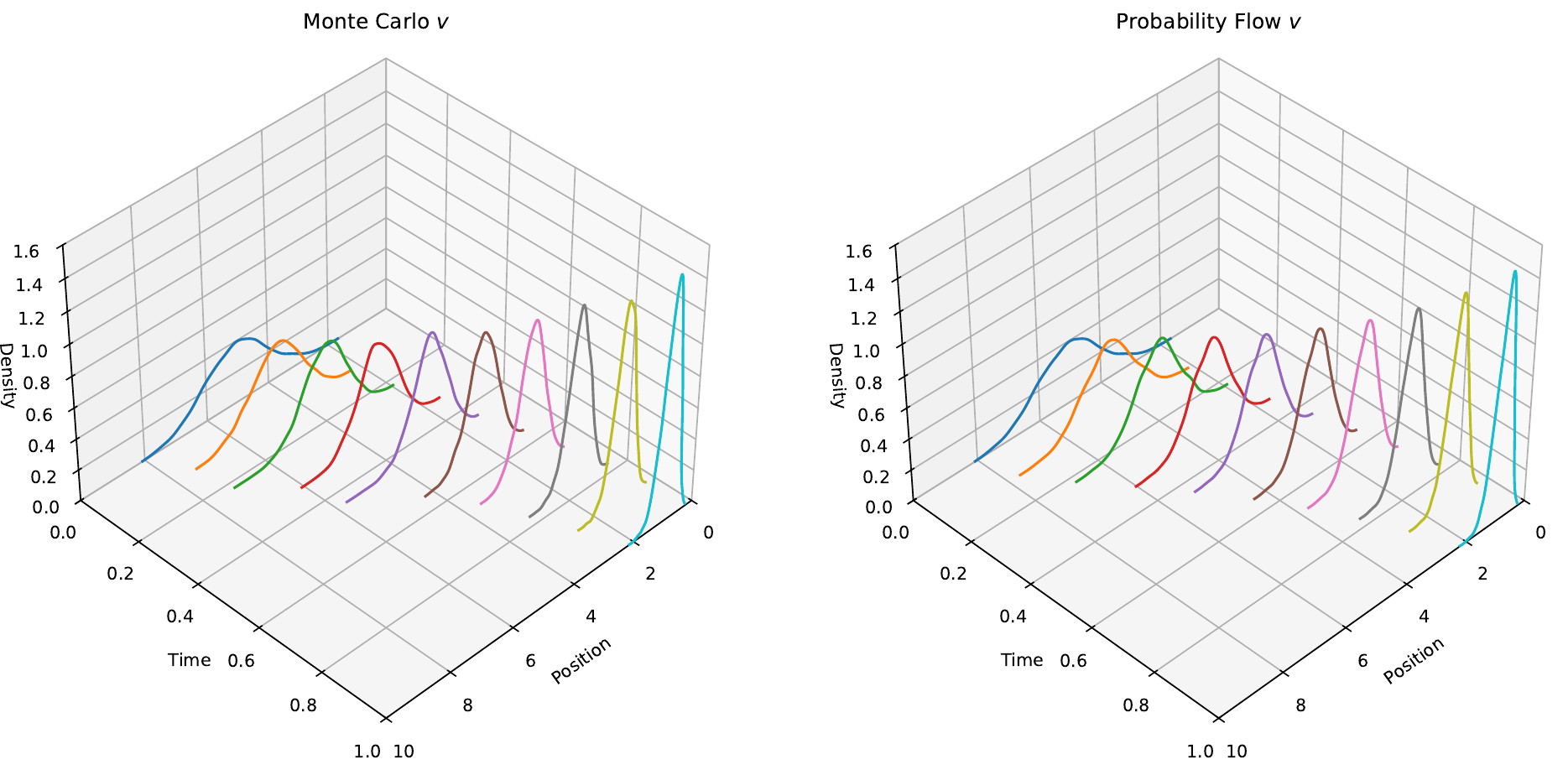}
   \centering
   \includegraphics[scale=0.19]{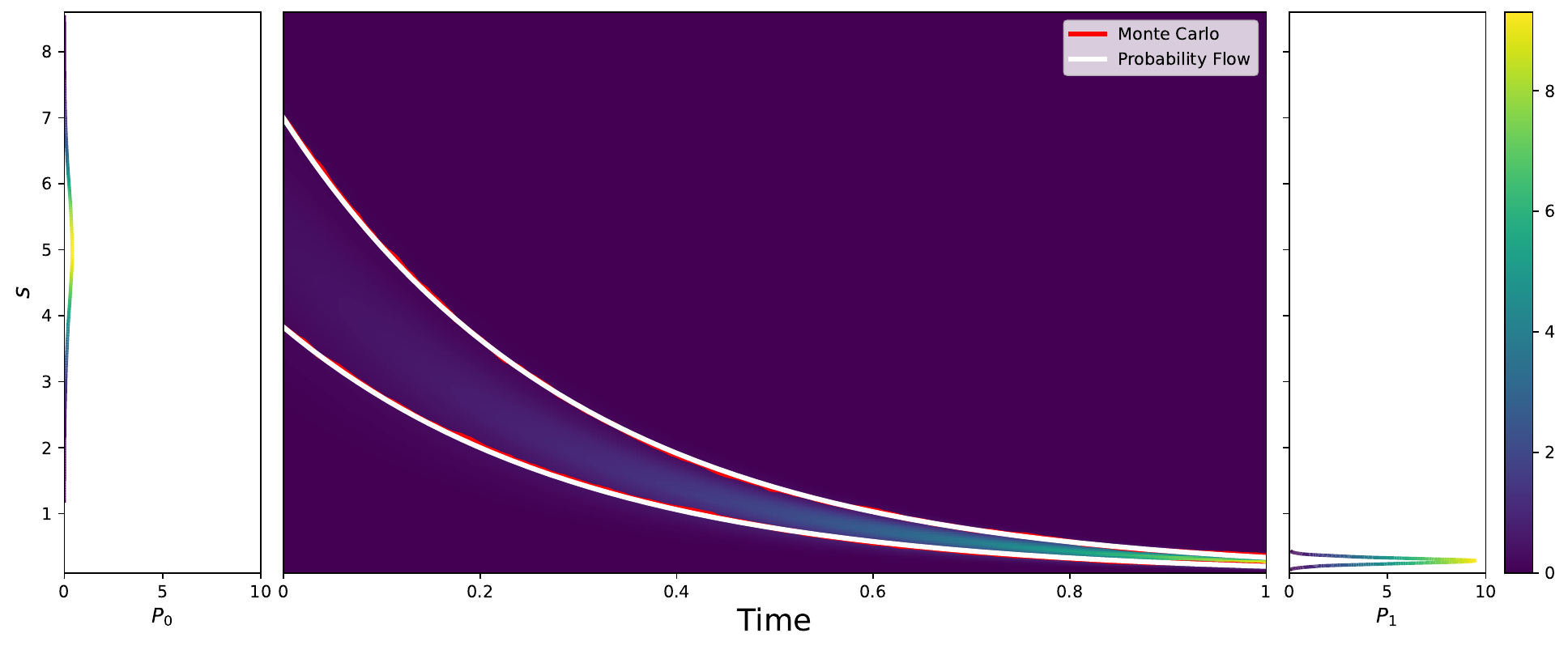}
   \centering
   \includegraphics[scale=0.19]{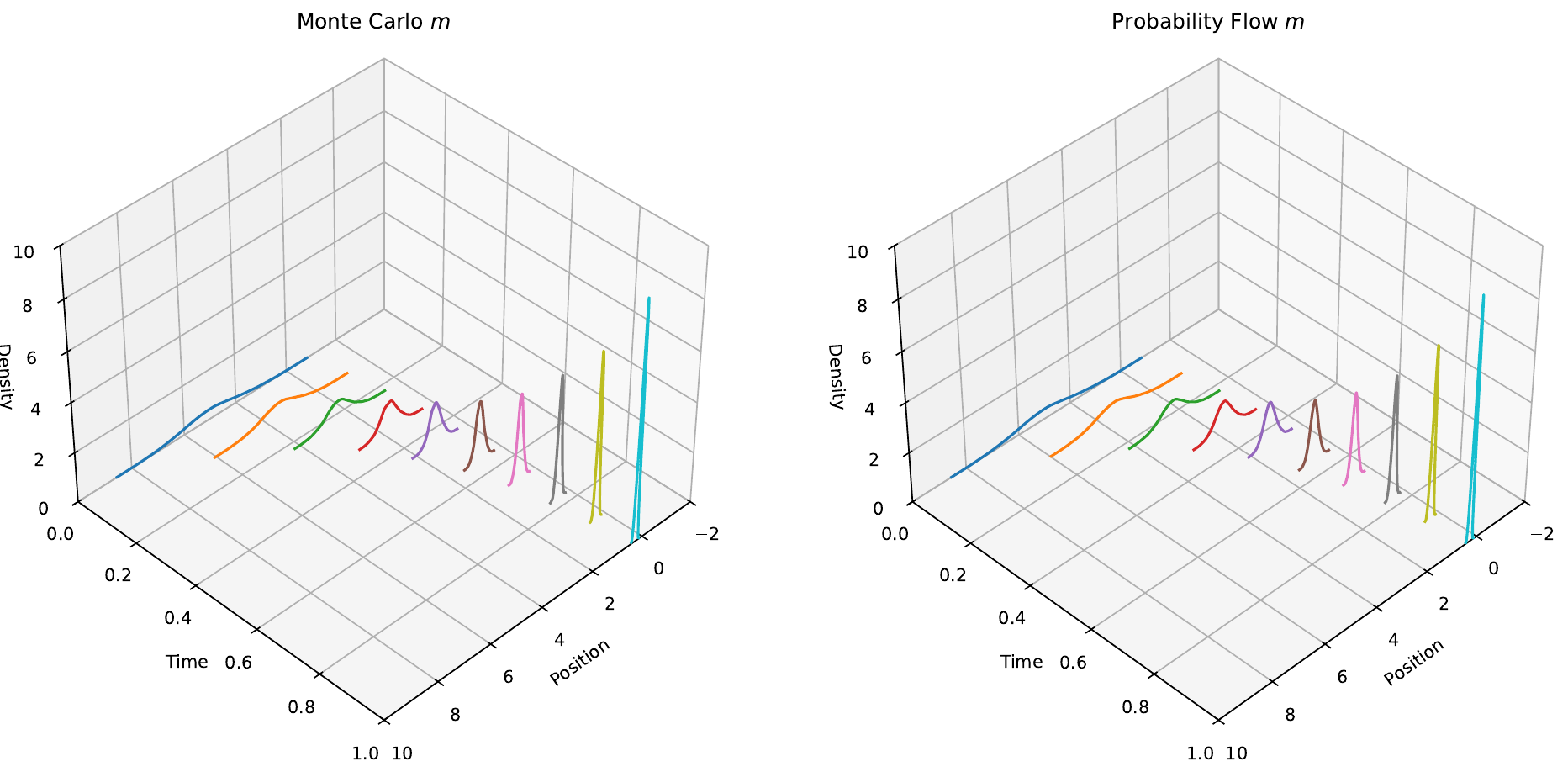}
    \caption{[Example \ref{example4}] Probability flows  along the axes: $s$ (top), $v$ (middle), and $m$ (bottom).}
    \label{fig:example3dim-probability}
\end{figure}


Figure \ref{fig:example3dim-probability} illustrates the evolution of probability flows for each marginal derived from Equation \eqref{eqn:example3}. Figure \ref{fig:example3dim-TV} presents the TV distance between the probabilities \(P^{\mathrm{MC}}\) and \(P^{\mathrm{S}}\), demonstrating that the distance remains on the order of \(10^{-2}\).  
These results highlight the capability of our proposed method to accurately and efficiently solve complex state-dependent jump-diffusion systems.

 \begin{figure}[ht]
    \centering
   \includegraphics[width=0.55\textwidth]{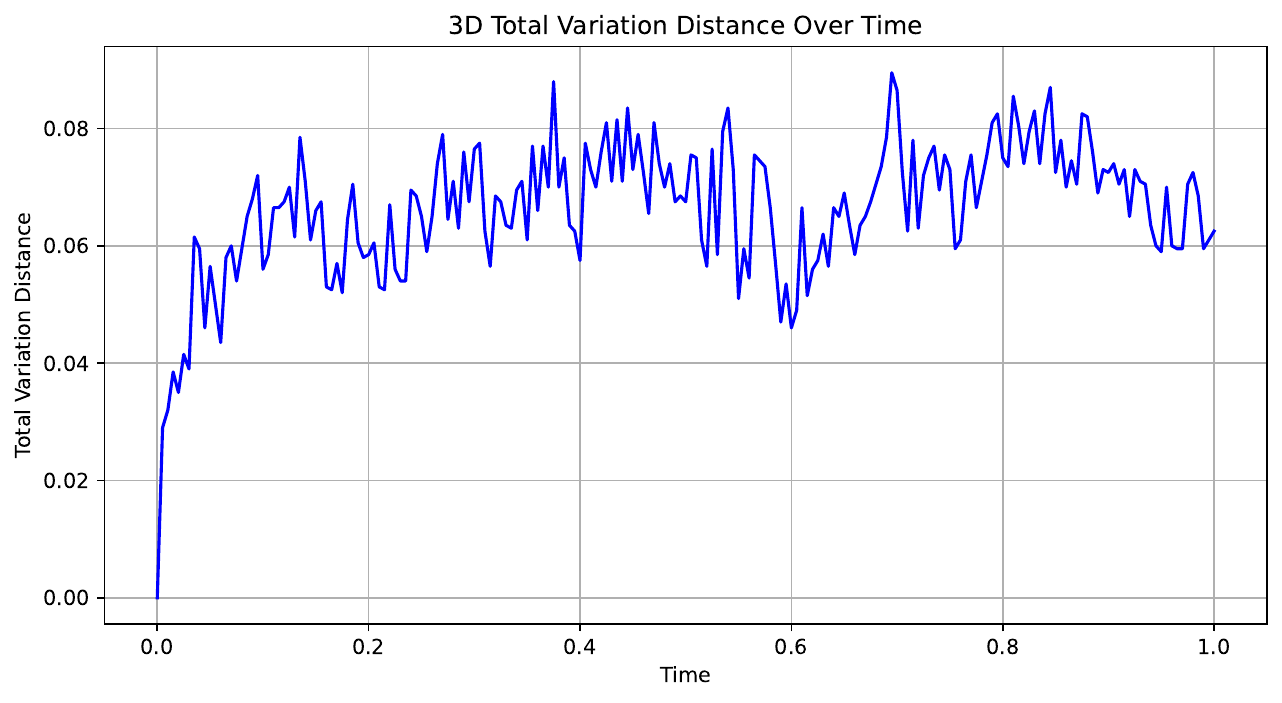}
   
    \caption{[Example \ref{example4}] The TV distance between $P^{\mathrm{MC}}$ and $P^{\mathrm{S}}$ for \eqref{eqn:example3}.}
    \label{fig:example3dim-TV}
\end{figure}

\end{example}

    \section{Conclusions}
    \label{sec:conclusions} In this paper, we introduced the concept of the
    \emph{generalized score} function, which unifies the conventional diffusion
    score and the L\'{e}vy score, thereby accommodating a broad class of non-Gaussian
    processes. Leveraging this generalized score function, we developed a score-based
    particle algorithm to efficiently solve the nonlinear L\'{e}vy-Fokker-Planck
    equation through sampling of probability flow data. Furthermore, we demonstrate
    that the Kullback-Leibler divergence between the numerical solution and the ground
    truth is bounded by a specific loss function, derived from a fixed-point
    perspective. The rigorous numerical analysis was conducted to validate the theoretical
    properties of our algorithm. Finally, several numerical examples were presented
    to illustrate the effectiveness and accuracy of the proposed method.

    \appendix

    \section*{Data Availability}
    The datasets generated and/or analyzed during the current study are
    available in the GitHub repositorry, accessible via the following link:
    \href{https://github.com/cliu687/sbtm-levy}{https://github.com/cliu687/sbtm-levy}.
    All scripts and supporting files necessary to reproduce the findings are openly
    provided in this repositorry.

    \bibliographystyle{siamplain}
    \bibliography{references}

@book{applebaum2009levy,
  title={{L}{\'e}vy processes and stochastic calculus},
  author={Applebaum, David},
  year={2009},
  publisher={Cambridge university press}
}

@book{tome2014book,
  title={Stochastic Dynamics and Irreversibility},
  author={Tom{\'e}, T. and Oliveira, M.J.},
  isbn={9783319117690},
  lccn={2014955664},
  series={Graduate Texts in Physics},
  year={2014},
  publisher={Springer International Publishing}
}

@article{boffi2023probability,
  title={Probability flow solution of the {F}okker--{P}lanck equation},
  author={Boffi, Nicholas M and Vanden-Eijnden, Eric},
  journal={Machine Learning: Science and Technology},
  volume={4},
  number={3},
  pages={035012},
  year={2023},
  publisher={IOP Publishing}
}

@article{maoutsa2020interacting,
  title={Interacting particle solutions of {F}okker--{P}lanck equations through gradient--log--density estimation},
  author={Maoutsa, Dimitra and Reich, Sebastian and Opper, Manfred},
  journal={Entropy},
  volume={22},
  number={8},
  pages={802},
  year={2020},
  publisher={MDPI}
}

@article{santambrogio2015optimal,
  title={Optimal transport for applied mathematicians},
  author={Santambrogio, Filippo},
  journal={Birk{\"a}user, NY},
  volume={55},
  number={58-63},
  year={2015},
  publisher={Springer}
}

@book{villani2009optimal,
  title={Optimal transport: old and new},
  author={Villani, C{\'e}dric and others},
  volume={338},
  year={2009},
  publisher={Springer}
}

@article{lu2024score,
  title={Score-based transport modeling for mean-field {F}okker-{P}lanck equations},
  author={Lu, Jianfeng and Wu, Yue and Xiang, Yang},
  journal={Journal of Computational Physics},
  volume={503},
  pages={112859},
  year={2024},
  publisher={Elsevier}
}

@inproceedings{shen2022self,
  title={Self-consistency of the {F}okker--{P}lanck equation},
  author={Shen, Zebang and Wang, Zhenfu and Kale, Satyen and Ribeiro, Alejandro and Karbasi, Amin and Hassani, Hamed},
  booktitle={Conference on Learning Theory},
  pages={817--841},
  year={2022},
  organization={PMLR}
}

@article{shen2024entropy,
  title={Entropy-dissipation informed neural network for {M}ckean-{V}lasov type {PDE}s},
  author={Shen, Zebang and Wang, Zhenfu},
  journal={Advances in Neural Information Processing Systems},
  volume={36},
  year={2024}
}

@inproceedings{li2023selfconsistent,
title={Self-Consistent Velocity Matching of Probability Flows},
author={Lingxiao Li and Samuel Hurault and Justin Solomon},
booktitle={Thirty-seventh Conference on Neural Information Processing Systems},
year={2023}
}

@article{kanazawa2020loopy,
  title={Loopy {L}{\'e}vy flights enhance tracer diffusion in active suspensions},
  author={Kanazawa, Kiyoshi and Sano, Tomohiko G and Cairoli, Andrea and Baule, Adrian},
  journal={Nature},
  volume={579},
  number={7799},
  pages={364--367},
  year={2020},
  publisher={Nature Publishing Group UK London}
}

@article{liu2023stochastic,
  title={Stochastic {M}cKean--{V}lasov equations with {L}{\'e}vy noise: {E}xistence, attractiveness and stability},
  author={Liu, Huoxia and Lin, Judy Yangjun},
  journal={Chaos, Solitons \& Fractals},
  volume={177},
  pages={114214},
  year={2023},
  publisher={Elsevier}
}

@article{olivera2024microscopic,
  title={Microscopic derivation of non-local models with anomalous diffusions from stochastic particle systems},
  author={Olivera, Christian and Simon, Marielle},
  journal={arXiv preprint arXiv:2404.03772},
  year={2024}
}

@article{de2024multidimensional,
  title={Multidimensional stable driven {M}cKean--{V}lasov {SDE}s with distributional interaction kernel: a regularization by noise perspective},
  author={de Raynal, P-E Chaudru and Jabir, J-F and Menozzi, St{\'e}phane},
  journal={Stochastics and Partial Differential Equations: Analysis and Computations},
  pages={1--54},
  year={2024},
  publisher={Springer}
}

@article{barthelemy2008levy,
  title={A {L}{\'e}vy flight for light},
  author={Barthelemy, Pierre and Bertolotti, Jacopo and Wiersma, Diederik S},
  journal={Nature},
  volume={453},
  number={7194},
  pages={495--498},
  year={2008},
  publisher={Nature Publishing Group UK London}
}

@article{song2018neuronal,
  title={Neuronal messenger ribonucleoprotein transport follows an aging {L}{\'e}vy walk},
  author={Song, Minho S and Moon, Hyungseok C and Jeon, Jae-Hyung and Park, Hye Yoon},
  journal={Nature communications},
  volume={9},
  number={1},
  pages={1--8},
  year={2018},
  publisher={Nature Publishing Group}
}

@article{liang2021exponential,
  title={Exponential ergodicity for {SDE}s and {M}cKean--{V}lasov processes with {L}{\'e}vy noise},
  author={Liang, Mingjie and Majka, Mateusz B and Wang, Jian},
  journal={Annales de l'Institut Henri Poincare (B) Probabilites et statistiques},
  volume={57},
  number={3},
  pages={1665--1701},
  year={2021}
}

@article{Ilin2024TransportBP,
  title={Transport based particle methods for the {F}okker-{P}lanck-{L}andau equation},
  author={Vasily Ilin and Jingwei Hu and Zhenfu Wang},
  journal={ArXiv},
  year={2024},
  volume={abs/2405.10392}
}

@article{huang2024score,
  title={A score-based particle method for homogeneous {L}andau equation},
  author={Huang, Yan and Wang, Li},
  journal={arXiv preprint arXiv:2405.05187},
  year={2024}
}

@article{jabin2018quantitative,
  title={Quantitative estimates of propagation of chaos for stochastic systems with ${W^{-1, \infty}}$ kernels},
  author={Jabin, Pierre-Emmanuel and Wang, Zhenfu},
  journal={Inventiones mathematicae},
  volume={214},
  pages={523--591},
  year={2018},
  publisher={Springer}
}

@article{golse2016dynamics,
  title={On the dynamics of large particle systems in the mean field limit},
  author={Golse, Fran{\c{c}}ois},
  journal={Macroscopic and large scale phenomena: coarse graining, mean field limits and ergodicity},
  pages={1--144},
  year={2016},
  publisher={Springer}
}

@article{chevallier2017estimation,
  title={Estimation of {L}{\'e}vy-driven {O}rnstein--{U}hlenbeck processes: Application to modeling of {CO}2 and fuel-switching},
  author={Chevallier, Julien and Goutte, St{\'e}phane},
  journal={Annals of Operations Research},
  volume={255},
  number={1},
  pages={169--197},
  year={2017},
  publisher={Springer}
}

@article{kristensen2024closed,
  title={Closed-form approximations of moments and densities of continuous--time {M}arkov models},
  author={Kristensen, Dennis and Lee, Young Jun and Mele, Antonio},
  journal={Journal of Economic Dynamics and Control},
  pages={104948},
  year={2024},
  publisher={Elsevier}
}

@article{ariga2021noise,
  title={Noise-induced acceleration of single molecule kinesin-1},
  author={Ariga, Takayuki and Tateishi, Keito and Tomishige, Michio and Mizuno, Daisuke},
  journal={Physical review letters},
  volume={127},
  number={17},
  pages={178101},
  year={2021},
  publisher={APS}
}

@article{huang2025probability,
  title={Probability flow approach to the Onsager--Machlup functional for jump-diffusion processes},
  author={Huang, Yuanfei and Zhou, Xiang and Duan, Jinqiao},
  journal={SIAM Journal on Applied Mathematics},
  volume={85},
  number={2},
  pages={524--547},
  year={2025},
  publisher={SIAM}
}

@article{boffi2024deep,
  title={Deep learning probability flows and entropy production rates in active matter},
  author={Boffi, Nicholas M and Vanden-Eijnden, Eric},
  journal={P NATL ACAD SCI USA},
  volume={121},
  number={25},
  pages={e2318106121},
  year={2024},
  publisher={National Acad Sciences}
}

@article{hyvarinen2005estimation,
  title={Estimation of non-normalized statistical models by score matching.},
  author={Hyv{\"a}rinen, Aapo and Dayan, Peter},
  journal={Journal of Machine Learning Research},
  volume={6},
  number={4},
  year={2005}
}

@inproceedings{NEURIPS2021_940392f5,
 author = {De Bortoli, Valentin and Thornton, James and Heng, Jeremy and Doucet, Arnaud},
 booktitle = {Advances in Neural Information Processing Systems},
 pages = {17695--17709},
 publisher = {Curran Associates, Inc.},
 title = {Diffusion {S}chr\"{o}dinger Bridge with Applications to Score-Based Generative Modeling},
 volume = {34},
 year = {2021}
}

@InProceedings{chen_improved_2022,
  title = 	 {Improved Analysis of Score-based Generative Modeling: User-Friendly Bounds under Minimal Smoothness Assumptions},
  author =       {Chen, Hongrui and Lee, Holden and Lu, Jianfeng},
  booktitle = 	 {Proceedings of the 40th International Conference on Machine Learning},
  pages = 	 {4735--4763},
  year = 	 {2023},
  volume = 	 {202},
}

@inproceedings{ho_2020_denoising,
    author = {Ho, Jonathan and Jain, Ajay and Abbeel, Pieter},
    booktitle = {Advances in Neural Information Processing Systems},
    pages = {6840--6851},
    publisher = {Curran Associates, Inc.},
    title = {Denoising Diffusion Probabilistic Models},
    volume = {33},
    year = {2020}
}

@inproceedings{
song2021denoising,
title={Denoising Diffusion Implicit Models},
author={Jiaming Song and Chenlin Meng and Stefano Ermon},
booktitle={International Conference on Learning Representations},
year={2021},
}

@article{gottwald2024stable1,
  title={Stable generative modeling using diffusion maps},
  author={Gottwald, Georg and Li, Fengyi and Marzouk, Youssef and Reich, Sebastian},
  journal={arXiv preprint arXiv:2401.04372},
  year={2024}
}

@article{liu2024diffusion,
  title={Diffusion-Model-Assisted Supervised Learning of Generative Models for Density Estimation},
  author={Liu, Yanfang and Yang, Minglei and Zhang, Zezhong and Bao, Feng and Cao, Yanzhao and Zhang, Guannan},
  journal={Journal of Machine Learning for Modeling and Computing},
  volume={5},
  number={1},
  year={2024},
  publisher={Begel House Inc.}
}

@article{gallon2024overview,
  title={An overview of diffusion models for generative artificial intelligence},
  author={Gallon, Davide and Jentzen, Arnulf and von Wurstemberger, Philippe},
  journal={arXiv preprint arXiv:2412.01371},
  year={2024}
}

@inproceedings{song_generative_2019,
    title = {Generative Modeling by Estimating Gradients of the Data Distribution},
    volume = {32},
    booktitle = {NeurIPS},
    publisher = {Curran Associates, Inc.},
    author = {Song, Yang and Ermon, Stefano},
    year = {2019},
}

@inproceedings{song2020sliced,
  title={Sliced score matching: A scalable approach to density and score estimation},
  author={Song, Yang and Garg, Sahaj and Shi, Jiaxin and Ermon, Stefano},
  booktitle={Uncertainty in Artificial Intelligence},
  pages={574--584},
  year={2020},
  organization={PMLR}
}

@article{song2020score,
  title={Score-based generative modeling through stochastic differential equations},
  author={Song, Yang and Sohl-Dickstein, Jascha and Kingma, Diederik P and Kumar, Abhishek and Ermon, Stefano and Poole, Ben},
  journal={International Conference on Learning Representations (ICLR)},
  year={2021}
}

@article{cao2024exploring,
  title={Exploring the optimal choice for generative processes in diffusion models: Ordinary vs stochastic differential equations},
  author={Cao, Yu and Chen, Jingrun and Luo, Yixin and Zhou, Xiang},
  journal={Advances in Neural Information Processing Systems},
  volume={36},
  year={2024}
}

@article{yoon2023score,
  title={Score-based generative models with {L}{\'e}vy processes},
  author={Yoon, Eun BI and Park, Keehun and Kim, Sungwoong and Lim, Sungbin},
  journal={Advances in Neural Information Processing Systems},
  volume={36},
  pages={40694--40707},
  year={2023}
}

@article{hu2024score,
  title={Score-f{PINN}: Fractional Score-Based Physics-Informed Neural Networks for High-Dimensional {F}okker-{P}lanck-{L}\'{e}vy Equations},
  author={Hu, Zheyuan and Zhang, Zhongqiang and Karniadakis, George Em and Kawaguchi, Kenji},
  journal={arXiv preprint arXiv:2406.11676},
  year={2024}
}

@article{caffarelli2007extension,
  title={An extension problem related to the fractional {L}aplacian},
  author={Caffarelli, Luis and Silvestre, Luis},
  journal={Communications in partial differential equations},
  volume={32},
  number={8},
  pages={1245--1260},
  year={2007},
  publisher={Taylor \& Francis}
}

@article{huang2025entropy,
  title={Entropy production in non-{G}aussian active matter: A unified fluctuation theorem and deep learning framework},
  author={Huang, Yuanfei and Liu, Chengyu and Miao, Bing and Zhou, Xiang},
  journal={arXiv preprint arXiv:2504.06628},
  year={2025}
}
\end{document}